\documentclass[10pt,amssymb,amsfonts,amscd]{amsart}
\usepackage{amscd}
\usepackage{amssymb}
\usepackage{amsfonts}
\newcounter{THNO}[section]
\newcounter{SNO}[THNO]
\renewcommand{\theTHNO}{\arabic{section}.\arabic{THNO}}

\renewcommand{\theSNO}{\arabic{section}.\arabic{THNO}.\arabic{SNO}}

\newcounter{multieq}[equation]

\newcounter{tmp}

\def\sec#1{\refstepcounter{section}\par\vspace{0.9cm}
\par\noindent
{\large\bf\arabic{section}. \ #1    }
\nopagebreak\par\vspace{0.3cm}}
\def\sect#1{\par\vspace{0.5cm}
\par\noindent
{\large\bf \ #1    }
\nopagebreak\par\vspace{0.3cm}}
\def\lm#1{\refstepcounter{SNO}\par
\medskip
\par\noindent\begingroup \it
\leftskip=0em\hspace{0em}{\bf \theSNO\, #1.}}
\def\elm{\par\endgroup}

\def\sssec{\refstepcounter{SNO}\par
\medskip
\par\noindent
\leftskip=0em\hspace{0em}{\bf \theSNO. \,}}
\def\th#1{\refstepcounter{THNO}\par
\medskip
\par\noindent\begingroup \it
\leftskip=0em\hspace{0em}{\bf \theTHNO\, #1.}}
\def\eth{\par\endgroup}
\def\ssec{\refstepcounter{THNO}\par
\medskip
\par\noindent
\leftskip=0em\hspace{0em}{\bf \theTHNO. \,}}
\def\pr{\smallskip\par\noindent{\sc Proof. }}

\def\ex#1{\refstepcounter{SNO}\par\medskip
\par\noindent\begingroup \rm
\leftskip=0em\hspace{0em}{\bf \theSNO\, #1.}}
\def\eex{\par\endgroup}

\def\db#1{ D^b({#1})}

\def\h#1,#2{{\rm Hom}({#1}\:,\; {#2})}
\def\H#1,#2,#3,#4{{\rm Hom}^{#1}_{#2}({#3}\:,\; {#4})}
\def\E#1,#2,#3,#4{{\rm Ext}^{#1}_{#2}({#3}\:,\; {#4})}
\def\lto{\longrightarrow}

\def\da{\big\downarrow}

\def\ss#1{\scriptsize{$#1$}}
\def\bl{\boldsymbol}
\def\ot{\otimes}
\def\bt{\boxtimes}
\def\ts{\times}
\def\op{\oplus}

\def\wh#1{\widehat{#1}}
\def\wt#1{\widetilde{#1}}
\def\R{{\Bbb{R}}}
\def\C{{\Bbb{C}}}
\def\Q{{\Bbb{Q}}}
\def\Z{{\Bbb{Z}}}

\def\Qu{{\Bbb{H}}}

\def\S{{\bl{S}}}

\def\O{{\mathcal{O}}}
\def\E{{\mathcal{E}}}
\def\F{{\mathcal{F}}}
\def\G{{\mathcal{G}}}
\def\H{{\mathcal{H}}}
\def\D{{\mathcal{D}}}
\def\J{{\mathcal{J}}}
\def\L{{\mathcal{L}}}
\def\cP{{\mathcal{R}}}
\def\La{{\varLambda}}
\def\Ga{{\varGamma}}
\def\Si{{\varSigma}}
\def\te{{\vartheta}}
\def\End{{\rm E}{\rm n}{\rm d}}
\def\Hom{{\rm H}{\rm o}{\rm m}}
\def\GL{{\rm G}{\rm L}}
\def\SO{{\rm S}{\rm O}}
\def\SU{{\rm S}{\rm U}}
\def\U{{\rm{U}}}
\def\SL{{\rm S}{\rm L}}
\def\Pic{{\rm P}{\rm i}{\rm c}}
\def\Sp{{\rm S}{\rm p}}
\title{\bf  Mirror symmetry for abelian varieties }
\author{Vasily Golyshev }
\address{\hspace*{-0.5cm} Institute for Problems of Information Transmission,
Russian Academy of Sciences,
19 Bolshoi Karetnyi, Moscow 101447, Russia}
\email{golyshev@mccme.ru}

\author{ Valery Lunts}
\thanks{the second author
was supported by the NSF}
\address{\hspace*{-0.5cm} Department of Mathematics, Indiana University,
Bloomington, IN 47405, USA}
\email{vlunts@indiana.edu}

\author{Dmitri Orlov}
\thanks{the third author was partly supported by the RFFI--99--01--01144}
\address{\hspace*{-0.5cm} Algebra Section, Steklov Mathematical Institute,
Russian Academy of Sciences,
8 Gubkin str., GSP-1,  Moscow 117966, Russia}
\email{orlov@mi.ras.ru}

\date{}
\begin{document}
\maketitle

\sec{Introduction.}

\medskip
\noindent{\bf 0.1.} We define the relation of mirror symmetry on the
class of pairs (complex abelian variety $A$ + an element
of the
complexified ample cone of $A$) and study its properties. More
precisely, let $A$ be a complex abelian variety, $C_A^a\subset NS_A(\R)$
-- the ample cone of $A$ and put
$$C_A^{\pm}:=NS_A(\R)\pm iC_A^a,$$
$$C_A:=C^+_A\sqcup C_A^-\subset NS_A(\C).$$
If $\omega _A\in C_A$ we call $(A,\omega _A)$ an {\it algebraic pair}.
In this work we define the notion of mirror symmetry for algebraic
pairs . The definition is given below in 0.4.1. Our notion of mirror
symmetry is defined purely in the language of algebraic geometry and algebraic
groups. We explain in 0.9 below that our construction is compatible
with the pictures of mirror symmetry of Kontsevich \cite{Ko} (see also \cite{PZ}) and
Strominger, Yau, Zaslow \cite{SYZ} (see also \cite{Mo}, \cite{Gro1}, \cite{Gro2}).

\medskip
\noindent{\bf 0.2.} Let us discuss some properties of this notion.

\smallskip
\noindent{1.} Not every algebraic pair has a mirror symmetric one (9.5.1),
but for ``general'' abelian varieties $A$ and any $\omega _A\in C_A$ the
pair $(A,\omega _A)$ does have a symmetric one (9.6.3). Also for any
abelian variety $A$ there exists $\omega _A\in C_A$ such that the
algebraic pair $(A,\omega _A)$ has a mirror symmetric one
(9.6.1).

\smallskip

\noindent{2.} Suppose two algebraic pairs $(B, \omega_{B})$ and $(C, \omega_{C})$ are both mirror symmetric to
the same algebraic pair $(A, \omega _A)$ then the derived
categories of coherent sheaves on $B$ and $C$ are equivalent (as
triangulated categories) (9.2.6).
Thus in particular for any algebraic pair
$(A,\omega _A)$ there are only finitely many isomorphism classes of
abelian varieties which are mirror symmetric to $(A,\omega _A)$ (9.2.3).

Conversely, if some algebraic pairs $(B, \omega_B)$ and $(A, \omega_A)$ are mirror symmetric
and an abelian variety $C$ is such that the derived categories of coherent sheaves on $C$ and $B$ are equivalent,
then there is $\omega_{C}$ such that the pairs $(C, \omega_C)$ and $(A, \omega_A)$ are mirror symmetric too.

\smallskip
\noindent{3.} Suppose that an algebraic pair $(B, \omega_B)$ is mirror symmetric
to the algebraic pair $(A,\omega _A).$ Starting with $(A,\omega _A)$
we cannot construct $B,$ but we can construct the product $B\times \wh{B},$
where $\wh{B}$ is the dual abelian variety.

\smallskip
Given an abelian variety $A$ consider its total cohomology $H^*(A,\Q).$
This space has a ``horizontal'' and ``vertical'' structures which are
discussed below in 0.3. More precisely, $H^*(A,\Q)$ is acted upon by
two reductive algebraic $\Q$--groups -- the ``horizontal'' and the
``vertical'', -- and these groups commute. The ``horizontal''
group is the Hodge group or the special Mumford-Tate group (its action
preserves each cohomology space $H^k(A,\Q)$). The vertical group
is defined below in 0.3.

\smallskip
\noindent{4.} Suppose that a mirror symmetry between algebraic
pairs $(A,\omega _A)$ and $(B,\omega _B)$ is given. Then there exists a
natural isomorphism
$$\beta :H^*(A,\Z)\stackrel{\sim}{\lto}H^*(B,\Z)$$
such that $\beta _{\Q}$ ``interchanges'' the horizontal and vertical
structures on  $H^*(A,\Q)$ and $H^*(B,\Q)$ respectively (9.3.3). There
exists a canonical choice (up to $\pm $) of such a $\beta .$ This isomorphism
will either preserve or switch the parity of the cohomology groups
depending on the parity of the dimension of $A$ and $B.$ For example,
if $A$ and $B$ are elliptic curves then $\beta $ induces isomorphisms
$$\beta :H^1(A,\Z)\stackrel{\sim}{\lto} H^0(B,\Z)\oplus H^2(B,\Z),$$
$$\beta :H^0(A,\Z)\oplus H^2(A,\Z)\stackrel{\sim}{\lto}H^1(B,\Z).$$

\smallskip
The mirror symmetry works ``best'' for abelian varieties $A$ which are
obtained by a version of the {\bf G}--construction of Gerritzen (section 10).

\smallskip
\noindent{5.} If the abelian variety $A$ is obtained by the {\bf G}--
construction,
then for any $\omega _A\in C_A$ the algebraic pair $(A,\omega _A)$
has a mirror symmetric pair $(B,\omega _B),$ such that $B$ is also
obtained by the {\bf G}--construction. Moreover, in this case we can
describe the
isomorphism $\beta $ explicitly. Namely, $\beta $ is given by an element
in $H^*(A\times B,\Z)$ which is the Chern character of some natural line
bundle on a certain real subtorus of $A\times B$ of dimension $3n$
(where $n=dimA=dimB).$

\medskip
\noindent{\bf 0.3.} Let us describe the two $\Q$--algebraic groups attached to any abelian variety $A,$ which act
on the
total cohomology $H^*(A,\Q)$ and commute with each
other.

\medskip
\noindent{\bf 0.3.1.} The first group is the classical {\it Hodge group} or
the {\it special Mumford-Tate group}. Let us recall its definition.

Put $\Ga =\Ga _A =H_1(A,\Z),$ $V=\Ga _{\R}=H_1(A,\R).$ The complex structure
on $A$ induces the operator of complex structure $J_A$ on $V.$
Consider the homomorphism of algebraic $\R$-groups
$$h_A:\S^1\lto \GL(V),\ \ \ \ h_A(e^{i\theta})={\rm cos}(\theta)\cdot Id
+{\rm sin}(\theta )\cdot J_A,$$
so that $h_A(e^{i\pi/2})=J_A.$
Then
$Hdg_{A,\Q}$ is defined as the smallest $\Q$--algebraic subgroup $G$ of
$\GL(H_1(A,\Q))$ such that $h_A(\S^1)\subset G(\R).$ By functoriality
$Hdg_{A,\Q}$
acts on $H^1(A,\Q)$ and on the total cohomology
$H^*(A,\Q)=\La ^{\cdot}H^1(A,\Q).$ We call this group ``horizontal'' because
it preserves each subspace $H^k(A,\Q).$

\medskip
\noindent{\bf 0.3.2.} Perhaps the main point of our work is the consideration
of the second (``vertical'') algebraic group which we view as the mirror
image of the Hodge group. This second group is the Zariski closure
$\overline{Spin(A)}$ in $\GL(H^*(A,\Q))$ of a certain discrete subgroup
$Spin(A)\subset \GL(H^*(A,\Q)).$

Let $D^b(A)$ be the bounded derived category of coherent sheaves on $A.$
Consider its group $Auteq(D^b(A))$ of exact autoequivalences.
We construct in 4.3 a natural representation
$$\rho _A:Auteq(D^b(A))\lto  \GL(H^*(A,\Z))$$
and denote its image by $Spin(A).$ The elements of $Spin(A)$ act
by algebraic correspondences, hence commute with the Hodge group and
preserve the Hodge verticals $\oplus _{p-q=fixed}H^{p,q}(A,\C).$ Hence
the same is true for the algebraic group $\overline{Spin(A)}.$
Let us discuss some properties of the groups $Spin(A)$ and
$\overline{Spin(A)}.$

\medskip
\noindent{\bf 0.3.3.} Let $\wh{A}$ be the dual abelian variety and $\Ga_{\wh{A}}$ be its first homology lattice that can be identified with ${\Ga}_{A}^*:=\Hom(\Ga_{A}, \Z).$
Consider
the lattice
$$\La :=\Ga_{\vphantom{\wh{A}}A}\oplus \Ga_{\wh{A}}$$
with the canonical bilinear symmetric form
$$Q((a,b),(c,d))=b(c)+d(a).$$

Let $\SO (\La ,Q)$ be the corresponding special orthogonal group. Let
$Spin(\La ,Q)$ be the corresponding spinorial group. We have the
canonical exact sequence
$$0\lto \Z /2\Z \lto Spin(\La ,Q) \lto \SO (\La, Q).$$
Note that the group $Spin(\La ,Q)$ acts naturally on the total
cohomology group
$\La ^{\cdot}\Ga ^*_A=H^*(A,\Z).$

For any abelian variety $A,$ there are two representations of
$Auteq(D^{b}(A)).$ One, in $H^*(A, \Z),$ is the representation $\rho_{A}$ considered above, in (0.3.2). The existence of the other, in $\La$ is implied by the explicit description of
$Auteq(D^{b}(A))$ given in \cite{Or2}.
In fact, a commutative diagram exists:
$$
\begin{array}{ccc}
Spin(A) & \hookrightarrow & Spin(\La ,Q) \\
\downarrow & & \downarrow \\
\U (A) & \hookrightarrow & \SO (\La ,Q)
\end{array}
$$
where the horizontal arrows are compatible with the actions on $H^*(A,\Z)$
and $\La $ respectively. The group $\U (A)$ was first introduced independently by
S.Mukai \cite{Muk2} and by A.Polishchuk \cite{Po2}, so we call it the
Mukai-Polishchuk
group.

Let $\overline{\U (A)}\subset \GL (\La _{\Q})$ be the $\Q-$algebraic subgroup
which is the Zariski closure of $\U (A)$ in $\GL (\La _{\Q}).$ Thus the groups
$\overline{Spin(A)}$ and $\overline{\U (A)}$ are isogeneous.

\medskip
\noindent{\bf 0.3.4.} The algebraic groups $\overline{Spin(A)}$ and
$\overline{\U (A)}$ are semisimple (7.2.1).

\medskip
\noindent{\bf 0.3.5.} Consider the set $C_A=C_A^+\sqcup C_A^-\subset
NS_A(\C )$ as in 0.1. Both $C_A^+$ and $C_A^-$ can be considered as Siegel domains
of the first kind: the Lie group $\overline{\U (A)}(\R )$ acts naturally
on $C_A$ preserving $C_A^+$ and $C_A^-.$ Moreover, $C_A^+$ and $C_A^-$
are single $\overline{\U (A)}(\R )$--orbits and the stabilizer $K_{\omega }$
of a point $\omega \in C_A$ is a maximal compact subgroup of the
semisimple Lie group $\overline{\U (A)}(\R ).$ Further, for each
$\omega \in C_A$ there is a natural choice of an element $I_{\omega }\in
K_{\omega }$ (\ref{pred}) such that $I_{\omega }$ defines a complex structure on the
space $\La _{\R }$ (0.3.3).

\medskip
\noindent{\bf 0.4.} Now we are ready to give the main definition of mirror
symmetry
for algebraic pairs.

Let $(A,\omega _A)$ be an algebraic pair, $\La_A=\Ga_{\vphantom{\wh{A}}A}\oplus \Ga_{\wh{A}}$
with the symmetric form $Q_A$ as in 0.3.3. Note that the Hodge group
$Hdg_{\vphantom{\wh{A}}A,\Q }=Hdg_{A\ts\wh{A},\Q}$ is naturally a subgroup of $\SO (\La_{A,\Q},Q_{A,\Q}).$
Moreover, it commutes with $\overline{\U (A)}.$  Thus we obtain two
commuting complex structures on $\La _{A,\R}$:
$J_{A\ts\wh{A}}\in Hdg_{\vphantom{\wh{A}}A,\Q}(\R)=Hdg_{A\ts\wh{A},\Q}(\R)$ and $I_{\omega_{A}}\in \overline{\U (A)}(\R).$

\medskip
\noindent{\bf 0.4.1  Definition.} {\it We call algebraic pairs $(A,\omega _A)$
and $(B, \omega _B)$ {\sf mirror symmetric} if there is given an isomorphism
of lattices
$$\alpha :\La _A\stackrel{\sim}{\lto }\La _B,$$
which identifies the forms $Q_A$ and $Q_B$ and satisfies the following conditions:
$$\alpha _{\R}\cdot J_{A\ts\wh{A}}=I_{\omega_{B}}\cdot \alpha _{\R},$$
$$\alpha _{\R}\cdot I_{\omega_{A}}=J_{B\ts\wh{B}}\cdot \alpha _{\R}.$$}

\medskip
\noindent{\bf 0.4.2.} Note that if we identify $\La _A$ and $\La _B$ by means
of $\alpha $ then
$$Hdg_{A,\Q }\subseteq \overline{\U (B)},\qquad
Hdg_{B,\Q}\subseteq \overline{\U (A)}.$$

\medskip
\noindent{\bf 0.5.} Actually we work in a more general context. Namely,
we consider {\it weak pairs} $(A,\omega _A),$ where $A$ is a {\it complex
torus}
and $\omega _A=\eta _1+i\eta _2\in NS_A(\C)$ is such that $\eta _2$
is nondegenerate. (Thus an algebraic pair is in particular a weak pair.)
 And we define the notion of mirror symmetry for weak pairs. However we
prove that if two weak pairs are mirror symmetric and one of them is algebraic
then so is the other.

\medskip
\noindent{\bf 0.6.} For an abelian variety $A$ the algebraic group
$\overline{Spin(A)}$ or rather its Lie algebra has a different description.
It turns out to be isomorphic (as a Lie subalgebra of ${\frak gl}(H^*(A,\Q)))$
to the {\it Neron-Severi} Lie algebra ${\frak g}_{NS}(A)$ defined in \cite{LL}.

\medskip
\noindent{\bf 0.6.1.} Let us recall the definition of ${\frak g}_{NS}(X)$ for
a smooth complex projective variety $X.$
If $\kappa \in H^{1,1}(X)\cap H^2(X,\Q)$ is an ample class, then
cupping with it defines an operator $e_{\kappa }$ in the total cohomology
$H^*(X)=H^*(X,\C)$ of degree 2 and the hard Lefschetz theorem asserts
that for $s=0,1,...,n,$ $e^s_{\kappa }$ maps $H^{n-s}(X)$ isomorphically
onto $H^{n+s}(X).$ As is well known, this is equivalent to the existence
of a (unique) operator
$f_{\kappa}$ on $H^*(X)$ of degree $-2$ such that the commutator
$[e_{\kappa },f_{\kappa }]$ is the operator $h$ which on $H^k(X)$ is
multiplication by $k-n.$
The elements $e_{\kappa },f_{\kappa },h$ make up a Lie subalgebra ${\frak g}
_{\kappa }$ of ${\frak gl}(H^*(X))$ isomorphic to $sl(2).$ Define the
{\it Neron-Severi} Lie algebra ${\frak g}_{NS}(X)$ as the Lie subalgebra
of ${\frak gl}(H^*(X))$ generated by ${\frak g}_{\kappa}$'s with $\kappa $ an
ample class. This Lie subalgebra is defined over $\Q$ and is evenly
graded by the adjoint action by the semisimple element $h.$ The Lie
algebra ${\frak g}_{NS}(X)$ is semisimple.

\medskip
\noindent{\bf 0.6.2.} Note that the equality
$$Lie\overline{Spin(A)}\ \ =\ \ {\frak g}_{NS}(A)\quad \quad \quad (*)$$
implies one of the standard conjectures of Grothendieck for $A.$ Namely, the
algebraicity of the operator $f_{\kappa }$ (0.5.1) that has been proved for abelian
varieties by Kleiman \cite{Kl}.

\medskip
\noindent{\bf 0.6.3.} {\it Question.} Let $X$ be a smooth complex projective
variety. Assume that the canonical sheaf $K_X$ is trivial. Does the analogue
of (*) hold for $X$?

\smallskip
Note that if $K_X$ or $K^{-1}_X$ is
ample then the analogue of (*) cannot hold. Indeed, by a theorem in \cite{BO}
the group $Auteq(D^b(X))$ is then generated by $Aut(X),$ the shift
operator [1], and by operators of tensoring with a line bundle on $X.$
Thus (if $Aut(X)$ is trivial) $Auteq(D^b(X))$ is abelian, whereas the
Lie algebra ${\frak g}_{NS}(X)$ is semisimple.

\medskip
\noindent{\bf 0.7.} Our picture of the mirror symmetry provides an explanation
of the phenomenon in Hodge theory, which was noticed long ago (see for
example p.68 in \cite{Gr}).
Namely, consider a variation of Hodge structures which degenerates along a
divisor. Then the logarithm of the monodromy operator has properties
which are similar to the properties of a Lefschetz operator on the
cohomology. But in our construction this logarithm of the monodromy
is transformed by the mirror symmetry to a Lefschetz operator indeed.

\medskip
\noindent{\bf 0.8.} A similar picture exists for hyperkahler manifolds
and will be described in our next paper.

\medskip
\noindent{\bf 0.9.} Our notion of mirror symmetry for abelian varieties
is compatible with other (conjectural) pictures of mirror symmetry. We will
explain the compatibility with the ideas of Kontsevich \cite{Ko} on one hand and
with the $T$--duality picture of Strominger, Yau and Zaslow \cite{SYZ} on the
other hand.

\medskip
\noindent{\bf 0.9.1.} In \cite{Ko} Kontsevich proposes the following thesis.
If a symplectic manifold $(V, \omega)$ is mirror symmetric to an algebraic
variety $W,$ then the Fukaya category $F(V,\omega )$ of $(V,\omega )$ is
equivalent to the derived category $D^b(W)$ of coherent sheaves on $W.$
This has been proved for elliptic curves in \cite{PZ}.

Assume that the algebraic pairs $(A,\omega _A)$ and $(B,\omega _B)$ are
mirror symmetric (0.4.1). Let $\omega _A=\eta _1+i\eta _2 \in C_A.$ Then
$\eta _2$ is a Kahler form on $A$ and in particular $(A,\eta _2)$ is a
symplectic manifold. (The class $\eta _1$ plays the role of a B-field).
 The Hodge group $Hgd_{A,\Q}$ acts on $H_1(A)$ preserving
the form $\eta _2.$ Hence its arithmetic subgroup
$$Hdg_{A}(\Z):=\{g\in Hdg_{A,\Q}|\ \ g(\Lambda _A)\subset \Lambda _A\}$$
acts by symplectomorphisms of the symplectic manifold $(A,\eta _2).$ Thus
$$Hdg_A(\Z)\subset AuteqF(A,\eta _2).$$

But our definition of mirror symmetry implies the natural inclusion of
groups
$$Hdg_A(\Z)\subset \U(B)$$
(0.3.3, 0.4.2), where the group $\U(B)$ captures the essential part of
$AuteqD^b(B).$ Thus we do not establish the equivalence of categories
$F(A,\eta _2)$ and $D^b(B)$ as suggested by Kontsevich, but rather compare
their groups of autoequivalences.

\medskip
\noindent{\bf 0.9.2.} Let us show the compatibility with the mirror symmetry
proposed in \cite{SYZ} and more specifically in \cite{Mo} (see also \cite{Gro1}, \cite{Gro2}).
Let algebraic pairs $(A,\omega _A)$ and $(B,\omega _B)$ be mirror symmetric.
Assume that the pair $(A,\omega _A)$ is obtained by the {\bf G}--construction
(10.2.1) (this is so for ``most'' pairs). Then there exists a real subtorus
$T\subset A\times B$ of dimension $3n$ ($n=\dim_{\C}A=\dim_{\C}B$) with the
following property: The fibres of the projection of $T$ onto $B$ (resp. $A$)
are special Lagrangian subtori of $A$ (resp. $B$) of dimension $n$
(10.3.2, 10.7). The existence of such a correspondence is part of the
``geometric mirror symmetry'' of Morrison (Definition 4 in \cite{Mo}).
Moreover, there is a natural complex line bundle $L$ on $T$ such that
its Chern class $c_1(L)$ considered as an element in $H^*(A\times B)
\simeq Hom(H^*(A),H^*(B))$ is the isomorphism $\beta $ mentioned in 0.2.

\medskip
\noindent{\bf 0.10.} Let us briefly discuss the contents of each section.

Section 1 contains some background material on abelian varieties and
complex tori. In section 2 we recall the Hodge group of a complex torus
and its basic properties. Section 3 contains a review of the Clifford
algebra (of a certain quadratic form which we define there) and its spinorial
representation. This Clifford algebra plays the key role in the definition
of the isomorphism $\beta $ discussed in 0.2.

Section 4 is devoted to the bounded derived category $D^b(X)$ of
coherent sheaves on a smooth projective variety $X.$ We define the group
$Auteq(D^b(X))$ of autoequivalences of $D^b(X)$ and construct a canonical
representation
$$\rho _X:Auteq(D^b(X))\lto \GL (H^*(X,\Q)).$$
We then show that if $X=A$ is an abelian variety then $Im\rho _A\subset
\GL (H^*(A,\Z)).$ The Mukai-Polishchuk group is introduced and then
the main effort is applied to prove the commutativity of the diagram in 0.3.3
(see Proposition 4.3.7).

In section 5 we discuss various $\Q$--algebraic subgroups of
$\GL (H_1(A\times \wh{A},\Q)),$ which enter the picture of mirror symmetry.
In particular we recall the algebraic group $\U_{A,\Q}$ which was introduced
by A. Polishchuk and present his result on the group of
$\R$--points of $\U_{A,\Q}.$ This group $\U_{A,\Q}$ is closely related to
the algebraic groups $\overline{Spin(A)}$ and $\overline{\U(A)}.$
Namely, the semisimple group $\overline{\U(A)}$ is a subgroup of the
reductive group $\U_{A,\Q}$ which consists (up to isogeny) of all
noncompact factors of $\U_{A,\Q}$ (7.2.1).

In section 6 we recall the Neron-Severi Lie algebra ${\frak g}_{NS}(X)$
of a smooth projective variety $X$ and describe ${\frak g}_{NS}(A)$
for an abelian variety $A$ following \cite{LL}. Then in section 7 we establish
the relationship between this Lie algebra and the group
$Auteq(D^b(A)).$

In section 8 we define a natural action of the Lie group $\U_{A,\Q}(\R)$
on the Siegel domain $C^+_A$ (and $C^-_A$). We then show that this action
is transitive and the stabilizer of a point is a maximal compact subgroup
of $\U_{A,\Q}(\R).$ To each $\omega \in C_A$ we associate a natural
element $I_{\omega }\in \U_{A,\Q}(\R),$ which is the main ingredient in
the definition of mirror symmetry.

Finally, in section 9 we define the notion of mirror symmetry for complex
tori and abelian varieties and study its properties.

Section 10 contains a variant of the {\bf G}--construction and its
relations with the mirror symmetry.

\medskip
\noindent{\bf 0.10.} It is our pleasure to thank Joseph Bernstein, Alexei Bondal, Yuri Manin,
Andrei Tyurin and Yuri Zarhin for useful discussions.
We also were inspired by the pleasant
atmosphere at the Moscow Steklov Institute, the Moscow Independent
University, and cafe houses on the Old Arbat street.

\sec{Some preliminaries on complex tori and abelian varieties.}

\noindent{\bf 1.1.}
Let $\Ga \cong \Z^{2n}$ be a lattice, $V=\Ga \otimes \R \cong \R^{2n}$
and $J\in \End(V),$ s.t. $J^2=-1.$ That is $J$ is a complex structure on $V.$
This way we obtain an $n$-dimensional complex torus
$$A=(V/\Ga,J).$$
Note the canonical
isomorphisms
$$\Ga =H_1(A,\Z),\quad V=H_1(A,\R).$$
Sometimes we will add the subscript ``A'' to the symbols $\Ga ,$ $V,$ $J.$

Given another complex torus  $B=(V_B/\Ga_B,J_B),$ the group $\Hom(A,B)$
consists of
homomorphisms
$f:\Ga _A \to \Ga _B$
such that
$$
J_B \cdot f_\R = f_{\R} \cdot J_A   : V_A
\lto V_B.
$$
Thus the abelian group
 $\Hom (A,B)$ can be considered as a subgroup of $\Hom(\Ga _A, \Ga _B).$

\medskip
\noindent{\bf 1.2.}
One has the dual torus $\wh A$ defined as follows. Put $\Ga^* = \Hom
_{\Z}(\Ga,\Z),$ $ V^* =\Ga^*\otimes \R=\Hom (V,\R)$ and
$\wh J:V^*\stackrel{\sim}{\to} V^*,$ s.t. $(\wh{J}w)(v)=w(-Jv)$ for
$v\in V,$ $w\in V^*.$ Then by definition
$$\wh A=(V^*/\Ga^*, {\wh{J}}).$$

\medskip
\noindent{\bf 1.3.}
Denote by  $\Pic_A$ the Picard group of $A.$ Let
$\Pic^0_A\subset \Pic_A$ be
the subgroup
of line bundles with the trivial Chern class. It has a natural structure of a complex torus.
 Moreover, there exists a natural
isomorphism of complex tori
$$\wh{A} \cong \Pic^0_A.$$
Every line bundle $L$ on $A$ defines a morphism $\varphi _L:A\to \wh A$
by the formula
$$\varphi _L(a)=T_a^*L\otimes L^{-1}.$$
(Here $T_a:A\to A$ is the translation by $a.$) We have $\varphi _L=0$ iff
$L\in \Pic^0_A$ and $\varphi _L$ is an isogeny if $L$ is ample. Thus
the correspondence $L\mapsto \varphi _L$ identifies  the Neron-Severi
group $NS_A:=\Pic_A/\Pic^0_A$ as a subgroup in $\Hom(A,\wh A).$ Also
$NS_A$ is naturally a subgroup of $H^2(A,\Z)$:
 to a line bundle $L$ there corresponds its first Chern class that can be considered as
a skew-symmetric bilinear
form  on $\Ga .$ Put $c_1(L)=c.$ Then the
morphism $\phi _L$ is given by the map
$$V_{\vphantom{\wh{A}}A}\to V_{\wh A},\ \ \ \ v\mapsto c(v,\cdot ).$$

We will identify $NS_A$ either as a subgroup of $\Hom(A,\wh A)$ or
$\Hom(\Ga_{\vphantom{\wh{A}}A}, \Ga_{\wh A})$ or as a set of (integral) skew-symmetric forms
$c$ on $\Ga _A$ such that the extension $c_{\R}$ on $V_A$ is $J$-invariant.
Put $NS_A(\R)=NS_A\otimes \R.$ We will denote by
$NS_A(\R)^0\subset NS_A(\R)$ the subset consisting of nondegenerate forms.

\medskip
\noindent{\bf 1.4.}
There is
a unique line bundle
$P_{A}$ on the product
$A \times \widehat A$ such that

i) for any point
$\alpha \in \widehat A$ its restriction
$P_\alpha$ on
$A\times \{\alpha\}$
represents an element of
$\Pic^0_{A},$ corresponding to
$\alpha$ with respect to the fixed isomorphism
$\Pic^0_{A} =\widehat A;$

ii) the restriction
$P|_{\{0\} \times \widehat A}$ is trivial.

Such
$P_{A}$ is called Poincare line bundle.

Given a morphism of complex tori
$f: A \to B,$ the dual morphism
$\widehat f: \widehat B \to \widehat A
,$ is defined.
Pointwise, for
$\alpha\in \wh{A}$
and
$\beta\in \wh{B}$
one has
$\wh{f}(\beta)=\alpha$
if and only if
the line bundle
$f^{*}P_{\beta}$
on
$A$
coincides with
$P_{\alpha}.$

The double dual torus $\wh{\wh A}$ is naturally identified with $A$
by means of the Poincare line bundle on $A\times \wh A$ and
$\wh A\times \wh{\wh A}.$
In other words,there exists a unique isomorphism
$\kappa_{A}: A\stackrel{\sim}{\to} \wh{\wh{A}}$
such that under
the isomorphism
$1\times\kappa_A: \wh{A}\times A\stackrel{\sim}{\lto} \wh{A}\times \wh{\wh{A}}$
the Poincare line bundle pulls back to the Poincare line bundle:
$
(1\times\kappa_A)^{*}P_{\widehat{A}} \cong
P_{\vphantom{\wh{A}}A}.
$

Thus
$\widehat{}$ is an antiinvolution on the category of
complex tori,
i.e. a contravariant functor whose square is
isomorphic to the identity:
$ \kappa: Id \stackrel{\sim}{\to}\widehat{\widehat{}}.$
It is known that $\wh \varphi _L=\varphi _L$ for any
$L\in \Pic_A$ (see next remarks).

\medskip
\noindent{\bf 1.5 Remark.}
Let $A$ be a complex torus. Consider the
Poincare line bundles $P_{\vphantom{\wh{A}}A}$ and $P_{\wh A}$ on $A\times \wh{A}$ and $\wh A\times
\wh{\wh A}$ respectively.
%
%
%
%
Choose a basis $l_1,...,l_{2n}$ of $\Ga_{\vphantom{\wh{A}}A}$ and the dual basis
$l_1^*,...,l^*_{2n}$ of $\Ga_{\wh A}={\Ga}^*_{A}$ (1.2). Let $l^{**}_1,...,l^{**}_{2n}$ be
the basis of $\Ga_{\wh{\wh A}}={\Ga}^*_{\wh{A}}$ dual to $l^*_1,...,l^*_{2n}.$
The isomorphism $\kappa_{A}:A\stackrel{\sim}{\lto}\wh{\wh{A}}$ induces an identification of $\Ga_{\vphantom{\wh{\wh{A}}}A}$ and $\Ga_{\wh{\wh A}}$
that takes $l_i$ to $-l^{**}_i$ (!), due to the fact that the
forms  $c_1(P_{\vphantom{\wh{A}}A})$ and $c_1(P_{\wh A})$ are skew-symmetric.


\medskip
\noindent{\bf 1.6 Remark.}
Let $f:A\to B$ be a morphism of complex tori and $\wh{f}:\wh{A}\to \wh{B}$ be the dual morphism
Choose basises $l_1,...,l_{2n}$ , $m_{1},..., m_{2n}$
for $\Ga_A,$ $\Ga_B$ and the dual basises
$l_1^*,...,l_{2n}^*$ , $m_{1}^*,..., m_{2n}^*$
for $\Ga_{\wh{A}},$ $\Ga_{\wh{B}}.$
Denote by $F$ and $\wh{F}$ the matrices of linear maps
$f:\Ga_{\vphantom{\wh{A}}A}\to \Ga_{\vphantom{B}B}$ and
$\wh{f}: \Ga_{\wh{B}}\to \Ga_{\wh{A}}$ in these basises.
The matrices $F$ and $\wh{F}$ are transposes of each other,
i.e. $\wh{F}=F^{t}.$

Now, let $f: A\to \wh{A}.$ Using the isomorphism $\kappa,$ we will consider $\wh{f}$ as morphism
form $A$ to $\wh{A}$ too.
If, as above, $F$ and $\wh{F}$ are the matrices of linear maps
$f, \wh{f}:\Ga_{\vphantom{\wh{A}}A}\to \Ga_{\wh{A}},$ then they are skew-transposes of each other,
i.e. $\wh{F}=-F^{t}.$ This  immediately follows from the previous remark.

Thus, in particular, $\wh \varphi _L=\varphi _L$ for any $L\in \Pic_A.$
%
%
%
%
%
%
%
%

\medskip
\noindent{\bf 1.7.}
A complex torus $A=(V/\Ga ,J)$ is algebraic, i.e. an abelian variety,
iff there exists $c\in NS_A$ such that the symmetric bilinear form
$c_{\R}(J\cdot,\cdot)$ on $V$ is positive definite.

For the rest of this section 1 we assume that $A$ is an abelian variety.

The endomorphism ring
$\End(A)$
is a finitely generated
$\Z$--module,
and we denote , as usual,
$\End^0(A):=\End(A) \otimes \Q.$
By the Poincare reducibility theorem the variety
$A$
is isogeneous to a product of simple abelian varieties. So
$\End^0(A)$
is a semisimple finite
$\Q$--algebra. Put $\Ga _{A,\Q}=\Ga _A \otimes \Q.$

Let $L\in \Pic _A$ be ample. Then the induced map
$$\varphi _L:\Ga _{\vphantom{\wh{A}}A,\Q} \to \Ga _{\wh{A},\Q}$$
is an isomorphism. The map
$$':\End^0(A)\to\End^0(A),\ \ \ \ a\mapsto \varphi _L^{-1}\cdot
\wh a \cdot \varphi _L$$
is a positive (anti-)involution of $\End^0(A)$ called the Rosati involution.

\medskip
\noindent{\bf 1.8.}
Assume that
$A$
is simple, so that
$\End^0(A)$
is a division algebra. Finite division
$\Q$--algebras
with positive involutions were classified by Albert.
Below is his classification which consists of four cases.

Let us introduce some notation. Put
$F=\End^0(A),$ $K$ - the center of
$F.$
Fix a Rosati involution
$'$ of $F.$
The involution
that it induces on
$K$ is
independent of
$'$ :
for any embedding of
$K$
in
$\C$
it is given by complex conjugation. The fixed
subfield
$K_0 = K'$
is totally real. Put
$[F:K]=d^2, [K_0:\Q]=e_0.$
\begin{list}{\Roman{tmp}.}%
{\usecounter{tmp}}
\item The case of totally real multiplication.\\
Here $F=K=K_0.$
\item The case of totally indefinite quaternion multiplication.\\
Here
$K=K_0$
and
$F$ is a
$K_0$-form of
$M(2) $:
there is
an
$\R$--algebra
isomorphism
$F\otimes_\Q \R \cong M(2,\R)^{e_0}$
such that the involution corresponds to the transpose in every summand.
\item The case of totally definite quaternion multiplication.\\
Here also
$K=K_0$
and
$F$
is a
$K_0$--form
of the quaternion algebra
$\Qu$
over
$K_0$:
there is an
$\R$--algebra
isomorphism
$F \otimes_\Q \R \cong \Qu^{e_0}$
such that the involution corresponds to the
quaternion conjugation in every summand.
\item The case of totally complex multiplication.\\
The field
$K$
has no real embedding (so is a totally imaginary quadratic
extension of
$K_0$) and
$F$
is a
$K$--form of
$M(d)$ :
there is an
$\R$--algebra
isomorphism
$F \otimes_\Q \R \cong M(d,\C)^{e_0}$
such that the involution corresponds to the
conjugate transpose in every summand.
\end{list}

\sec{The Hodge group $Hdg_{A,\Q}.$}

\noindent{\bf 2.1.}
Let
$A=(V/\Ga,J)$
be a complex torus.
Consider the homomorphism of algebraic $\R$-groups
$$h_A:\S^1\lto \GL(V),\ \ \ \ h_A(e^{i\theta})={\rm cos}(\theta)\cdot Id
+{\rm sin}(\theta )\cdot J,$$
so that $h_A(e^{i\pi/2})=J.$ This defines a representation of $\S^1$ on the
exterior algebra $\La^{\cdot}H_1(A,\R),$ hence on the total (co-)homology of
$A.$
Decomposing the representation of $\S^1$ on the k-th cohomology space with respect to the
characters of
$\S^1$
$$
\S^1\lto \C^*,
\quad
e^{i\theta}\mapsto e^{i(q-p)\theta}, \quad p+q=k,
$$
we get the usual Hodge decomposition
$$
H^k (A, \C)=\bigoplus_{p+q=k} H^{p,q}(A, \C)
$$

\medskip
\noindent{\bf 2.2 Definition.}
{\it The {\sf  Hodge group (or special Mumford-Tate group)}
$Hdg_{A, \Q}$ of
$A$
is defined as a minimal algebraic $\Q$-subgroup
$G$ of
${\rm GL}(H_1 (A, \Q))$
such that
$h_A (\S^1)\subset G(\R).$}

\smallskip

The group
$Hdg_{A, \Q}$
acts naturally on the homology
$H_k (A, \Q)$
and the cohomology
$H^k (A, \Q).$

It is clear that
$Hdg_{A, \Q}$
acts trivially on the Hodge spaces
$H^{p,p}(A, \Q):= H^{p,p}(A, \C)\cap H^{2p}(A, \Q).$

Consider the dual torus $\wh A.$ We have canonical identifications
$$Hdg_{\vphantom{\wh{A}}A,\Q}=Hdg_{{\wh A},\Q}=Hdg_{A\times {\wh A},\Q}.$$
The representation of $Hdg_{\vphantom{\wh{A}}A,\Q}$ on $\Ga _{{\wh A},\Q}$ is the
contragradient of its representation on $\Ga _{A,\Q}.$
Depending upon a context, we will view
$Hdg_{A,\Q}$ as a subgroup of $\GL(\Ga_{\vphantom{\wh{A}}A,\Q})$
or  $\GL(\Ga _{{\wh A},\Q})$
or $\GL (\Ga _{A,\Q}\oplus \Ga _{{\wh A},\Q}).$

\medskip
\noindent{\bf 2.3.}
The following facts about
$Hdg_{A, \Q}$
are known
(see \cite{Mum2}, \cite{De}).

\medskip
\noindent{\bf 2.3.1 Theorem.} {\it Assume that $A$ is an abelian variety.
\begin{list}{\alph{tmp})}%
{\usecounter{tmp}}
\item
$Hdg_{A, \Q}$
is a connected reductive algebraic $\Q$-group.
\item
 $h_A (-1)$
is in the center of
$Hdg_{A, \Q}.$
\item
 The involution
$Ad(h_{A}(i))$ on
$Hdg_{A, \Q}(\R)^{0}$
is a Cartan involution,
i.e. its fixed subgroup is a maximal compact subgroup
$K$ of
$Hdg_{A, \Q}(\R)^{0}.$
\item
 The symmetric space
$Hdg_{A, \Q}(\R)^{0}/K$
is of hermitian type.
\item
 The reductive group $Hdg_{A,\Q}$ has no simple factors of exceptional
type.
\end{list}}

\sec{The Clifford algebra and the spinor representation.}

\noindent{\bf 3.1.}
Fix a lattice $\Ga \cong \Z^{2n}$ for some $n\geq 1.$ (In the future
applications $\Ga $ will be $\Ga _A$ for a complex torus $A$).
Put $\La =\Ga \oplus \Ga ^*$ with the canonical symmetric bilinear form
$Q:\La\times\La\lto\Z$ defined by
$$Q((a_1,b_1),(a_2,b_2))=b_1(a_2)+b_2(a_1).$$
This form is even and unimodular. Moreover, $Q\simeq U^{2n},$ where
$U=\left(
\begin{smallmatrix}
0 & 1\\
1 & 0
\end{smallmatrix}
\right).$

Consider the Clifford algebra $Cl(\La,Q)$ defined as follows
$$Cl(\La,Q):=T(\La)/<x\otimes x=\frac{1}{2}Q(x,x)>,$$
where $T(\La)$ is the tensor algebra of the lattice $\La.$

\medskip
\noindent{\bf 3.2.}
Choose a basis $l_1,...,l_{2n}$ of $\Ga $ and the dual basis $x_1,...,x_{2n}$
of $\Ga ^*.$ Let $I\subset Cl(\La,Q)$ be the left ideal generated by the
product ${\bf l}:=l_1\cdots l_{2n}.$ The left multiplication in
$Cl(\La,Q)$ defines a homomorphism of algebras
$$\sigma :Cl(\La,Q)\lto\End_{\Z}(I).$$

\medskip
\noindent{\bf 3.2.1 Proposition.} {\it
\begin{list}{\alph{tmp})}%
{\usecounter{tmp}}
\item
$Cl(\La,Q)$
is a free abelian group of rank
$2^{4n}$
with the basis consisting of monomials
$$
x_{j_1}\cdots x_{j_l}l_{i_1}\cdots l_{i_k},
\quad 1\le k,s\le 2n,
\quad i_1<\cdots<i_k , j_1<\cdots <j_s .
$$
In particular, $\La $ is naturally a subgroup of $Cl(\La,Q).$
\item
$Cl(\La,Q)$
is a
$\Z_2$-graded algebra,
$Cl(\La,Q)=Cl^{+}(\La,Q)\op Cl^{-}(\La,Q),$
where
$Cl^{+}(\La,Q)$
(resp.
$Cl^{-}(\La,Q)$)
is the span of monomials with even (resp. odd) number of elements.
$Cl^{+}(\La,Q)$
is a subalgebra of
$Cl(\La,Q).$

\item
 The ideal
$I$
is independent of the choice of the basis
$l_1,...,l_{2n}$
of
$H_1 (A, \Z).$

\item
 The ideal
$I$
has a
$\Z$-basis consisting of the monomials
$$
x_{j_1}\cdots x_{j_s}l_1\cdots l_{2n},
\quad
1\le s \le 2n,
\quad
j_1<\cdots <j_s .
$$
Hence $I\cong \La ^{\cdot}\Ga ^*$ as abelian groups.

\item
 The homomorphism
$\sigma: Cl(\La,Q)\lto {\rm \End}_{\Z}(I)$
is an isomorphism, hence
$Cl(\La,Q)\cong{\rm M}_{2^{2n}}(\Z).$
\end{list}}
\pr
Assertions a), b) are standard and easy

\noindent
c) is obvious since $l_i,$ $l_j$ anticommute
for $i\neq j.$

\noindent
d) follows from a).

\noindent
Let us prove e). Since $Cl(\La,\Q)$
and $\End_Z(I)$ are both free $Z$-modules of rank $2^{4n}$ it suffices
to prove that $\sigma $is a surjection. For a subset
$p=\{i_1,...,i_k \}\subset [1,2n],$ $i_1\leq...\leq i_k,$ denote by ${\bf l}_p$
(resp. ${\bf x}_p$) the ordered monomial $l_{i_1}\cdots l_{i_k}$
(resp. $x_{i_1}\cdots x_{i_k}$). Let $p'=[1,2n]-p$ be the complementary set.
One immediately checks that for subsets $p,h,m\subset [1,2n]$
$${\bf x}_h{\bf l}{\bf x}_{p'}({\bf x}_m{\bf l})=
\left\{
\begin{array}{cc}
\pm {\bf x}_h{\bf l} & {\rm if}\ \ \ \ p=m\\
0 & {\rm otherwise}
\end{array}\right\},$$
where ${\bf l}=l_1 \cdots l_{2n}.$
This implies the surjectivity of $\sigma $ in view of d) above
(see also \cite{Ba}).
$\Box$

\medskip
\noindent{\bf 3.3 Corollary.} {\it Let $\La =M_1\oplus M_2$ be a decomposition
such that $M_1,M_2$ are (maximal) $Q$-isotropic sublattices. Let
$m_1,...,m_{2n}$ be a basis of $M_1,$ ${\bf m}=m_1 \cdots m_{2n}$ and
$I'=Cl(\La,Q){\bf m}$ be the corresponding left ideal. There exists a
unique (up to  $\pm 1$) isomorphism of left $Cl(\La,Q)$-
modules $I$ and $I'.$}

\smallskip
\pr
The form $Q$ identifies $M_2$ with the dual $M_1^*$ of $M_1.$ Moreover,
the form $Q'$ on $\La =M_1\oplus M_1^*$ defined by
$$Q'((m_1,n_1),(m_2,n_2))=n_1(m_2)+n_2(m_1)$$
coincides with $Q,$ hence gives rise to the same Clifford algebra.
Therefore we may apply the previous proposition (parts d),e)) to the
ideal $I'$ instead of $I.$ Since $\End_Z(I)=Cl(\La,Q)=\End_Z(I')$ there
exists a unique (up to a sign) isomorphism $I\cong I'$ of left
$Cl(\La,Q)$-modules. This proves the corollary.

\medskip
\noindent{\bf 3.4.}
Let us recall the spinorial group of the Clifford algebra.

There exists a unique involution
$'$
on
$Cl(\La,Q)$
that is the identity on
$\La.$
In terms of the monomial basis as in Proposition 3.2.1a) it is
defined by rule
$$
(x_{j_1}\cdots x_{j_s}l_{i_1}\cdots l_{i_k})'
=l_{i_k}\cdots l_{i_l}x_{j_s}\cdots x_{j_1}
$$

\medskip
\noindent{\bf 3.4.1 Definition.}
The {\sf  spinorial group} $Spin(\La,Q)$
{\it is the
multiplicative
subgroup of such elements $z$ in $Cl^+(\La,Q)$ that
\begin{enumerate}
\item $z\La z^{-1}=\La,$
\item $N(z):=zz'=1.$
\end{enumerate}}

\medskip
\noindent{\bf 3.4.2 Proposition.} {\it For $z \in  Spin(\La, Q)$
denote by $r_z$ the conjugation by $z$ restricted to $\La.$ Then
$r_z \in SO(\La,Q).$ The kernel of the map $r$ = \{ 1,-1 \}. In other words,
there is an exact sequence
$$
  0 \to Z/2\Z \to Spin(\La,Q) \to SO(\La,Q).
$$}

\smallskip
{\pr} See \cite{PR}.

\medskip
\noindent{\bf 3.5.} Let us extend the scalars from $\Z$ to $\Q.$ That is
we consider $\Gamma _{\Q},$ $\La _{\Q},$ $Q_{\Q},$ $\SO (\La _{\Q},Q_{\Q}),$
etc. The conditions 1, 2 in the Definition 3.4.1 above define a
$\Q$-algebraic group $Spin _{\Q}(\La _{\Q},Q_{\Q}).$ The map $r$ as in
Proposition 3.4.2 induces an exact sequence of $\Q$-algebraic groups (\cite{PR})
$$0\lto \Z /2\Z \lto Spin _{\Q}(\La _{Q},Q_{\Q}) \lto \SO(\La _{\Q},Q_{\Q})
\lto 1.$$

\medskip
\noindent{\bf 3.5.1.} Since $Spin _{\Q}(\La _{\Q},Q_{\Q})$ is a subgroup
of the multiplicative group $Cl^+(\La _{\Q}, Q_{\Q})^*,$ it acts on the
space $I_{\Q}$ (3.2). This is called the {\it spinorial} representation
of the group $Spin _{\Q}(\La _{\Q},Q_{\Q}).$ It is a direct sum of two
(nonequivalent) irreducible semispinorial representations $I_{\Q}=I_{\Q}^{ev}
\oplus I_{\Q}^{odd},$ where $I_{\Q}^{ev}$ (resp. $I_{\Q}^{odd}$) is the
direct sum of monomials $x_{j_1}\cdots x_{j_{s}}l_1\cdots l_{2n}$ with $s$ even
(resp. odd) (3.2.1 d)). We have a similar $Spin _{\Q}(\La _{\Q},Q_{\Q})$--
decomposition $I_{\Q}'=I_{\Q}^{' ev}\oplus I_{\Q}^{' odd}$ (3.3). The
isomorphism of $Cl(\La ,Q)$-modules $\beta :I\stackrel{\sim}{\to} I'$ induces an
isomorphism of $Spin _{\Q}(\La _{\Q},Q_{\Q})$--modules
$\beta _{\Q}:I_{\Q}\stackrel{\sim}{\to}I'_{\Q}.$ Thus
$$\beta _{\Q}:I_{\Q}^{ev}\stackrel{\sim}{\to} I_{\Q}^{' ev},\quad I_{\Q}^{ev}
\stackrel{\sim}{\to} I_{\Q}^{' odd}$$
or
$$\beta _{\Q}:I_{\Q}^{ev}\stackrel{\sim}{\to} I_{\Q}^{' odd},\quad I_{\Q}^{ev}
\stackrel{\sim}{\to} I_{\Q}^{' ev}.$$

\medskip
\noindent{\bf 3.5.2 Corollary.} {\it The isomorphism of $Cl(\La ,Q)$-
modules $\beta :I\stackrel{\sim}{\lto }I'$ (3.3) satisfies
$$\beta (I^{ev})=I^{' ev},\quad \beta (I^{odd})=I^{' odd}$$
or
$$\beta (I^{ev})=I^{' odd},\quad \beta (I^{odd})=I^{' ev}$$}

\medskip
\noindent{\bf 3.5.3 Definition.}
{\it We call $\beta $ {\sf even} (resp.
{\sf odd}) if $\beta (I^{ev})=I^{' ev}$ (resp. $\beta (I^{ev}=I^{' odd})$).}

\medskip
\noindent{\bf 3.5.4 Proposition.} {\it The isomorphism $\beta $ is even
(resp. odd) if $\dim_{\Q}(\Ga _{\Q}\cap M_{1\Q})$ is even (resp. odd).}

\smallskip
\pr Let $d=\dim_{\Q}(\Ga _{\Q}\cap M_{1\Q}).$
Consider the obvious $\Q$-versions of 3.2.1 and 3.3.
We may assume that
$l_1=m_1,...,l_d=m_d.$ Put ${\bf m}':=m_{d+1}\cdots m_{2n}.$ Then the
right multiplication $R_{{\bf m}'}$ by ${\bf m}'$ induces an
isomorphism $R_{{\bf m}'}:I_{\Q}\stackrel{\sim}{\lto}I'_{\Q}$ of
left $Cl(\La _{\Q},Q_{\Q})$--modules. Hence $R_{{\bf m}'}$ is a scalar
multiple of $\beta _{\Q}.$ But $R_{{\bf m}'}({\bf l})=
l_{d+1}\cdots l_{2n}{\bf m}'.$
Hence the parity of $\beta $ is equal to the parity of $d.$ This proves
the proposition.
$\Box$

\medskip
\noindent{\bf 3.6 Remark.} In the notations of 3.2 the choice of a basis
$x_1,...x_{2n}$ of $\Ga $ induces an isomorphism of abelian groups
$$I\cong \La ^{\cdot} \Ga ^*$$
(see Prop. 3.2.1 d)). A different choice of a basis of $\Ga $ may change this
isomorphism by $-1.$ Thus we get a canonical left $Cl(\La ,Q)$-module
structure on $\La ^{\cdot }\Ga ^*.$ Similarly in the notation of Corollary
3.3 the group $\La ^{\cdot}M_2$ is canonically a left $Cl(\La, Q)$-module.
Moreover, this corollary asserts that there exists a unique (up to $\pm 1$)
isomorphism of $Cl (\La ,Q)$-modules
$$\La ^{\cdot }\Ga ^*\cong \La ^{\cdot }M_2.$$

\medskip
\noindent{\bf 3.6.1.} More generally, we can introduce the following
equivalence relation $\sim $ on the set of maximal $Q_{\Q}$-isotropic
subspaces of $\La _{\Q}.$ Namely, let $L\subset \La _{\Q}$ be a maximal
isotropic subspace. Then the $\Q$-versions of 3.2.1 and 3.6 provide
a canonical $Cl(\La _{\Q},Q_{\Q})$-module structure on the exterior
algebra $\La ^{\cdot }L^*$ of the dual space $L^*.$ Given two
maximal isotropic subspaces $L_1,$ $L_2$ we get a unique (up to a scalar)
isomorphism $\beta _{\Q}:\La ^{\cdot }L_1^*\stackrel{\sim}{\lto}
\La ^{\cdot }L_2^*$ of $Cl(\La _{\Q},Q_{\Q})$-modules (by the $\Q$-version of
3.3). We say that $L_1\sim L_2$ if $\dim(L_1\cap L_2)$ is even. Then by 3.5.4
$L_1\sim L_2$ iff $\beta _{\Q}(\La ^{ev}L_1^*)=\La ^{ev}L_2^*.$ Thus
$\sim $ is indeed an equivalence relation, which divides the set of
maximal isotropic subspaces into two equivalence classes.

\medskip
\noindent{\bf 3.7.} We will be interested in the following situation. Let
$A$ be a complex torus. Put $\Ga =\Ga _{A},$ $\Ga ^* =\Ga _{\wh{A}},$
$\La _A=\Ga_{\vphantom{\wh{A}}A} \op\Ga _{\wh{A}}$ with the symmetric bilinear form $Q_A$
on $\La _A$ as in 3.1. Then by Remark 3.4 above $\La ^{\cdot }\Ga _{\wh{A}}$
is naturally a $Cl(\La _A,Q_A)$-module. Note the canonical isomorphisms
$\Ga _{\wh{A}}=H_1(\wh{A},\Z)=H^1(A,\Z).$ Thus the total
cohomology of $A$  $\La ^{\cdot}H^1(A,\Z)=H^*(A,\Z)$ is naturally a
$Cl(\La _A,Q_A)$-module.

Assume that $B$ is another complex torus and there exists an isomorphism
$$\alpha :\La _A\stackrel{\sim}{\lto}\La _B$$
which identifies the forms $Q_A$ and $Q_B.$ Let us identify
$$\La _A =\La =\La _B,\quad Q_A =Q=Q_B$$
$$Cl(\La _A,Q_A)=Cl(\La ,Q)=Cl(\La _B,Q_B)$$
by means of $\alpha .$ Then by Remark 3.6 there exists a unique
(up to $\pm 1$) isomorphism of $Cl(\La ,Q)$--modules
$$\beta :H^*(A,\Z)\stackrel{\sim}{\lto}H^*(B,\Z).$$
By Corollary 3.5.2 this isomorphism either preserves even and odd
cohomology groups or interchanges them. By Proposition 3.5.4 the parity of
$\beta $ (3.5.3) is equal to the parity of the dimension
$\dim_{\Q}(\Ga _{A,\Q}\cap \Ga _{B,\Q}).$

\sec{Derived categories and their groups of autoequivalences.}

\noindent{\bf 4.1 Categories of coherent sheaves and functors between them.}

\smallskip
\noindent{\bf 4.1.1.} Let
$X$
be an algebraic variety over an arbitrary field
$k.$
By
$\O_X$
denote the structure sheaf on
$X.$
Let
$coh(X)$
(resp. $Qcoh(X)$)
be the category of (quasi)--coherent sheaves on
$X.$
Recall that a  quasicoherent (coherent)
sheaf is
a sheaf of
$\O_X$-modules,
which locally on
$X$
has a (finite) presentation
by free
$\O_X$-modules.

It is more convenient to work with derived categories
instead of abelian categories.
Let us denote by
${\db{X}}$
the bounded derived category of
$coh(X).$

Any derived category
$\D$
has a structure of a triangulated category.
This means  that there are fixed
\begin{list}{\alph{tmp})}%
{\usecounter{tmp}}
\item a translation functor
$[1]: \D\lto\D$
that is an additive autoequivalence,
\item a class of distinguished (or exact) triangles:
$$
X\stackrel{u}{\lto}Y\stackrel{v}{\lto}Z\stackrel{w}{\lto}X[1]
$$
that satisfies a certain set of axioms (see \cite{Ve}).
\end{list}

An additive functor
$F : \D\lto \D'$
between two triangulated categories
is called exact
if it commutes with the translation functors,
and it takes every distinguished triangle in
$\D$
to
a distinguished  triangle in
$\D'.$

A natural example of an exact functor is related to a
map of algebraic varieties.
Every such map $f: X\to Y$ induces the functor
of direct image $f_{*}: Qcoh(X)\to Qcoh(Y).$
Since a category
$Qcoh(X)$ has enough injectives
there exists  the right derived functor
${\bf R}f_{*}: D^b(Qcoh(X))\lto D^b(Qcoh(Y)).$
This gives us an example of  exact functor between
derived categories.
If the map
$f$
is proper then for any $i$ the functor
${\bf R}^{i}f_{*}$ takes a coherent sheaf to a coherent one.
It is known that
${D^b(coh(X))}$
is equivalent to the full subcategory
$D^{b}(Qcoh(X))_{coh}\subset D^{b}(Qcoh(X))$
 objects of which have coherent cohomologies.
Hence any proper map $f$ induces
the exact functor
${\bf R}f_{*}: D^b(coh(X))\lto D^b(coh(Y)).$
Further, the functor $f_{*}$ has a left adjoint functor
$f^{*}: coh(Y)\lto coh(X).$
If the map $f$  has finite Tor--dimension,
then there exists
a left derived functor
${\bf L}f^{*}: D^b(coh(Y))\lto D^b(coh(X)),$
which is left adjoint for
${\bf R}f_{*}.$
This is another example of an exact functor.
(see \cite{EGAIII}, \cite{SGA6} Ex.I,II,III).

\medskip
\noindent{\bf 4.1.2.} From now on we will assume that all varieties are smooth
and projective and will denote
$D^b(coh(X))$ by
${\db{X}}.$
In this case every map
$f: X\to Y$ induces the functors
${\bf R}f_{*}: {\db{X}}\lto D^b(Y)$
and
${\bf L}f^{*}: D^b(Y)\lto D^b(X).$
Also each complex
${\F}\in {\db{X}}$ induces an exact functor
$\stackrel{\bf L}{\ot}\F : {\db{X}}\lto{\db{X}}.$
Using these functors one can introduce a larger
class of exact functors.

Let $X$ and $Y$ be  smooth projective varieties of dimension
$n$ and $m$ respectively.
Consider the projections
$$
\begin{array}{ccccc}
X&\stackrel{p}{\longleftarrow}&X\times Y&\stackrel{q}{\lto}&Y\\
\end{array}
$$
Every object
${\E}\in {\db{X\times Y}}$
defines an exact functor
$\Phi_{\E}:{\db{X}}\lto{\db{Y}}$
by the following formula:
\begin{equation}\label{dfun}
\Phi_{\E}(\cdot):={\bf R}q_{*}({\E}\stackrel{\bf L}{\otimes}
p^*(\cdot)).
\end{equation}
Obviously, the same object defines another functor
$\Psi_{\E}:{\db{Y}}\lto{\db{X}}$
by the similar formula
$$
\Psi_{\E}(\cdot):={\bf R}p_{*}({\E}\stackrel{\bf L}{\otimes}
q^*(\cdot)).
$$
The functor
$\Phi_{\E}$
has left and right adjoint functors
$\Phi_{\E}^*$
and
$\Phi_{\E}^{!}$
respectively,
defined by the rule:
\begin{equation}
\begin{array}{lll}\label{adj}
\Phi_{\E}^{*}&\cong&\Psi_{\E^{\vee}\ot q^{*}K_{Y}[m]}\\
\Phi_{\E}^{!}&\cong&\Psi_{\E^{\vee}\ot p^{*}K_{X}[n]}
\end{array}
\end{equation}
Here $K_X$ and $K_Y$ are the canonical sheaves
on $X$ and $Y$ respectively, and
$\E^{\vee}$
is notation for
${\bf R}\underline{{\H}om}(\E, \O_{X\ts Y}).$

\medskip
\noindent{\bf 4.1.2.1 Example.}
Let $f: X\to Y$ be a morphism of varieties.
Denote by $\Ga(f)$ a subvariety of $X\ts Y$ that is the graph of
$f.$ It is clear that the functor
${\bf R}f_* : {\db{X}}\lto{\db{Y}}$ is isomorphic to
$\Phi_{\O_{\Ga(f)}},$ where $\O_{\Ga(f)}$ is the structure sheaf of
the graph $\Ga(f).$ In the same way the functor
${\bf L}f^*:{\db{Y}}\lto{\db{X}}$ is isomorphic to
$\Psi_{\O_{\Ga(f)}}.$


\medskip
\noindent{\bf 4.1.2.2 Example.}\label{tm}
Let
$\delta: \Delta\hookrightarrow X\ts X$
be the embedding  of the diagonal.
First, since the $\Delta$ is the graph of the identity morphism
the identity functor from the
derived category
${\db{X}}$ to itself is represented by the structure sheaf
$\O_{\Delta}.$
Moreover, for any object $\F\in {\db{X}}$ the functor
$\stackrel{\bf L}{\ot}\F : {\db{X}}\lto{\db{X}}$
is isomorphic to
$\Phi_{\delta_{*}\F}$

\medskip
\noindent{\bf 4.1.2.3.}
Thus there is a reasonably large class of exact functors
between derived categories of smooth projective varieties
that consists of functors having the form
$\Phi_{\E}$ for some complex
$\E$ on the product.
This class is closed under composition of functors.
Actually, let
$X, Y, Z$ be three (smooth projective) varieties
and let
$$
\Phi_{I}:{\db{X}}\lto{\db{Y}},\qquad
\Phi_{J}:{\db{Y}}\lto{\db{Z}}
$$
be two functors, where
$I$ and $J$ are objects of
${\db{X\ts Y}}$ and
${\db{Y\ts Z}}$
respectively.
Denote by
$p_{XY}, p_{YZ}, p_{XZ}$ the projections of
$X\ts Y\ts Z$ on the corresponding pair of factors.

\medskip
\noindent{\bf 4.1.2.4 Lemma.}{\rm \cite{Mu1}}\label{comp}
 {\it In above notations, the composition
$\Phi_{J}\circ\Phi_{I}$ is isomorphic to
$\Phi_{K},$ where
$K\in{\db{X\ts Z}}$ is given by the formula:
$$
K\cong {\bf R}p_{XZ*}(p_{YZ}^*(J)\ot p_{XY}^*(I))
$$}

\medskip
Presumably, the  class of exact functors
described above embraces all exact functors between
bounded derived categories of coherent sheaves on
smooth projective varieties.
We do not know if it is true or not.
However it is definitely true for exact equivalences.

\medskip
\noindent{\bf 4.1.2.5 Theorem.}{\rm (\cite{Or1})}\label{rep} {\it
Let $X$ and $Y$ be smooth projective varieties.
Suppose
$F: {\db{X}}\stackrel{\sim}{\lto}{\db{Y}}$ is an exact equivalence.
Then there exists a unique (up to an isomorphism) object
$\E\in{\db{X\ts Y}}$ such that the functors
$F$ and $\Phi_{\E}$ are isomorphic.}

\medskip

\noindent{\bf 4.1.2.6 Definition.}
{\it The group of isomorphism classes of exact equivalences
$F: {\db{X}}\stackrel{\sim}{\lto} {\db{X}}$
is called the {\sf group of autoequivalences} of
${\db{X}}$ and is denoted by
$Auteq({\db{X}}).$}

\medskip
\noindent{\bf 4.1.3.} Assume for simplicity that $k=\C.$
Analogously to derived categories and functors between them
one can consider cohomologies and maps between them.
The latter is much simpler.

For every element
$\xi\in H^{*}(X\ts Y, \Q)$
let us define linear maps
$$
{v}_{\xi}: H^{*}(X, \Q)\lto H^{*}(Y, \Q),\qquad
{w}_{\xi}: H^{*}(Y, \Q)\lto H^{*}(X, \Q)
$$
by the formula :
\begin{equation}\label{dcor}
v_{\xi}(\cdot)=q_{*}({\xi}\cup p^*(\cdot)),\qquad
w_{\xi}(\cdot)=p_{*}({\xi}\cup q^{*}(\cdot)).
\end{equation}
It follows from the K\"unneth formula that any linear map between
cohomologies have
this form for some
$\xi.$

The composition formula for cohomological correspondence
is well-known and
is similar to the composition of functors.
Let as above
$X, Y, Z$ be three varieties
and let
$\xi$ and $\eta$ be  elements of
$H^{*}(X\ts Y)$ and
$H^{*}(Y\ts Z)$
respectively.

\medskip
\noindent{\bf 4.1.3.1 Lemma.}\label{compcor} {\it
 The composition
$v_{\eta}\circ v_{\xi}$ is isomorphic to
$v_{\zeta},$ where
$\zeta\in H^{*}(X\ts Z)$ is given by the formula:
$$
\zeta = p_{XZ*}(p_{YZ}^*(\eta)\cup p_{XY}^*(\xi)).
$$}

\medskip
\noindent{\bf 4.1.4.}
There is a natural correspondence that attaches to a functor
$\Phi_{\E}:{\db{X}}\lto{\db{Y}}$ a linear map
of vector spaces
$\phi_{\E}: H^*(X, \Q)\lto H^*(Y, \Q).$
To describe it note that each exact functor
$F:{\db{X}}\lto {\db{Y}}$
induces a homomorphism
$\overline{F}: K(X)\lto K(Y),$
where $K(X)$ and $K(Y)$ are the Grothendieck groups
 of the categories
${\db{X}}$ and ${\db{Y}}.$
On the other side,
there is a map
$$
ch: K(X)\lto H^*(X, \Q)
$$
that is called
the Chern character. The map $ch$ is  a ring homomorphism,
if to recall that $K(X)$ and $H^*(X, \Q)$ are endowed
with multiplications induced by $\ot$--product  and
$\cup$--product respectively.

\medskip
\noindent{\bf 4.1.4.1 Definition.}
{\it For every object ${\E}\in{\db{X\ts Y}}$ put
\begin{equation}\label{ccor}
\begin{array}{c}
\phi_{\E}(\cdot):=v_{ch(\E)\cup p^*(td_X)}(\cdot)=
q_{*} (ch(\E)\cup p^{*}(td_X)\cup p^*(\cdot)),\\
\psi_{\E}(\cdot):=w_{ch(\E)\cup q^*(td_Y)}(\cdot)=
p_{*} (ch(\E)\cup q^{*}(td_Y)\cup q^*(\cdot)),
\end{array}
\end{equation}
where
$td_X$ and $td_Y$
are the Todd classes of
$X$ and $Y.$}

\medskip
\noindent{\bf 4.1.4.2 Lemma.}\label{ch} {\it
 For given
${\E}\in {\db{X\ts Y}}$
the following diagram
$$
\begin{array}{ccc}
{K(X)}&\stackrel{\overline{\Phi}_{\E}}{\lto}&{K(Y)}\\
\llap{\ss{ch}}\da&&\da\rlap{\ss{ch}}\\
H^* (X, \Q)&\stackrel{\phi_{\E}}{\lto}&H^* (Y, \Q)
\end{array}
$$
is commutative.}

\medskip
\noindent{\bf 4.1.4.3 Lemma.}\label{repr} {\it
 Assume that a functor
$\Phi_{K}:{\db{X}}\lto{\db{Z}}$
is isomorphic to
$\Phi_{J}\circ\Phi_{I}$
for some
$$
\Phi_{I}:{\db{X}}\lto{\db{Y}}, \qquad
\Phi_{J}:{\db{Y}}\lto{\db{Z}}.
$$
Then
$\phi_{K}=\phi_{J}\circ\phi_{I}.$
Consequently, the correspondence
$\Phi_{\E}\mapsto \phi_{\E}$
induces a representation
$$
\rho_{X}: Auteq({\db{X}})\lto \GL(H^*(X,\Q)).
$$}
\smallskip
Both lemmas immediately follows from the
Riemann-Roch-Grothendieck theorem.

\medskip
\noindent{\bf 4.1.4.4 Example.}
Let $f: X\to Y$ be a morphism. The functors
${\bf L}f^*\cong \Psi_{\O_{\Ga(f)}}$ and
${\bf R}f_{*}\cong \Phi_{\O_{\Ga (f)}}$ give the maps
$\psi_{\O_{\Ga (f)}}$ and $\phi_{\O_{\Ga (f)}}$
between
cohomologies. Also the morphism $f$ induces the map:
$
f^*:H^{i}(Y, \Q)\lto H^{i}(X, \Q).
$
 It can be  checked that
$$
\psi_{\O_{\Ga (f)}}=f^{*}.
$$
On the other side, there is a map of homologies:
$
f_{*}: H_{i}(X,\Q)\lto H_{i}(Y, \Q)
.$
Using the Poincare isomorphism
$H^{2n-i}(X, \Q)\cong H_{i}(X,\Q),$ the map $f_{*}$ can be considered
as a map between cohomologies:
$
f_{*}: H^{2n-i}(X, \Q)\lto H^{2m-i}(Y, \Q).
$
In this case
$\phi_{\O_{\Ga (f)}}$
does not coincide with
$f_{*},$
but they are connected by
the following relation:
$$
\phi_{\O_{\Ga (f)}}(x)\cup td_Y = f_{*}(x\cup td_X)
$$
In particular, if $td_X=td_Y=1,$ then these maps coincide.

\medskip
\noindent{\bf 4.1.4.5 Example.}
Let as in Example 4.1.2.2
$\delta:\Delta\hookrightarrow X\ts X$
be an embedding of the diagonal.
For any
$\F\in{\db{X}}$
a functor
$\stackrel{\bf L}\ot\F$ is isomorphic to
$\Phi_{\delta_*\F}.$
Since $ch$ is a ring homomorphism
Lemma 4.1.4.2 implies that
$$
\phi_{\delta_{*}\F}(\cdot)=(\cdot)\cup ch(\F).
$$

\medskip
\noindent{\bf 4.1.5.}
For any variety $X$ the group
$Auteq({\db{X}})$ contains a subgroup
$G(X)=(\Z\ts \Pic(X))\rtimes Aut(X)$
that is the semi-direct product of the normal subroup
$(\Z\ts\Pic(X))$ and $Aut(X)$ naturally acting on
the former.
Under the inclusion
$G(X)\subset Auteq({\db{X}})$
the generator of
$\Z$
goes to the translation functor
$[1],$
a line bundle
$\L\in\Pic(X)$
maps to the functor
$\ot \L=\Phi_{\delta_{*}\L}$
and an automorphism
$f:X\to X$ induces the autoequivalence
${\bf R}f_{*}\cong \Phi_{\O_{\Ga(f)}}.$

\medskip
\noindent{\bf 4.1.5.1 Lemma.}\label{simpeq} {\it
In the above notation
$\rho_{X}([1])=-Id_{H^*(X,\Q)}$
and
$\rho_{X}(\ot\L)=\cup ch(\L).$}

\smallskip
The proof follows from the Example 4.1.4.5 above and the equality
$ch(E[1])=-ch(E).$

\noindent{\bf 4.1.6.}
Now let us show that any equivalence
$\Phi_{\E}:{\db{X}}\stackrel{\sim}{\lto}{\db{Y}}$
defines a functor between
derived categories of coherent sheaves on
$X\ts X$ and $Y\ts Y,$ which will be called adjoined.
Consider two functors
$$
\Phi_{\E}:{\db{X}}{\lto}{\db{Y}},\qquad
\Psi_{\F}:{\db{Y}}{\lto}{\db{X}},
$$
 where
$\E$ and $\F$ are objects of
${\db{X\ts Y}}.$
By $\F\bt\E$ denote the object
$p_{13}^{*}(\F)\ot p_{24}^{*}(\E)$
from
${\db{X\ts X\ts Y\ts Y}},$
which is called the external tensor product of
$\E$ and $\F.$ It defines the functor
$$
\Phi_{\F\bt\E}:{\db{X\ts X}}\lto{\db{Y\ts Y}}.
$$

Let us take an object
${\G}\in {\db{X\ts X}}$ and denote
$\Phi_{\F\bt \E}(\G)\in{\db{Y\ts Y}}$ by
$\H$ for short.
These two objects define the functors
$$
\Phi_{\G}:{\db{X}}\lto{\db{X}},\qquad
\Phi_{\H}:{\db{Y}}\lto{\db{Y}}.
$$

\medskip
\noindent{\bf 4.1.6.1 Lemma.}\label{adcom} {\it
In above notation,  there is an isomorphism of functors
$\Phi_{\H}\cong \Phi_{\E}\Phi_{\G}\Psi_{\F}.$}

\smallskip
This follows from Lemma 4.1.2.4. $\Box$

\medskip
\noindent{\bf 4.1.6.2 Lemma.}\label{pr} {\it
If
$\Phi_{\E}$ and
$\Psi_{\F}$ are equivalences, then the functor
$\Phi_{\F\bt\E}$
is an equivalence too.}

\smallskip

The proof is straightforward (see for example \cite{Or2}).

\medskip
\noindent{\bf 4.1.6.3 Definition.}\label{adjfun}
{\it Assume that $\Phi_{\E}:{\db{A}}\lto {\db{B}}$ is an equivalence and
$\Psi_{\F}\cong\Phi_{\E}^{-1}.$ In this case
the functor
$\Phi_{\F\bt\E}$ and the map
$\phi_{\F\bt\E}$ will be  denoted  by $Ad_{\E}$ and by
$ad_{\E}$ respectively.}

\smallskip

Thus any exact equivalence
$\Phi_{\E}:{\db{X}}\stackrel{\sim}{\lto}{\db{Y}}$ defines an exact equivalence
$Ad_{\E}:{\db{X\ts X}}\stackrel{\sim}{\lto}{\db{Y\ts Y}}$
such that , by Lemma 4.1.6.1, there is an isomorphism:
\begin{equation}\label{ad}
\Phi_{Ad_{\E}(\G)}\cong \Phi_{\E}\Phi_{\G}\Phi_{\E}^{-1}.
\end{equation}

\vspace{1.3cm}
\noindent{\bf 4.2 Category of coherent sheaves on an abelian variety.}

\medskip
\noindent{\bf 4.2.1.}
Let
$A$
be an abelian variety of dimension $n$ over $\C.$
Denote by
$m:A\ts A\to A$
the multiplication morphism.
For any point
$a\in A$ by
$T_a$ denote the translation automorphism
$m(\cdot, a):A\to A.$

Let
$\wh{A}$
be dual abelian variety (1.2). It is canonically isomorphic to
$\Pic^0(A)$ (1.3).
It is well-known that there is a unique line bundle
$P$
on the product
$A\ts\wh{A}$
such that for any point
$\alpha\in \wh{A}$
the restriction
$P_{\alpha}$
on $A\ts \{\alpha\}$
represents an element of
$\Pic^0 A,$ corresponding to
$\alpha$
 and , in addition, the restriction
$P\Big|_{\{0\}\ts \wh{A}}$ should be trivial.
Such $P$ is called Poincare line bundle (1.4).
As in 1.4, 1.5 we identify $A$ with $\wh{\wh{A}}$ using
the Poincare line bundles on $A\times \wh{A}$ and
$\wh{A}\times \wh{\wh{A}}.$

\medskip
\noindent{\bf 4.2.2.}
The Poincare line bundle gives an example of an exact equivalence
between derived categories of coherent sheaves of two non-isomorphic varieties.
Let us consider the projections
$$
\begin{array}{ccccc}
A&\stackrel{p}{\longleftarrow}&A\times \wh{A}&\stackrel{q}{\lto}&\wh{A}\\
\end{array}
$$
and the functor
$\Phi_{P}:{\db{A}}\lto{\db{\wh{A}}},$
defined as in 4.1.2, i.e.
$\Phi_{P}(\cdot)={\bf R}q_{*}(P\ot p^{*}(\cdot)).$

The following proposition  was proved by Mukai.

\medskip
\noindent{\bf 4.2.2.1 Proposition.}{\rm (\cite{Mu1})} {\it
Let $P$ be an Poincare line bundle on
$A\ts\wh{A}.$ Then the functor
$\Phi_{P}:{\db{A}}\lto{\db{\wh{A}}}$
is an exact equivalence, and there is an isomorphism of functors:
$$
\Psi_{P}\circ\Phi_{P}\cong (-1_A)^*[n],
$$
where $(-1_A)$ is the inverse map on the group $A.$}

\medskip

\noindent{\bf 4.2.3.}
Let
$x_1,...,x_{2n}$
be a basis of
$H^1(A,\Z ).$
Let
$l_1,...,l_{2n}$
be the dual basis of
$H^{1}(\wh{A}, \Z).$
It is clear that
$p^{*}(x_1),...,p^{*}(x_{2n}), q^*(l_1),...,q^{*}(l_{2n})$
is a basis of
$H^{1}(A\ts\wh{A}, \Z).$

\medskip
\noindent{\bf 4.2.3.1 Lemma.}{\rm (\cite{Mum1})}\label{c1} {\it
The first Chern class of the Poincare line bundle
$P$ on
$A\ts \wh{A}$  satisfies the equality:
\begin{equation}\label{c1f}
c_1(P)=\sum_{i=1}^{2n} p^{*}(x_i)\cup q^{*}(l_i).
\end{equation}}

Now consider the map
$\phi_{P}: H^*(A, \Q)\lto H^*(\wh{A}, \Q),$ given by formula (\ref{ccor}).

\medskip
\noindent{\bf 4.2.3.2 Lemma.} {\it
For any $k$ the map
$\phi_{P}$ sends $k$-th cohomology to $(2n-k)$-th cohomology
and induces isomorphisms
$$
H^{k}(A, \Z)\stackrel{\sim}{\lto}H^{2n-k}(\wh{A}, \Z),
$$
 under which a monomial
$x_{i_1}\cup\cdots\cup x_{i_k}$ goes to
$(-1)^{\varepsilon}l_{j_1}\cup\cdots\cup l_{j_{2n-k}},$ where
$(j_1<\cdots< j_{2n-k})$ is the complement of
$(i_1<\cdots< i_{k})$ in the set $(1,...,2n),$ and
$\varepsilon=\sum_{t}j_t.$}

\pr
Take $\alpha\in H^k(A, \Z).$ Since
$ch(P)=\exp(c_1(P))$ and
$c_1(P)$
has the form (\ref{c1f}),
we obtain
$$
\phi_{P}(\alpha)=q_{*}(ch(P)\cup p^*(\alpha))=
\frac{1}{(2n-k)!}q_{*}(c_1(P)^{2n-k}\cup p^{*}(\alpha))
$$
This implies that
$\phi_{P}$ sends $k$-th cohomology to $(2n-k)$-th.
Moreover, it is easy to see that
$$
\frac{1}{s!}c_1(P)^{s}=\sum_{j_{1}<\cdots<j_s}p^*(x_{j_s}\cup\cdots\cup x_{j_1})
\cup
q^{*}(l_{j_1}\cup\cdots\cup l_{j_s})
$$
Substituting this expression in previous formula, we get
$$
\begin{array}{ll}
\phi_{P}(x_{i_1}\cup\cdots\cup x_{i_k})&=
q_{*}p^*(x_{i_1}\cup\cdots\cup x_{i_k}
\cup x_{j_{2n-k}}\cup\cdots\cup
x_{j_1})\cup l_{j_1}\cup\cdots\cup l_{j_{2n-k}}=\\
&(-1)^{\varepsilon}q_{*}p^*(x_{1}\cup\cdots\cup x_{2n})\cup
l_{j_1}\cup\cdots\cup l_{j_{2n-k}}=(-1)^{\varepsilon}
l_{j_1}\cup\cdots\cup l_{j_{2n-k}}
\end{array}
$$
where
$(j_1<\cdots<j_{2n-k})$
is the complement of
$(i_1<\cdots<i_k),$
and
$\varepsilon=\sum_{t}j_t.$
The monomials
$x_{i_1}\cup\cdots\cup x_{i_k}$
and
$l_{j_1}\cup\cdots\cup l_{j_{2n-k}}$
form basises of
$H^{k}(A, \Z)$
and
$H^{2n-k}(\wh{A}, \Z)$ respectively.
Hence the map $\phi_{P}$ gives an isomorphism between these lattices.
$\Box$
\bigskip

\noindent{\bf 4.2.4.}
For any point $(a, \alpha)\in A\ts \wh{A}$
one can introduces a functor   from
${\db{A}}$ to itself, defined by the rule:
\begin{equation}\label{fker}
\Phi_{(a,\alpha)}(\cdot):= T_a^* (\cdot)\ot P_{\alpha}.
\end{equation}
A functor
$\Phi_{(a,\alpha)}$ is represented by a sheaf
$S_{(a,\alpha)}=\O_{\Ga_a}\ot p_2^*(P_{\alpha})$ on
$A\ts A,$ where
$\Ga_a$ is the graph of the automorphism
$T_a : A\lto A.$
It is clear that the functor $\Phi_{(a,\alpha)}$ is
an equivalence. Since an automorphism $T_a$ acts identically
on the cohomology group and
$ch(P_{\alpha})=1$ the functors
$\Phi_{(a,\alpha)}$
belongs to kernel of the homomorphism
$\rho$ (4.1, 4.3).

The set of the functors
$\Phi_{(a,\alpha)},$ parametrized by
$A\ts\wh{A},$
can be  unified in one functor from
${\db{A\ts\wh{A}}}$ to ${\db{A\ts A}}$
that takes a skyscraper
$\O_{(a,\alpha)}$
to
$S_{(a,\alpha)}.$
(Note that this condition does not define the functor by uniquely,
because it can be changed by  tenzoring with any line bundle pulling back
from
$A\ts \wh{A}.$)
Let us define  such a functor
$\Phi_{S_A} : {\db{A\ts\wh{A}}}\lto {\db{A\ts A}}$
by the representing object
$$
{S}_A=(m\cdot p_{13}, p_4)^* \O_{\Delta}\ot p_{23}^* P_A
$$
on the product
$(A\ts\wh{A})\ts (A\ts A).$
Here
$(m\cdot p_{13}, p_4)$
is a morphism onto
$A\ts A$
that takes
$(a_1,\alpha,a_3,a_4)$ to
$(m(a_1,a_3),a_4).$
A direct check shows that
$\Phi_{S_A}$ takes the skyscraper sheaf
$\O_{(a,\alpha)}$ to
$S_{(a,\alpha)}.$

The functor
$\Phi_{S_A}$
can be consider as a composition of two
functors.
Denote by
$\cP_A=p^*_{14}\O_{\Delta}\ot p^{*}_{23}P
\in{\db{(A\ts\wh{A})\ts(A\ts A)}}.$
By
$\mu_A: A\ts A \lto A\ts A$
denote a morphism that takes
$(a_1, a_2)$ to $(a_1, m(a_1, a_2)).$
Consider two functors
$$
\Phi_{\cP_{A}}: D^b(A\ts\wh{A})\lto D^b(A\ts A),\qquad
{\bf R}\mu_{A*} :{\db{A\ts A}}\lto{\db{A\ts A}}.
$$

\medskip
\noindent{\bf 4.2.4.1 Lemma.}\label{cP} {\it
The functor
$\Phi_{S_A}$
is isomorphic to the composition
${\bf R}\mu_{A*} \circ\Phi_{\cP_A}.$}

\smallskip
The proof is  omitted. Also this statement can be used as another definition
of the functor
$\Phi_{S_A}.$

\medskip
\noindent{\bf 4.2.4.2 Lemma.} {\it
The functor
$\Phi_{S_A}$
is an equivalence.}

\smallskip
\pr
By Lemma 4.1.6.2 the functor
$\Phi_{\cP_A}$
is an equivalence.
Since
$\mu_A$
is an automorphism of
$A\ts A$
the functor
${\bf R}\mu_{A*}$
is an equivalence too. This implies that
$\Phi_{S_A}$ is also an equivalence.

\medskip
\noindent{\bf 4.2.5.} Let $A$ and $B$ be two abelian varieties. Fix
an equivalence $\Phi_{\mathcal{E}}:
{\db{A}}\stackrel{\sim}{\lto}{\db{B}}.$ Consider the functor
$\Phi_{S_B}^{-1} Ad_{\mathcal{E}} \Phi_{S_A}$ from
${\db{A\ts\wh{A}}}$ to ${\db{B\ts\wh{B}}}$ (4.1, 6.3). By
${\J}({\E})$ denote an object that represents this functor. Thus
there is the commutative diagram:
\begin{equation}\label{diag}
\begin{array}{ccc}
{\db{A\ts\wh{A}}}&\stackrel{\Phi_{S_A}}{\lto}&{\db{A\ts A}}\\
\llap{\ss{\Phi_{{\mathcal{J}}({\mathcal{E}})}}}\da&&\da\rlap{\ss{Ad_{\mathcal{E}}}}\\
{\db{B\ts\wh{B}}}&\stackrel{\Phi_{S_B}}{\lto}&{\db{B\ts B}}
\end{array}
\end{equation}
Let us describe the object
$\J(\E).$

\medskip

\noindent{\bf 4.2.5.1 Theorem.}{\rm(\cite{Or2})}\label{obj} {\it
There exists a group isomorphism
$f_{\E}: A\ts \wh{A} \lto B\ts \wh{B}$
and a line bundle
$L_{\E}$ on
$A\ts \wh{A}$
such that the object
${\J}({\E})$
is isomorphic to
$i_* (L_{\E}),$ where
$i$
is the embedding of
$A \ts \wh{A}$
in
$(A\ts\wh{A})\ts (B\ts\wh{B})$
as the
graph
$\Ga(f_{\E})$
of the isomorphism
$f_{\E}.$}

\smallskip

In particular, it immediately follows from the theorem that if two abelian
varieties
$A$
and
$B$
have equivalent derived categories of coherent sheaves,
then the varieties
$A\ts\wh{A}$ and
$B\ts \wh{B}$ are isomorphic.
\medskip

\noindent{\bf 4.2.5.2 Corollary.} {\it
An isomorphism
$f_{\E}$ takes
a point
$(a,\alpha)\in A\ts\wh{A}$
to a point
$(b,\beta)\in B\ts\wh{B}$
iff
the equivalences
$$
\Phi_{(a,\alpha)}:\db{A}\stackrel{\sim}{\lto}\db{A},\qquad
\Phi_{(b,\beta)}:\db{B}\stackrel{\sim}{\lto}\db{B},
$$
defined by formula (\ref{fker}), are connected by the relation
$$
\Phi_{(b,\beta)}\cong\Phi_{\E}\Phi_{(a,\alpha)}\Phi_{\E}^{-1}.
$$}

\pr
By the theorem, $\Phi_{\J(\E)}$ takes a skyscraper
$\O_{(a,\alpha)}$ to
$\O_{(b,\beta)},$ where
$(b,\beta)=f_{\E}(a,\alpha).$
By construction of the functor
$\Phi_{S_A}$ it takes a skyscraper
$\O_{(a,\alpha)}$  to
$S_{(a,\alpha)},$ where the object
$S_{(a,\alpha)}$ represents the functor
$\Phi_{(a,\alpha)}.$
Thus diagram (\ref{diag}) implies that the morphism
$f_{\E}$ takes
$(a,\alpha)$ to
$(b,\beta)$
iff
$S_{(b,\beta)}\cong Ad_{\E}(S_{(a,\alpha)}).$
Now the formula (\ref{ad}) implies that
$\Phi_{(b,\beta)}\cong \Phi_{\E}\Phi_{(a,\alpha)}\Phi_{\E}^{-1}.$
$\Box$

\medskip
\noindent{\bf 4.2.5.3 Proposition.}{\rm (\cite{Or2})} {\it Let
$B\cong A,$ the correspondence $\Phi_{\mathcal{E}} \mapsto
f_{\mathcal{E}}$ induces  a group homomorphism
$$
\gamma_A:Autoeq({\db{A}})\lto Aut(A\ts\wh{A}).
$$}

\pr
Let $\E_1, \E_2, \E_3$ be objects of ${\db{A\ts A}}$
such that $\Phi_{\E_i}$ are equivalences and
$\Phi_{\E_3}\cong \Phi_{\E_2}\circ\Phi_{\E_1}.$
This implies that
$Ad_{\E_3}\cong Ad_{\E_2}\circ Ad_{\E_1}.$
Further we have the sequence of isomorphisms:
$$
\Phi_{\J(\E_2)}\circ \Phi_{\J(\E_1)}\cong
\Phi_{S_A}^{-1}Ad_{\E_2}\Phi_{S_A}\Phi_{S_A}^{-1}Ad_{\E_1}\Phi_{S_A}
\cong
\Phi_{S_A}^{-1}Ad_{\E_3}\Phi_{S_A}
\cong
\Phi_{\J(\E_3)}
$$
By theorem 4.2.5.1 the object $\J(\E)$ is a line bundle over
the graph  of some automorphism $f_{\E}.$ Thus we obtain that
$f_{\E_3}= f_{\E_2}\cdot f_{\E_1}.$
$\Box$

\medskip
\noindent{\bf 4.2.5.4 Proposition.}{\rm (\cite{Or2})} {\it
The kernel of the homomorphism
$\gamma_A$
coincides with
$\Z\ts A\ts\wh{A},$
consisting of all equivalences
$\Phi_{(a,\alpha)}[i]= T_a^{*}(\cdot)\ot P_{\alpha}[i].$}

\vspace{0.5cm}
\noindent{\bf 4.3 The Mukai-Polishchuk group $\U (A)$ and the spinorial group
$Spin(A).$}
\medskip

In his recent preprint S.Mukai \cite{Muk2} was  interested in the group
of autoequivalences
$Auteq (D^{b}(A))$
of the bounded derived category
$D^{b}(A)$
of coherent sheaves on an abelian variety
$A$ and
he defined a certain related discrete group (which he called
$U(A\times \wh A)$).
Independently A.Polishchuk considered the same group,
which he called
${\rm S}{\rm L}_2 (A),$
and proved that
$\SL_2(A)$
acts naturally on
${\db{A}}$
"up to the shift functor"(see \cite{Po2}).
Let us recall the definition of this group for complex tori. We call
this group $\U (A).$

\medskip

\noindent{\bf 4.3.1 Definition.}
{\it Let $A$ be a complex torus.
Put
$$
\U (A):=
\left\{
g=\left(
\begin{array}{cc}
a&b\\
c&d
\end{array}
\right)
\in
\left(
\begin{array}{ll}
{\End}(A)&{\Hom}({\wh{A}}, A)\\
{\Hom}(A, \wh{A})&{\End}(\wh{A})
\end{array}
\right)
\left|
\quad
g^{-1}=
\left(
\begin{array}{rr}
\wh{d}&-\wh{b}\\
-\wh{c}&\wh{a}
\end{array}
\right)
\right\}
\right.
$$
and call it the {\sf  Mukai-Polishchuk  group}.}

\medskip

We may consider $\U (A)$ as a subgroup of $\GL(\Ga _A\oplus\Ga _{\wh A}).$
Consider the canonical symmetric bilinear form $Q$ on $\La :=
\Ga _A\oplus \Ga _{\wh A}$ as in 3.1 above.
Here is another description of the group $\U (A).$

\medskip
\noindent{\bf 4.3.2 Proposition.} {\it
There are equalities
$$
\U (A)=O(\La, Q)\cap Aut(A\ts\wh{A})=
SO(\La, Q)\cap Aut(A\ts\wh{A})
$$
of subgroups of $\GL({\La}).$}

\smallskip
\pr
The second equality follows from the fact that elements in $Aut(A\ts\wh A)$
preserve the complex structure on $V_A\oplus V_{\wh A},$ hence have a positive
determinant.

To prove the first equality it suffices to show that the group
$$
T:=
\left\{
g=\left(
\begin{array}{cc}
a&b\\
c&d
\end{array}
\right)
\in
\GL(\La)
\left|
\quad
g^{-1}=
\left(
\begin{array}{rr}
\wh{d}&-\wh{b}\\
-\wh{c}&\wh{a}
\end{array}
\right)
\right\}
\right.
$$
coincides with the orthogonal group $O(\La,Q).$ Choose a basis
$l_1,...,l_{2n}$ of $\Ga _A$ and the dual basis $x_1,...,x_{2n}$
of $\Ga _{\wh A}.$ Consider the basis $l_1,...,l_{2n},x_1,...,x_{2n}$
of the lattice $\La .$ Then the group $O(\La ,Q)$ consists of matrices
$$\left(
\begin{array}{cc}
\alpha & \beta \\
\gamma & \delta
\end{array}\right),$$
such that
$$\left(
\begin{array}{cc}
\alpha & \beta \\
\gamma & \delta
\end{array}\right)
 \left(
\begin{array}{cc}
0 & 1 \\
1 & 0
\end{array}\right)
\left(
\begin{array}{cc}
\alpha^t & \gamma ^t \\
\beta^t  & \delta ^t
\end{array}\right)
=
\left(
\begin{array}{cc}
0 & 1 \\
1 & 0
\end{array}\right).$$
The equality $T=O(\La,Q)$ now follows from Remark 1.6.
$\Box$

\medskip
\noindent{\bf 4.3.3 Proposition.} {\rm (\cite{Or2})} {\it  The image of the
homomorphism
$\gamma_A:Autoeq ({\db{A}})\lto Aut(A\ts\wh{A})$ coincides
with the Mukai-Polishchuk group
$\U (A).$
Moreover  there is an exact sequence of groups
$$
0\lto \Z\ts A\ts\wh{A}\lto Auteq({\db{A}})\lto \U (A)\lto 1.
$$}

%

\smallskip

This was essentially conjectured by Polishchuk in \cite{Po2}.

By Lemma 4.1.4.3 the correspondence
$\Phi_{\E}\mapsto \phi_{\E}$
induces a homomorphism
$$\rho_{A}: Auteq({\db{A}})\lto \GL(H^{*}(A, \Q)).$$
As a matter of fact, this representation preserves the integral
cohomology lattice.

\medskip
\noindent{\bf 4.3.4 Proposition.}\label{Iac} {\it
For any equivalence
$\Phi_{\E}: {\db{A}}\stackrel{\sim}{\lto}{\db{B}}$
the linear map
$\phi_{\E}$
preserves the integral cohomology, i.e.
$$
\phi_{\E}(H^{*}(A, \Z))=H^{*}(B, \Z).
$$}

\pr By definition of $\phi_{\E},$ one needs to show that
$ch({\mathcal{E}})$ belongs to integral cohomology $H^*(A\ts B,
\Z).$ This is equivalent to checking that
$ch({\mathcal{E}}^{\vee}\boxtimes{\mathcal{E}})$ belongs to the
integral cohomology $H^*((A\ts A)\ts (B\ts B), \Z).$ This object
represents the functor $Ad_{\E}.$ Since $Ad_{\mathcal{E}}\cong
\Phi_{S_B} \Phi_{{\mathcal{J}}({\mathcal{E}})} \Phi_{S_A}^{-1}$ it
is sufficient to show that $ch({S}_A),$ $ch({S}^{\vee}_B)$ and
$ch({\J}({\E}))$ are integral.  But , by construction, all these
objects are line bundles on abelian subvarieties. Thus, the
following lemma implies the proposition. $\Box$

\medskip
\noindent{\bf 4.3.4.1 Lemma.} {\it
For any line bundle $L$ on an abelian variety $A$
the Chern character
$ch(L)$
is integral, i.e belongs to
$H^*(A, \Z).$}

\smallskip
\pr
Denote by
$c_1(L)$
the first Chern class of
$L.$
The Chern character
$ch(L)$
is equal to
$exp(c_1(L)).$
Thus we must show that
$\frac{1}{k!}c_1(L)^k$
belong to integral cohomology for any $k.$
If $x_1,..., x_{2n}$ be a basis of
$H^{1}(A, \Z),$ then the first Chern class as an element of
$H^2(A, \Z)$ can be written as:
$$
c_1(L)=\sum_{i<j} d_{ij}\cdot (x_{i}\cup x_{j}),\qquad d_{ij}\in \Z
$$
Taking $k$-power, we obtain
$$
{c_1(L)}^{k}=k! \sum_{(i_1 < j_1);...;(i_k < j_k)\atop (i_1,j_1,...,i_k,j_k)
\subset
(1,...,2n)}  d_{i_1 j_1}\cdots d_{i_k j_k}
( x_{i_1}\cup x_{j_1}\cup\cdots\cup x_{i_k}\cup x_{j_k})
$$
That proves this lemma and the proposition.
$\Box$
\bigskip

This way, the correspondence $\Phi_{\E}\mapsto \phi_{\E}$
gives us the homomorphism
$$
\rho_{A}: Auteq({\db{A}})\lto \GL(H^{*}(A, \Z)),
$$
denoted by the same letter as in Lemma 4.1.4.3.

\medskip
\noindent{\bf 4.3.5 Definition.}
{\it The subgroup of
$\GL(H^{*}(A, \Z))$ that is the image of $\rho_{A}$ will
be called the {\sf  spinor group of }$A$ and will be denoted by
$Spin(A).$}

\medskip
\noindent{\bf 4.3.6.}
To justify this definition consider the Clifford algebra $Cl(\La ,Q)$ (3.1)
and
the homomorphism of algebras
$$
cor_A: Cl(\La, Q)\lto \End(H^*(A, \Z))\cong H^{*}(A\ts A, \Z)
$$
defined on the $\La\subset Cl(\La, Q)$ by rule
$$
cor_A ((l, x))(\cdot)=l(\cdot)+ x\cup(\cdot),
$$
where $(l, x)\in \Ga _A\op\Ga _A^{*}=\La$
and $l(\cdot)$ is the convolution with $l,$ i.e.
for a monomial
$\alpha=x_1\cup...\cup x_k$
it is defined as
$$
l(\alpha):=
\sum(-1)^{i-1} l(x_i)\cdot x_1\cup...\cup x_{i-1}\cup x_{i+1}\cup...\cup x_k .
$$
By Proposition 3.2.1e) the homomorphism
$cor_A$ is an isomorphism.
Clearly it sends
$Cl^{+}(\La, Q)$
onto
$H^{ev}(A\ts A, \Z).$
The homomorphism
$cor_{A}$
induces the action of the group
$Spin(\La, Q)$ (3.4.1)
on cohomology lattice
$H^{*}(A, \Z).$
Henceforth, we will consider $Spin(\La, Q)$ as
a subgroup of
$\GL(H^*(A, \Z))$
with respect to
$cor_{A}.$

There is the standard involution $'$ on
$Cl(\La).$
It is defined as a unique ring involution on
$Cl(\La)$
that is the identity on $\La.$
The isomorphism
$cor_A$ induces an involution on
$\End(H^{*}(A, \Z)),$ which will be denoted the same symbol $'.$
It is not hard to check that for any
$\xi\in H^{d}(A\ts A, \Z)$ there is the equality
$$
{v_{\xi}}^{'}=(-1)^{n+\frac{d(d-1)}{2}}w_{\xi},
$$
where $v_{\xi}$ and $w_{\xi}$ defined in 4.1.3 and $n=\dim A.$
In particular, this implies that
\begin{equation}\label{invol}
{\phi_{\E}}^{'}=\psi_{\E^{\vee}[n]}.
\end{equation}
for any object $\E\in D^{b}(A\ts A).$

Thus there are two correspondences
$\Phi_{\E}\mapsto f_{\E}$
and
$\Phi_{\E}\mapsto \phi_{\E},$
which gives two group homomorphisms
$$
\gamma_{A}: Auteq({\db{A}})\lto \U(A)\subset SO(\La, Q),
\qquad \rho_{A}: Auteq({\db{A}})\lto Spin(A)\subset \GL(H^{*}(A, \Z))
$$
such that
$\gamma_{A}(\Phi_{\E})=f_{\E}$ and
$\rho_{A}(\Phi_{\E})=\phi_{\E}.$
The following proposition relates these two homomorphisms.

\medskip
\noindent{\bf 4.3.7 Proposition.}\label{inspin}{\it
\begin{list}{\alph{tmp})}
{\usecounter{tmp}}
\item The group $Spin(A)$ is contained in $Spin(\La, Q)\subset \GL(H^{*}(A,
\Z)).$
\item Let
$\pi: Spin(\La, Q)\lto SO(\La, Q)$
be the canonical map.
Then
$\pi\cdot\rho_{A}=\gamma_{A}$
and
$\pi^{-1}(\U(A))=Spin(A).$
\end{list}}

\pr
Let
$\Phi_{\E}:{\db{A}}\stackrel{\sim}{\lto}{\db{A}}$
be an autoequivalence.
The following diagram is a
cohomological analog of diagram (\ref{diag}) with $B=A$:
\begin{equation}\label{diag2}
\begin{array}{ccc}
H^{*}(A\ts\wh{A}, \Z)&\stackrel{\phi_{S_A}}{\lto}&H^{*}(A\ts A, \Z)\\
\llap{\ss{\phi_{{\mathcal{J}}({\mathcal{E}})}}}\da&&\da\rlap{\ss{ad_{\E}}}\\
H^{*}(A\ts\wh{A}, \Z)&\stackrel{\phi_{S_A}}{\lto}&H^{*}(A\ts A, Z)
\end{array}
\end{equation}

By Theorem 2.4.5.1 the sheaf
$\J(\E)$
is a line bundle over the subvariety
$\Ga(f_{\E})$
that is the graph
of the isomorphism
$f_{\E}.$
This implies that
$\phi_{\J(\E)}$ sends $H^{4n-1}(A\ts\wh{A}, \Z)$
to itself.
Moreover, the restriction of
$\phi_{\J(\E)}$
on the $(4n-1)-th$ cohomology
coincides with
$$
f_{\E}: H_1 (A\ts \wh{A}, \Z)\lto H_1 (A\ts \wh{A}, \Z)
$$
under the Poincare isomorphism
$D:H_1(A\ts\wh{A}, \Z)=\La \stackrel{\sim}{\lto}H^{4n-1}(A\ts\wh{A}, \Z).$
Fix this identification of
$H^{4n-1}(A\ts \wh{A}, \Z)$
with
$\La.$

Let  us assume that  the restriction of
$\phi_{S_A}$ on $H^{4n-1}(A\ts \wh{A}, \Z)$
coincides with the restriction of $cor_A$ on $\La.$
Under this assumption, it follows from the commutativity of the diagram
(\ref{diag2}) that
$ad_{\E}$ takes
$cor_A(\La)$
to itself.
Hence $\phi_{\E} cor(\La) \phi_{\E}^{-1}= cor(\La).$
Further, since the inverse of an equivalence
$\Phi_{\E}$ is isomorphic to
$\Psi_{\E^{\vee}[n]}$
we get from (\ref{invol}) that
${\phi_{\E}}'=\phi_{\E}^{-1}.$
Therefore $N(\phi_{\E})=\phi_{\E}\cdot{\phi_{\E}}'=id.$
By Definition 3.4.1  $\phi_{\E}$ belongs to $Spin(\La, Q)$
and a) is proved.

Moreover, we have that action of $f_{\E}$ on $\La$
coincides with the action of $ad_{\E}$ on $cor_{A}(\La).$ This implies that
$\gamma_{A}=\pi \rho_{A}.$
Finally, since $Ker\pi=\Z/2\Z$ and
$\rho([1])=-Id$
the full inverse image
$\pi^{-1}(\U(A))$ not only contains but also coincides
with $Spin(A).$

Thus to complete the proof of the proposition  it is sufficient to prove
the following lemma.

\medskip
\noindent{\bf 4.3.7.1 Lemma.} {\it
The restriction of
the map $\phi_{S_A}D$ on $H_{1}(A\ts \wh{A}, \Z)= \La$
coincides with the restriction of the map $cor_A$ on $\La.$}

\smallskip
\pr
As above let $x_1,...,x_{2n}$ and $l_1,..., l_{2n}$ be the dual bases
of
$\Ga _A^*$
and
$\Ga _A$ respectively. For any
$\alpha\in H_{i}(A\ts \wh{A}, \Z)$ and $\beta\in H^{i}(A\ts \wh{A}, \Z)$
the following
equality holds
$$
\langle \beta, \alpha\rangle=\langle \beta\cup D(\alpha), [A\ts \wh{A}]\rangle,
$$
where
$\langle, \rangle$
is the canonical pairing of cohomology with homology and
$[A\ts\wh{A}]\in H_{4n}(A\ts \wh{A}, \Z)$ is the fundamental class.
It is easy to check that
$$
\begin{array}{l}
D(x_i)= (-1)^{i-1}
p^{*}(x_1\cup\cdots\cup x_{2n})\cup
q^{*}(l_1\cup\cdots\cup l_{i-1}\cup\l_{i+1}\cup\cdots\cup l_{2n}),\\
D(l_i)= (-1)^{i-1}
p^{*}(x_1\cup\cdots\cup x_{i-1}\cup x_{i+1}\cup\cdots\cup x_{2n})\cup
q^{*}(l_1\cup\cdots\cup l_{2n}).
\end{array}
$$

Denote $\mu _A=\mu $ (4.2.4). By Lemma 4.2.4.1 the functor
$\Phi_{S_A}$
is the composition
${\bf R}\mu_{*}\Phi_{\cP_{A}}.$
Hence
$\phi_{S_A}=\mu_{*}\phi_{\cP_{A}}.$
By construction of
$\Phi_{\cP_{A}},$ the map
$\phi_{\cP_{A}}$ takes
$q^{*}(l)\cup p^{*}(x)$ to
$p_{1}^{*}(\psi_{P}(l))\cup p_{2}^{*}(x).$
Combining Lemma 4.2.3.2
and Proposition 4.2.2.1
we get
\begin{eqnarray}
\phi_{\cP_{A}}D(x_i)&=&p_{1}^* (x_i)\cup p_{2}^* (x_1\cup\cdots\cup x_{2n})
\label{comp1}\\
\phi_{\cP_{A}}D(l_i)&=& (-1)^{i-1}p_{2}^* (x_1\cup\cdots\cup x_{i-1}\cup
x_{i+1}\cup\cdots\cup x_{2n})\label{comp2}
\end{eqnarray}
Further, to compute $\mu_{*}(\alpha)$ it is sufficient to
note that $\mu_{*}$ is a ring homomorphism because $\mu$ is an
isomorphism and, consequently, $\mu_{*}$
is the inverse of $\mu^{*}.$
Besides, since $\mu(a_1, a_2)=(a_1, m(a_1, a_2))$
we have
$$
\mu^{*}(p_1^{*}(\alpha))=p_{1}^{*}(\alpha), \qquad
\mu^{*}(p_2^{*}(\alpha))=m^{*}(\alpha)
$$
for any $\alpha\in H^{*}(A, \Z).$
Moreover, if $\alpha=x\in H^{1}(A, \Z)$ then
$$
\mu^{*}(p_{2}^{*}(x))=m^{*}(x)=p_{1}^{*}(x)+p_{2}^{*}(x)
$$
Therefore, for any $x\in H^{1}(A, \Z)$
there are equalities
$$
\mu_{*}(p_{1}^{*}(x))=p_{1}^{*}(x)\quad\mbox{and}\quad
\mu_{*}(p_{2}^{*}(x))=p_{2}^{*}(x) - p_{1}^{*}(x).
$$

To find  an element
$\phi_{S_A}(\lambda)$
for some
$\lambda\in H^{4n-1}(A\ts\wh{A}, \Z)$
it is sufficient to compute the map that is defined by this element, i.e.
$$
\phi_{S_A}(\lambda)(\cdot):=p_{2*}(\phi_{S_A}(\lambda)\cup p_{1}^{*}(\cdot)).
$$
Substituting (\ref{comp1}) in the last expression, we get
$$
\begin{array}{l}
\phi_{S_A}D(x_i)(\alpha)=p_{2*}( \mu_{*}(p_{1}^*(x_i)\cup
p_{2}^{*}(x_1 \cup\cdots\cup x_{2n}))\cup p^{*}_{1}(\alpha))=
m_{*}(p_{1}^{*}(x_i\cup \alpha)\cup p_{2}^{*}(x_1\cup\cdots\cup x_{2n}))\\
=m_{*}(m^{*}(x_i\cup \alpha)\cup p_{2}^{*}(x_1\cup\cdots\cup x_{2n}))=
 x_{i}\cup \alpha\cup m_{*}(p_{2}^{*}(x_1\cup\cdots\cup x_{2n}))=
 x_{i}\cup \alpha
\end{array}
$$
Hence $\phi_{S_A}D(x)$ coincides with
 $cor_{A}(x).$

Let us compute the action of
$\phi_{S_A}D(l_i)$ on the cohomology. We have
$$
\begin{array}{l}
\phi_{S_A}D(l_i)(x_{s_1}\cup\cdots\cup x_{s_k})=
p_{2*}(\mu_{*}((-1)^{i-1}p_{2}^*(\bigcup\limits_{j=1,\atop j\ne i}^{2n}x_{j}))
\cup p_{1}^{*}(x_{s_1}\cup\cdots\cup x_{s_k}))
=\\
(-1)^{i}p_{2*}(\bigcup\limits_{j=1,\atop j\ne i}^{2n}(p^{*}_{1}(x_j)-p^{*}_{2}(x_j))
\cup p_{1}^{*}(x_{s_1}\cup\cdots\cup x_{s_k}))
\end{array}
$$
If the set $(s_1,...,s_k)$ does not contain $i,$ then this expression equals
$0.$
Suppose that $s_1=i.$ In this case the last expression can be simplified
$$
\begin{array}{l}
\phi_{S_A}D(l_i)(x_{i}\cup x_{s_2}\cup\cdots\cup x_{s_k})=
p_{2*}(p^{*}_{1}(x_{1}\cup\cdots\cup x_{2n})\cup
p_{2}^{*}(x_{s_2}\cup\cdots\cup x_{s_k}))=x_{s_2}\cup\cdots\cup x_{s_k}
\end{array}
$$
Thus,
$\phi_{S_A}D(l_i)$ coincides with
$cor_A(l_i),$ because
$$
\begin{array}{l}
cor_{A}(l_i)(x_{s_1}\cup\cdots\cup x_{s_k})=0\quad\mbox{if}
\quad i\not\in(s_1,...,s_k),\\
cor_{A}(l_i)(x_{i}\cup x_{s_2}\cup\cdots\cup x_{s_k})=
x_{s_2}\cup\cdots\cup x_{s_k}.
\end{array}
$$
This concludes the proof of the lemma and, consequently, completes
the proof of the proposition.
$\Box$

\medskip

\noindent{\bf 4.3.8 Corollary.}\label{ac}{\it
$Ker \rho_A$
coincides with
$2\Z\ts A\ts\wh{A},$
consisting of all equivalences
$\Phi_{(a,\alpha)}[2i]= T_a^{*}(\cdot)\ot P_{\alpha}[2i].$}

\smallskip
\pr  It is clear that
$\Phi_{(a,\alpha)}[2i]$ are contained in $Ker \rho_{A}.$ On the other side,
it follows from Proposition 4.3.7 that
$Ker \rho_{A}$
is a subgroup of
$Ker \gamma_{A}\cong \Z\ts A\ts\wh{A}.$ Note that by Lemma 4.1.5.1
$\rho_{A}([1])=-Id.$ Corollary is proved.
$\Box$

\medskip
\noindent{\bf 4.3.9 Corollary.} {\it
There is the exact sequence of groups
$$
0\lto \Z/2\Z \lto Spin(A)\lto \U(A)\lto 1.
$$}

The last theorem of this chapter gives description
of abelian varieties that have equivalent derived categories of coherent sheaves.

\medskip
\noindent{\bf 4.3.10 Theorem.}(\cite{Or2}) { \it Let $B$ and $C$ be abelian
varieties. Then the derived categories $D^b(B)$ and $D^b(C)$ are
equivalent if and only if there exists an isomorphism
$$\gamma :B\times \wh{B}\stackrel{\sim }{\lto}C\times \wh{C}$$
which identifies the forms $Q_B$ and $Q_C$ on
$\Ga_{\vphantom{\wh{A}}B}\oplus \Ga_{\wh{B}}$ and $\Ga_{\vphantom{\wh{A}}C}\oplus \Ga_{\wh{C}}.$}.

\medskip
The ``if'' direction in the above theorem was first proved by A. Polishchuk
in \cite{Po1}.

%
%
%
%

\sec{The algebraic group $\U_{A,\Q}.$}

\noindent{\bf 5.1.}
Fix a complex torus $A.$ Put $\La :=\Ga_{\vphantom{\wh{A}}A}\oplus \Ga_{\wh A}$ and consider
the group $\GL(\La _{\Q}).$ We have the following $\Q$-algebraic subgroups of
$\GL(\La _{\Q}).$

1) $Hdg_{\vphantom{\wh{A}}A,\Q}=Hdg_{A\times \wh{A},\Q}.$

2) $Aut _{\Q}(A\times \wh A),$ which is defined as the group of invertible
elements in $\End^0(A\times \wh A).$

3) $O(\La _{\Q}, Q_{\Q}),$ $SO(\La _{\Q},Q_{\Q}),$ where $Q_{\Q}$ is the
extension
to
$\La _{\Q}$ of the canonical symmetric bilinear form on $\La $ (3.1).

4) $\U_{A,\Q},$ which is defined as follows
$$
\U_{A,\Q}=
\left\{
g=\left(
\begin{array}{cc}
a&b\\
c&d
\end{array}
\right)
\in
Aut_{\Q}(A\times \wh A)
\left|
\quad
g^{-1}=
\left(
\begin{array}{rr}
\wh{d}&-\wh{b}\\
-\wh{c}&\wh{a}
\end{array}
\right)
\right\}
\right.
$$
Thus the discrete group $\U(A)$ defined in 4.3.1
is the arithmetic subgroup of $\U_{A,\Q}$ consisting
of elements which preserve the lattice $\La.$

The following proposition summarizes the interrelations of these subgroups.

\medskip
\noindent{\bf 5.2 Proposition.} {\it

1) $Hdg_{A,\Q} \subset \SO (\La _{\Q},Q_{\Q}).$

2) $Aut_{\Q}(A\times \wh A)$ is the centralizer of $Hdg_{A,\Q}$ in
$\GL(\La _{\Q}).$

3) $\U_{A,\Q}=Aut_{\Q}(A\times \wh A)\cap \SO(\La _{\Q},Q_{\Q}).$

4)  $\U_{A,\Q}$ is the centralizer of $Hdg_{A,\Q}$ in
$\SO (\La _{\Q},Q_{\Q}).$}

\smallskip
\pr
1) The operator $J_A$ of complex structure on $\La _{\R}$ has determinant 1
and preserves the form $Q.$ Thus $h_A(\S^1)$ lies in the group of
$\R$-points of
$\SO (\La _{\Q},Q_{\Q}).$ Since $Hdg_{A,\Q}$ is the $\Q$-closure of
$h_A(\S^1),$
we have
$Hdg_{A,\Q}\subset \SO(\La _{\Q},Q_{\Q}).$

2) By definition $\End^0(A\times \wh A)$ is the centralizer of $Hdg _{A,\Q}$
in $\End(\La _{\Q}).$ Thus $Aut_{\Q}(A\times \wh A)$ is the centralizer of
$Hdg_{A,\Q}$ in $\GL(\La _{\Q}).$

3) The proof is the same as that of proposition 4.3.2 above.

4) This follows from 1), 2), 3).
$\Box$

\medskip
\noindent{\bf 5.3.} Let us study the $\Q$-algebraic group $\U_{A,\Q}$
in case $A$ is an
abelian variety.

Clearly, the group $U_{A,\Q}$ depends only on the isogeny class of $A.$
If $A=A_1^{m_1}\times ... \times A_k^{m_k},$ where $A_i$ are pairwise
nonisogeneous simple abelian varieties, then
$$\U_{A,\Q}=\prod _i\U_{A_i^{m_i},\Q}.$$

Below we describe the $\Q$-algebraic
group  $\U_{A^m,\Q}$
and its group of $\R$-points $\U_{A^m,\Q}(\R)$
for a simple abelian variety
$A.$

\medskip
\noindent{\bf 5.3.1.} Let us use the notations of section 1.8.
In particular denote
$F=\End^0(A).$ Let $\varphi \in \Hom(A,\wh A)$ be a polarisation and
$':F\to F$ be the corresponding Rosati involution. We identify naturally
$M(m, F)=\End ^0(A^m),$ which makes $\Ga _{A^m,\Q}$ a left $M(m, F)$
module. Consider the diagonal polarisation $\sigma =(\varphi,...\varphi )$
of $A^m.$ Then the corresponding Rosati involution on $M(m, F)$ is
$$Z\mapsto {}^tZ'.$$

Consider the following homomorphism of algebras
$$\tau : \End ^0(A\times \wh A)\lto \End (\Ga _{A^m,\Q}\oplus\Ga _{A^m,\Q})$$
$$\tau : \left(\begin{array}{cc}
B & C\\
D & E
\end{array}\right)
\mapsto
\left(\begin{array}{cc}
B & C\sigma\\
\sigma ^{-1}D & \sigma ^{-1}E\sigma
\end{array}\right).$$
Obviously, the image of $\tau $ is contained in $M(2m, F).$ We have
$$
\tau(\U_{A^m,\Q})=
\left\{
g=\left(
\begin{array}{cc}
a&b\\
c&d
\end{array}
\right)
\in
\GL(2m,F)
\left|
\quad
g^{-1}=
\left(
\begin{array}{rr}
{}^t\wh{d}&-{}^t\wh{b}\\
-{}^t\wh{c}&{}^t\wh{a}
\end{array}
\right)
\right\}
\right.
$$
Thus $\tau(\U_{A^m,\Q})$ is the group of isometries of the skew-hermitian
form on $F^m\oplus F^m$
$$\vartheta ((z_1,...,z_m,z_{-1},...,z_{-m}),
(w_1,...,w_m,w_{-1},...,w_{-m}))=\sum_{k=1}^{m}z_kw'_{-k}-z_kw_k'.$$

\medskip
\noindent{\bf 5.3.2.}
Denote by $U^*_{2m,F}\subset \GL(2m,F)$ the $K_0$-algebraic group of
isometries of the skew-hermitian form $\vartheta.$ Then the
$\Q$-algebraic group $\U_{A^m,\Q}$ is obtained from $U^*_{2m,F}$ by
restriction of scalars from $K_0$ to $\Q.$ Fix an embedding $K_0
\hookrightarrow \R.$ The corresponding group of real points
$U^*_{2m,F}(\R)$ was computed by A. Polishchuk (\cite{Po2}). His result
according to the four cases of Albert's classification is the
following
\begin{list}{\Roman{tmp}.}%
{\usecounter{tmp}}
\item $Sp_{2m}(\R),$
\item $Sp_{4m}(\R),$
\item $U^*_{2m}(\Qu )$ -- the group of automorphisms of $\Qu ^{2m}$
preserving the standard skew-hermitian form,
\item $U(md,md).$
\end{list}

\medskip
\noindent{\bf 5.3.3 Corollary.} {\it The reductive Lie group $U^*_{2m,F}(\R)$
is semisimple unless we are in case IV. This group has no compact factors
unless we are in case IV (then it is isogeneous to the product
$\S^1\times SU(md,md)$) or in case III and $m=1$ (then it is isogeneous to
the product $SU(2)\times \SL(2,\R)$).}

\medskip

\noindent{\bf 5.3.4.} We have
$$\U_{A^m,\Q}(\R)\cong \prod _{e_0} U^*_{2m,F}(\R),$$
hence the above provides the corresponding description of the group
$\U _{A^m,\Q}(\R).$

\medskip
\noindent{\bf Corollary.} {\it The reductive group $\U_{A^m,\Q}(\R)$
is semisimple unless we are in case IV. This group has no compact
factors unless we are in case IV (then the compact part is isogeneous to
the product of $e_0$ copies of $\S^1$) or in case III and $m=1$ (then
the compact part is isogeneous to the product of $e_0$ copies of $SU(2)$).}

\medskip
\noindent{\bf 5.3.5 Corollary.} {\it For any abelian variety $A$ the
Lie group $\U _{A,\Q}(\R)$ is connected.}

\sec{ The Neron-Severi Lie algebra ${\frak g}_{NS}(X).$}

\noindent{\bf 6.1.}
Let $X$ be a smooth projective variety of dimension $n.$ Let us
recall the definition of the Neron-Severi Lie algebra ${\frak g}_{NS}(X)$
from \cite{LL}. If $\kappa \in H^{1,1}(X)\cap H^2(X,\Q)$ is an ample class, then
cupping with it defines an operator $e_{\kappa }$ in the total cohomology
$H^*(X)=H^*(X,\C)$ of degree 2 and the hard Lefschetz theorem asserts
that for $s=0,1,...,n,$ $e^s_{\kappa }$ maps $H^{n-s}(X)$ isomorphically
onto $H^{n+s}(X).$ As is well known, this is equivalent to the existence
of a (unique) operator
$f_{\kappa}$ on $H^*(X)$ of degree $-2$ such that the commutator $[e_{\kappa },f_{\kappa }]$ is the operator $h$ which on $H^k(X)$ is multiplication by $k-n.$
The elements $e_{\kappa },f_{\kappa },h$ make up a Lie subalgebra ${\frak g}
_{\kappa }$ of ${\frak gl}(H^*(X))$ isomorphic to $sl(2).$ Define the
{\it Neron-Severi} Lie algebra ${\frak g}_{NS}(X)$ as the Lie subalgebra
of ${\frak gl}(H^*(X))$ generated by ${\frak g}_{\kappa}$'s with $\kappa $ an
ample class. This Lie subalgebra is defined over $\Q$ and is evenly
graded by the adjoint action by the semisimple element $h.$

In the above definition we could replace an ample class $\kappa $ by a
Kahler class or by any element $\kappa \in H^2(X),$ such that $e_{\kappa }$
satisfies the conclusion of the hard Lefschetz theorem. The resulting Lie
subalgebras of ${\frak gl}(H^*(X))$ are called the {\it Kahler} Lie algebra
and
the {\it total} Lie algebra respectively and denoted by ${\frak g}_{K}(X)$
and ${\frak g}_{tot}(X)$ respectively. Obviously we have the inclusions
${\frak g}_{NS}(X)\subset {\frak g}_{K}(X)\subset {\frak g}_{tot}(X).$

The above Lie algebras preserve (infinitesimally) the following form
on the cohomology $H^*(X).$
$$\chi (\alpha,\beta):=(-1)^q\int_X\alpha\cup \beta,$$
where $\alpha $ is homogeneous of degree $n+2q$ or $n+2q+1.$ The main fact
about these Lie algebras is that they are semisimple (\cite{LL},p.369).

\medskip
\noindent{\bf 6.2.}
In case $X$ is a complex torus (resp. an abelian variety) the Lie algebras
${\frak g}_{tot}(X),$ ${\frak g}_{K}(X)$ (resp. ${\frak g}_{NS}(X)$) were
computed explicitly in \cite{LL}, p.381. Let us recall the result.

\medskip
\noindent{\bf 6.2.1.}
Let $X=(V_X/\Ga _X,J_X)$ be a complex torus.
Put $\La =\Ga_{\vphantom{\wh{A}}X}\oplus \Ga_{\wh X}$ with the canonical symmetric bilinear
form $Q$ (3.1).
There exists a natural isomorphism of Lie
algebras
$${\frak g}_{tot}(X)\cong {\frak so}(\La _{\Q},Q _{\Q}),$$
such that the semisimple element $h\in{\frak g}_{tot}(X)$ corresponds to
the element
$-1_{\Ga _X}\oplus1_{\Ga _{\wh X}}\in {\frak so}(\La _{\Q},Q_{\Q}).$
The natural
representation of ${\frak g}_{tot}$ on $H^*(X,\Q)$ under the above isomorphism
is the spinorial representation of ${\frak so}(\La _{\Q},Q_{\Q})$ (which is
the direct sum of two semi-spinorial ones corresponding to $H^{ev}$ and
$H^{odd}$ respectively).

Let $J=(J_{\vphantom{\wh{A}}X},J_{\wh X})$ be the complex structure on $V_{\vphantom{\wh{A}}X}\oplus V_{\wh X}.$
Then the
form $Q_{\R}$ is $J$-invariant, hence defines a hermitian form on $\La _{\R}.$
Let $\frak{su}(\La _{\R})$ denote the Lie algebra of the corresponding
special unitary Lie group. Then the isomorphism mentioned above
${\frak g}_{tot}(X)\cong {\frak so}(\La _{\Q},Q_{\Q})$ identifies the Lie
subalgebras ${\frak g}_{K}(X;\R)\cong {\frak su}(\La _{\R}).$

\medskip
\noindent{\bf 6.2.2.}
Finally, let $X$ be an abelian variety. The Lie algebra ${\frak g}_{NS}(X)$
depends only on the isogeny type of $X.$ If
$$X=X_1^{m_1}\times ...\times X_k^{m_k},$$
then
$${\frak g}_{NS}(X)={\frak g}_{NS}(X_1^{m_1})\times ... \times
{\frak g}_{NS}(X_k^{m_k}).$$

We therefore assume that $X$ is a power $A^m$ of a simple abelian variety
$A.$ Let us stick to notations of section 1.8. In particular, $F=\End^0(A).$
Choose a polarization $\varphi \in \Hom(A,\wh A)$ of $A.$ Denote by $'$ the
corresponding Rosati involution on $F.$ Consider the diagonal polarisation
$\sigma=(\varphi,...,\varphi)$ of $A^m.$ If we naturally identify
$M(m, F)=\End^0(A^m)$ then the Rosati involution defined by $\sigma$ is
$Z\mapsto {}^tZ'.$

Define a $K_0$--Lie subalgebra of $M(2m,F)$ by
$${\frak slu}(2m,F,')=\{
\left(\begin{array}{cc}
B & C\\
D & -{}^tB'
\end{array}\right)
| B,C,D\in M(m,F);\ \ C={}^tC',\ \ D={}^tD'\}.$$
This is the Lie algebra of infinitesimal isometries of the skew-hermitian
form $\vartheta $ on $F^m\oplus F^m,$ which was defined in 5.3.1 above.
It is a
reductive $K_0$--Lie algebra whose center is the space of scalars $\lambda \in
K$ with $\lambda'=-\lambda.$ So ${\frak slu}(2m,F,')$ is semisimple unless
we are in the case of totally complex multiplication (case IV).
We grade this Lie
algebra  by means of the semisimple element
$$u_m:=\left(
\begin{array}{cc}
-1_m & 0\\
0 & 1_m
\end{array}\right)
\in {\frak slu}(2m,F,'),$$
so that $B,$ $C$ and $D$ parametrize the summands of degree $0,$ $-2$
and $2$ respectively. Let ${\frak g}(2m,F,')$ denote the $K_0$--Lie subalgebra
of ${\frak slu}(2m,F,')$ generated by summands of degree 2 and -2; let
${\frak u}(m,F,')$ denote the union of $\GL(m,F)$--conjugacy classes in
$M(m,F)$ made up of anti-invariants with respect to the involution
$Z\mapsto {}^tZ'$ and identify ${\frak u}(m,F,')$ with the
subspace of ${\frak slu}(2m,F,')$ in an obvious way. We have the following
result.

\medskip

\noindent{\bf 6.2.2.1 Lemma.} (\cite{LL}, Lemma 3.9) {\it We have
$${\frak slu}(2m,F,')={\frak g}(2m,F,')\times {\frak u}(m,F,').$$
The summand ${\frak u}(m,F,')$ is trivial except in the following cases

1) $m=1$ and $F$ is totally definite quaternion (case III): then
${\frak g}(2,F,')
\cong sl(2,K_0)$ and ${\frak u}(1,F,')$ can be identified with pure
quaternians in $F$ (i.e. the $'$-antiinvariants in $F$) or

2) $K$ is totally complex (case IV): then ${\frak g}(2m,F,')$ consists of
the matrices
for which $B$ has its $K$--trace in $K_0,$ whereas ${\frak u}(m,F.')$ can be
identified with the purely imaginary scalars in $K$ (i.e. the $\lambda \in
K$ with $\bar{\lambda}=-\lambda$). }

\medskip

\noindent{\bf 6.2.2.2 Remark.} It was noted in \cite{LL} that in the exceptional
cases the
connected Lie subgroup of $\GL(m,F\otimes _{\Q}\R)$ with the Lie algebra
${\frak u}(m,F,')\otimes _{\Q}\R$ is a product of $e_0$ copies of $\U(1),$
resp $\SU(2),$ and hence is compact.

\smallskip

In view of 6.2.1 above we may consider ${\frak g}_{NS}(X)$ as a Lie subalgebra
of $\End (\La _{\Q})=\End(\Ga_{\vphantom{\wh{A}}X,\Q}\oplus \Ga _{\wh X,\Q}).$
Recall the homomorphism of
algebras $\tau $ as in 5.3.1 above.
$$\tau : \End ^0(\Ga_{\vphantom{\wh{A}}X,\Q}\oplus \Ga _{\wh X,\Q})\lto \End
(\Ga_{X,\Q}\oplus \Ga_{X,\Q})$$
$$\tau : \left(\begin{array}{cc}
B & C\\
D & E
\end{array}\right )
\mapsto
\left(\begin{array}{cc}
B & C\sigma\\
\sigma ^{-1}D & \sigma ^{-1}E\sigma
\end{array}\right).$$
Note that $\Ga _{X,\Q}\oplus \Ga _{X,\Q}$ is naturally a $M(2m,F)$-
module.

\medskip
\noindent{\bf 6.2.2.3 Proposition.} (\cite{LL}, Prop.3.10) {\it The above
homomorphism
$\tau$ identifies the Lie algebras ${\frak g}_{NS}(X)$ and
${\frak g}(2m,F.').$ In particular, we get an isomorphism
$${\frak slu}(2m,F,')\cong {\frak g}_{NS}(X)\times {\frak u}(m,F,').$$}

\sec{Relation between the group $Auteq(D^b(A))$ and the Neron-Severi
Lie algebra ${\frak g}_{NS}(A)$ for an abelian variety $A.$ }

\noindent{\bf 7.1.}
Let $A$ be an abelian variety. We know that the group
$Auteq(D^b(A))$ acts on the total cohomology $H^*(A,\Z)$ of $A$ and
we denoted its image in $\GL(H^*(A,\Z))$ by $Spin(A)$ (4.3.5).
Put $\La =\Ga_{\vphantom{\wh{A}}A}\oplus \Ga _{\wh A}$ and let $Q$ be the canonical
symmetric bilinear form on $\La $ (3.1).
Recall the exact sequence of group homomorphisms
$$0 \lto \Z/2\Z \lto Spin(A) \lto \U(A) \lto 1$$
(4.3.9) and the commutative diagram (4.3.7)
$$\begin{array}{ccc}
Spin(A) & \hookrightarrow & Spin(\La,Q) \\
\downarrow & & \downarrow \\
\U(A) & \hookrightarrow & \SO(\La,Q)
\end{array},$$
where horizontal arrows are compatible with the actions on $H^*(A,\Z)$ and
$\La $  respectively.

\medskip
\noindent{\bf 7.1.1 Definition.}
{\it Let
$\overline{Spin(A)}$ (resp. $\overline{\U(A)}$)
be the {\sf  Zariski closure} of $Spin(A)$ (resp. $\U(A)$) in
$\GL(H^*(A,\Q))$ (resp $\GL (\La _{\Q})$).}

 These are  $\Q$-algebraic groups.
These groups are isogeneous, hence have isomorphic Lie algebras.

\medskip
\noindent{\bf 7.1.2 Remark.} Since the group $Spin(A)$ acts on $H^*(A,\Q)$
by algebraic correspondences its action commutes with $Hdg_{A,\Q}.$
Hence the same is true for $\overline{Spin(A)}.$

\medskip

On the other hand we know that ${\frak g}_{NS}(A)$
can be identified as a Lie subalgebra of ${\frak so}(\La _{\Q},Q_{\Q})$
so that the
representation of ${\frak g}_{NS}(A)$ in $H^*(A,\Q)$ is the restriction of the
spinorial representation of ${\frak so}(\La _{\Q},Q_{\Q})$ (6.2.1).

\medskip
\noindent{\bf 7.2 Theorem.} {\it

a) The Lie algebra of $\overline{Spin(A)}$ is ${\frak g}_{NS}(A)$ considered
as a Lie subalgebra of ${\frak gl}(H^*(A,\Q)).$

b) The Lie algebra of $\overline{\U(A)}$ is ${\frak g}_{NS}(A)$
considered as a Lie subalgebra of ${\frak gl}(\La _{\Q}).$}

\smallskip
\pr
The discussion in 7.1 above implies that a) and b) are equivalent.
Let us prove b).

Both $\overline{\U(A)}$ and ${\frak g}_{NS}(A)$ depend only on the
isogeny type of $A.$ If $A$ is a product of pairwise nonisogeneous
abelian varieties then both $\overline{\U(A)}$ and ${\frak g}_{NS}(A)$
are products of the corresponding factors. Thus it suffices to prove that
${\frak g}_{NS}(A^m)$ is the Lie algebra of the algebraic group
$\overline{\U (A^m)}$ for a simple abelian variety $A.$

By definition $\U (A^m)$ is an arithmetic subgroup of the
reductive $\Q$-algebraic group $\U_{A^m,\Q}$ (5.1).
By Corollary 5.3.4 all abelian factors of $\U_{A^m,\Q}$ are
compact.
Thus by the density theorem (\cite{PR}) the $\Q$-algebraic group
$\overline{\U (A^m)}$ consists (up to isogeny) of all noncompact
factors of $\U_{A^m,\Q}.$ Let ${\frak g} (A^m)$ be the Lie
algebra of $\U_{A^m,\Q}.$ As follows from the results in 5.3.1, 6.2.2,
6.2.2.3  that it is
a product of Lie algebras
$${\frak g} (A^m)={\frak g}_{NS}(A^m)\times {\frak u}_{A^m},$$
where ${\frak u}_{A^m}=\tau^{-1}({\frak u}(m,F,'))$ in the notation of
5.3.1, 6.2.2  above. Moreover, it follows from the results in 5.3.4,
6.2.2.1, 6.2.2.2 above that the
$\Q$-algebraic subgroup of $\U_{A^m,\Q}$ corresponding to
${\frak u}_{A^m}$ consists (up to isogeny) of all compact factors of
$\U_{A^m,\Q}.$ Thus ${\frak g}_{NS}(A^m)$ is the Lie algebra of
the group $\overline{\U (A^m)}.$ This proves the theorem.
$\Box$

\medskip
\noindent{\bf 7.2.1 Corollary.} {\it

a) The algebraic groups $\overline{Spin(A)}$
and $\overline{\U(A)}$ are semisimple.

b) The group $\overline{\U(A)}$ consists (up to isogeny) of all
noncompact factors of $\U_{A,\Q}.$}

\smallskip
\pr a) Indeed, this follows from the above theorem, since the Lie algebra
${\frak g}_{NS}(A)$ is semisimple. b) This was obtained in the proof of the
above theorem.
$\Box$

\medskip
\noindent{\bf 7.3.}  Let $A$ be an abelian variety.
Consider ${\frak g}_{NS}(A;\R)$ as a Lie subalgebra of
${\frak gl}(\La _{\R}).$ Let $\varphi \in NS_A(\R)^0\subset
\Hom (V_A,V_{\wh A}).$ Then $\varphi,\varphi ^{-1}\in {\frak g}_{NS}
(A;\R)$ and
$$[\varphi, \varphi ^{-1}]=h=\left(\begin{array}{cc}
                                   -1_A & 0 \\
                                     0  & 1_{\wh A}
                                   \end{array}\right).$$
The elements $\varphi, \varphi ^{-1}, h$ make up a Lie subalgebra
${\frak g}_{\varphi }$ of ${\frak g}_{NS}(A;\R)$ isomorphic to $sl(2).$
By definition the Lie algebra ${\frak g}_{NS}(A;\R)$ is generated by these
subalgebras.

Consider the group of $\R$-points $\overline{\U (A)}(\R)$ of the
$\Q$-algebraic group $\overline{\U(A)}.$ It is a semisimple Lie
subgroup of the reductive Lie group $\U_{A,\Q}(\R)$ (5.3.4) which
consists (up to isogeny) of all noncompact factors of $\U_{A,\Q}(\R).$
By the above theorem the Lie algebra of $\overline{\U (A)}(\R)$ is
${\frak g}_{NS}(A;\R).$ For $\varphi \in NS_A(\R)^0$ denote by
$G_{\varphi }\subset \overline{\U (A)}(\R)$ the connected Lie subgroup
corresponding to the Lie subalgebra ${\frak g}_{\varphi }.$ We have
$G_{\varphi }\cong \SL (2;\R).$

\sec{Action of $\U_{A,\Q}(\R)$ on a Siegel domain.}

\noindent{\bf 8.1.}
Let
$A$
be a complex torus.
Let us define a rational
(i.e. not defined everywhere) action of the group $\U_{A,\Q}(\R)$ (5.1)
on the complex space
$NS_{A}(\C)=NS_{A}\ot \C\subset \Hom(A, \wh{A})\ot \C.$
It is given by the formula
\begin{equation}\label{action}
\left(
\begin{array}{cc}
a&b\\
c&d
\end{array}
\right)
\omega
:=
(c+d\omega)(a+b\omega)^{-1},
\quad
\left(
\begin{array}{cc}
a&b\\
c&d
\end{array}
\right)
\in \U_{A,\Q}(\R),
\quad
\omega\in NS_{A}(\C)\subset \Hom(A, \wh{A})\ot \C.
\end{equation}
Here multiplication is understood as composition of maps.

In case
$A$ is an abelian variety
$NS_{A}(\C)$ contains Siegel domains of the first kind on which this action
is
well-defined.
Namely, let
$C_{A}^{a}\subset NS_{A}(\R)$
be the ample cone of $A,$
which is defined as
$\R^{+}$-linear combinations of ample classes in
$NS_A .$
It is an open subset in $NS_A(\R).$
Put
$$
C^{\pm}_{A}:= NS_{A}(\R)\pm i C^{a}_{A}\subset NS_{A}(\C),$$
$$
C_{A}:= C^{+}_{A}\amalg C^{-}_{A}.
$$
Thus
$C_{A}\ne \emptyset$
if and only if the complex torus
$A$
is an abelian variety.

\medskip
\noindent{\bf 8.2 Theorem.} {\it
Let $A$ be an abelian variety.

a) The action of
$\U_{A,\Q}(\R)$
on
$C_{A}$
is well defined.

b) This action on $C^{+}_{A}$
(resp. $C^{-}_{A}$) is transitive.

c) The stabilizer of a point in
$C_{A}$
is a maximal compact subgroup of the Lie group
$\U_{A,\Q}(\R).$}

\medskip

\pr This theorem is proved in Appendix.

\medskip
\noindent{\bf 8.3 Remark.}  Let $A$ be an abelian variety. The action of
$\U_{A,\Q}(\R)$ on $C_A$ restricts to an action of the Lie subgroup
$\overline{\U (A)}(\R).$ Since this subgroup consists (up to
isogeny) of all noncompact factors of $\U_{A,\Q}(\R)$ its action of
$C_A^+$ (resp. $C_A^-$) is also transitive. Thus $C_A^+$ (resp. $C_A^-$) is
a hermitian symmetric space for the semisimple Lie group $\overline{\U (A)}
(\R).$

\medskip
\noindent{\bf 8.4.} Let $A$ be a complex torus.
As above, denote by $NS_A(\R)^0\subset NS_A(\R)$ the subset consisting of nondegenerate forms.
Assume that
$NS_{A}(\R)^0 \neq \emptyset.$
Put
$$
NS_{A}(\C)^0:=
\left\{
\varphi_1+i\varphi_2\in NS_{A}(\C)\,\left|
\;
\varphi_2  \in NS_A(\R)^0 \right\}
\right.
$$
Thus, $NS_A(\C)^0$ is a nonempty Zariski open subset in $NS_A(\C).$
(Note that if
$A$ is an abelian variety then
$C^{+}_{A}$
(resp. $C^{-}_{A}$)
is a connected component of
$NS_{A}(\C)^{0}).$

Given
$\omega=\varphi_1 + i\varphi_2 \in NS_{A}(\C)^{0}$
consider the following element
\begin{equation}\label{pred}
I_{\omega}:=
\begin{pmatrix}
\varphi_2^{-1}\varphi_1&-\varphi_2^{-1}\\
\varphi_2 +\varphi_1\varphi_2^{-1}\varphi_1&-\varphi_1\varphi_2^{-1}
\end{pmatrix}
=
\begin{pmatrix}
1&0\\
\varphi_1&1
\end{pmatrix}
\begin{pmatrix}
0&-\varphi_2^{-1}\\
\varphi_2&0
\end{pmatrix}
\begin{pmatrix}
1&0\\
-\varphi_1&1
\end{pmatrix}
\in
\U_{A,\Q}(\R)
\end{equation}
and the morphism of algebraic $\R$-groups
$$
\mu_{\omega}: \S^1\lto \U_{A,\Q}(\R)
$$
given by the formula
$\mu_{\omega}(e^{i\theta})=
\cos\theta
\cdot
Id
+
\sin\theta \cdot
I_{\omega},$
so that
$\mu_{\omega}(e^{i\frac{\pi}{2}})= I_{\omega}.$

The correspondence $\omega \mapsto I_{\omega}$ is injective.

\medskip
\noindent{\bf 8.4.1 Properties of $I_{\omega}$
and
$\mu_{\omega}.$}

\begin{enumerate}%
\item
$
I_{\overline{\omega}}=-I_{\omega}
,$ hence
$\mu_{\bar{\omega}}=\mu_{\omega}\cdot(-Id_{\S^{1}}).$\\
This is obvious.
\item
 $I_{\omega}^2=-Id.$\\
 This is a direct computation.

\smallskip
By $K_{\omega}\subset \U _{A,\Q}(\R)$
denote the stabilizer of
$\omega$ with respect to the action (\ref{action}),
and
$Z(K_{\omega})$
denote its center. One checks directly that
$$
K_{\omega}=\left\{
\left(
\begin{array}{cc}
a&b\\
c&d
\end{array}
\right)
\in
\U_{A,\Q}(\R)\,
\left|
\;
\begin{array}{l}
\varphi_2 a +\varphi_1 b\varphi_2 +\varphi_2 b\varphi_1=d\varphi_2;\\
\varphi_1 a +\varphi_1 b\varphi_1 -\varphi_2 b\varphi_2=c+ d\varphi_1
\end{array}
\right\}
\right.
$$
\item
$\mu_{\omega}(\S^1)\subseteq Z(K_{\omega}).$\\
This is also a direct computation.

\item
 $Z_{\U_{A,\Q}(\R)}(I_{\omega })=Z_{\U_{A,\Q}(\R)}(\mu_{\omega }(\S^{1}))
=K_{\omega}.$ \\
The first equality is obvious and the second is a direct computation.

\smallskip
Let
$T_{\omega}NS_{A}(\C)$ be the tangent space to
$NS_{A}(\C)$
at the point
$\omega.$
Each element
$g\in K_{\omega}$
defines a linear operator on
$T_{\omega}NS_{A}(\C).$

\item
 $\mu_{\omega}(e^{i\frac{\pi}{4}})$
defines the complex structure
(i.e. multiplication by $i$ )
on $T_{\omega}NS_{A}(\C).$\\
Indeed, given
$\gamma+i \delta\in NS_{A}(\C)$ it suffices to check that
$$
\lim_{\epsilon\in\R , \epsilon\to 0}
\frac{1}{\epsilon}
\left[
\left(
\begin{array}{cc}
1+\varphi_{2}^{-1}\varphi_{1}& -\varphi_{2}\\
\varphi_2 + \varphi_{1}\varphi_{2}^{-1}\varphi_{1}& 1-\varphi_{1}\varphi_{2}^{-1}
\end{array}
\right)
\left(
\varphi_1 +i\varphi_2 +\epsilon(\gamma+i\delta)
\right)
-(\varphi_1+i\varphi_2)
\right]
=i(\gamma+i\delta)
$$
This is a straightforward computation, which we omit.

\item
If $A$ is an abelian variety and $\omega\in C_{A},$
then $Ad_{I_{\omega}}$ is a Cartan involution in
$\U_{A,\Q}(\R),$
i.e.

$\U_{A,\Q}(\R) ^{Ad_{I_{\omega}}}$ is a maximal compact subgroup $K_{\omega}$
of
$\U_{A,\Q}(\R).$\\
Indeed, this follows from properties 2, 3, 4 above and from
Theorem 8.2 c).

\item
If $A$ is an abelian variety, then
 $\mu_{\omega}(\S^1)\subset \overline{\U (A)}(\R)$ (7.3).\\
Let us prove this.
If
${\omega}' =i\varphi_2,$ then
$
\mu_{\omega}=g_{\varphi_1}\mu_{\omega'}g^{-1}_{\varphi_1}
$
where
$g_{\varphi_1}=
\left(
\begin{array}{ll}
1&0\\
\varphi_1&1
\end{array}
\right)
\in G_{\varphi _1}
$ (7.3).
So we may assume that
$\omega=i\varphi$ and hence
$$
\mu_{\omega}(e^{i\theta})=
\cos\theta
\left(
\begin{array}{ll}
1&0\\
0&1
\end{array}
\right)
+\sin\theta
\left(
\begin{array}{cc}
0&-\varphi^{-1}\\
\varphi&0
\end{array}
\right).
$$
But then $\mu _{\varphi }(\S ^1)\subset G_{\varphi }.$ (7.3). $\Box$

\end{enumerate}

\sec{Mirror symmetry.}

Here we introduce the notion of mirror symmetry for abelian varieties and
complex tori.


\medskip
\noindent{\bf 9.1 Definition.}
{\it An {\sf algebraic pair} (resp. a {\sf weak pair})
$(A, \omega _A)$ consists of an abelian variety
(resp. a complex torus)
$A$ and an element
$\omega_{A}\in C_{A}$ (8.1)
(resp. $\omega_{A}\in NS_{A}(\C)^{0}$ (8.4)).}

\medskip
Thus a algebraic pair is a special case of a weak pair.

Let
$(A, \omega_A )$
be a weak pair.
Recall the canonical symmetric bilinear  form
$Q_A$ (3.1)
on the lattice
$\Ga_{\vphantom{\wh{A}}A}\op \Ga _{\wh{A}}.$
The group
$Hdg_{A,\Q}(\R)\ts \U_{A,\Q}(\R)$
acts on the space
$V_{\vphantom{\wh{A}}A}\op V_{\wh{A}},$
preserving the form
$Q_{A, \R}.$ In fact $\U_{A,\Q}(\R)$ is the centralizer of
$Hdg _{A,\Q}(\R)$
in $\SO (V_{\vphantom{\wh{A}}A}\op V_{\wh{A}}, Q_{A,\R})$ (Prop. 5.2.4).
We have the elements
$$
J_{A\ts\wh{A}} \in Hdg_{\vphantom{\wh{A}}A,\Q}(\R)=Hdg_{A\ts\wh{A}, \Q}(\R),
\qquad
I_{\omega_{A}}\in \U_{A,\Q}(\R)
$$
(see 2.1, 8.4) which commute and define two complex structures  on $V_{\vphantom{\wh{A}}A}
\oplus V_{\wh{A}}.$

\medskip
\noindent{\bf 9.2 Definition.}
{\it We say that algebraic pairs  (resp. weak pairs)
$(A, \omega_A),\;(B, \omega_B)$
are {\sf mirror symmetric}
if there is an isomorphism
$$
\alpha: \Ga_{\vphantom{\wh{A}}A}\op \Ga _{\wh{A}} \stackrel{\sim}{\lto}
 \Ga_{\vphantom{\wh{A}}B}\op \Ga _{\wh{B}}
$$
which identifies bilinear forms
$Q_A$ and
$Q_B$ and satisfies the following conditions
$$
\begin{array}{llll}
\alpha_{\R}\cdot J_{A\ts\wh{A}}&=&I_{\omega_{B}} \cdot \alpha _{\R},\\
\alpha_{\R}\cdot I_{\omega_{A}}&=& J_{B\ts\wh{B}} \cdot \alpha _{\R}.
\end{array}
$$}

\medskip
\noindent{\bf 9.2.1 Remark.} Let the weak pairs $(A,\omega _A)$ and
$(B,\omega _B)$
be
mirror symmetric. We will identify the lattices
$\Ga_{\vphantom{\wh{A}}A}\op \Ga_{\wh{A}}$
and
$\Ga_{\vphantom{\wh{A}}B}\op \Ga_{\wh{B}}$
by means of
$\alpha.$
That is we have
$$
\Ga_{\vphantom{\wh{A}}A}\op \Ga _{\wh{A}} =
\La
=\Ga_{\vphantom{\wh{A}}B} \op \Ga _{\wh{B}}
$$
with the bilinear form
$Q_A = Q = Q_B.$

Then the $\Q$-algebraic groups $Hdg_{A,\Q}\times \U_{A,\Q}$ and
$Hdg _{B,\Q}\times \U_{B,\Q}$ are subgroups of $\SO(\La _{\Q},Q_{\Q}).$
By the definition of the Hodge group we have the following group
inclusions
$$Hdg_{A,\Q}\subseteq \U_{B,\Q},
\qquad
Hdg_{B,\Q}\subseteq \U_{A,\Q}.$$
If $A$ and $B$ are abelian varieties, then stronger inclusions hold (see
8.4.1 (7))
$$Hdg _{A,\Q}\subseteq \overline{\U(B)},
\qquad
Hdg _{B,\Q}\subseteq \overline{\U(A)}.$$

\medskip
\noindent{\bf 9.2.2 Remark.} Let $(A,\omega _A)$ be a weak pair.
Assume that weak pairs $(B,\omega _B)$ and $(C,\omega _C)$ are
mirror symmetric to $(A,\omega _A).$ If $B$ is an abelian variety, then
so is $C.$ Indeed, composing the $\alpha $'s for $B$ and $C$ we obtain
an isomorphism of complex tori
$$B\times \wh{B}\cong C\times \wh{C}.$$
Thus $C$ is an abelian variety (as a subtorus of an abelian variety).

\medskip
\noindent{\bf 9.2.3 Remark.} Let $(A,\omega _A)$ be a weak pair. The
collection of isomorphism classes of abelian varieties $B$ for which there
exists $\omega _B \in NS_B(\C)^0$ such that the pairs $(A,\omega _A)$ and
$(B,\omega _B)$ are mirror symmetric, is finite. Indeed, fix one such $B.$
Then for any other such $C$ there exists an isomorphism of abelian
varieties
$$B\times \wh{B}\cong C\times \wh{C},$$
hence $C$ is isomorphic to an abelian subvariety of $B\times \wh{B}.$ The
finiteness follows from the following theorem of Lenstra, Oort and Zarhin.

\medskip
\noindent{\bf Theorem.}(\cite{LOZ}) {\it Let $A$ be an abelian variety. There are finitely
many isomorphism classes of abelian varieties which admit an embedding
as an abelian subvariety in $A.$}

\medskip
\noindent{\bf 9.2.4 Proposition.} {\it Let $(A,\omega _A)$ and $(B,\omega _B)$
be two mirror symmetric weak pairs. Then the set of $\alpha$'s that establish
a mirror symmetry between the pairs is a torsor over the stabilizer of
$\omega _A$ in $\U(A).$ In particular, if the pair $(A,\omega_{A})$ is algebraic
then this set is finite. }

\smallskip
\pr
Let $\alpha _1,\ \alpha _2$ be two isomorphisms that establish a mirror
symmetry of $(A,\omega _A)$ and $(B,\omega _B).$  Then $\gamma :=
\alpha _2^{-1}\alpha _1$ is an automorphism of the torus $A\times
\wh{A}$ which preserves the form $Q_A.$ Thus by Proposition 4.3.2
$\gamma \in \U(A).$ Also $\gamma I_{\omega_{A}}\gamma ^{-1}=I_{\omega_{A}}.$
Therefore by property 4 in 8.4.1 $\gamma \omega _A=\omega _A.$ This
proves the first assertion. If $\omega _A\in C_A$ then its stabilizer
in $\U_{A,\Q}(\R)$ is compact (Theorem 8.2 c). This implies the second
assertion since the group $\U(A)$ is discrete.
$\Box$

\medskip
\noindent{\bf 9.2.5 Proposition.} {\it Two weak pairs $(B,\omega _1)$
and $(B,\omega _2)$ are both mirror symmetric to the same
weak pair $(A,\omega _A)$ if and only if
there exists $g\in \U(B)$ such that $g\omega _1 =\omega _2.$}

\smallskip
\pr
$\Rightarrow$
Similarly to the previous proof we find $g\in \U(B)$ such that
$gI_{\omega _1}g^{-1}=I_{\omega _2}.$ This means that $g\omega _1=\omega _2.$

$\Leftarrow$
If $\alpha$ establishes a mirror symmetry of $(A, \omega_A)$ and $(B, \omega_1),$ then $g\alpha$ induces a mirror symmetry
of $(A, \omega_A)$ and $(B, \omega_2).$
$\Box$

\medskip
\noindent{\bf 9.2.6 Proposition.} {\it Let $A,\ B,\ C$ be abelian varieties.

a) If some algebraic pairs $(B,\omega _B),$ $(C,\omega _C)$ are both
mirror symmetric to an algebraic pair $(A,\omega _A),$ then the derived
categories $D^b(B)$ and $D^b(C)$ are equivalent.

b) If some algebraic pairs $(B,\omega _B)$ and $(A,\omega _A)$ are mirror
symmetric and the derived categories $D^b(B)$ and $D^b(C)$ are equivalent,
then there exists $\omega _C\in C_C$ such that the pairs $(C,\omega _C)$
and $(A,\omega _A)$ are mirror symmetric.}

\smallskip
\pr
This follows immediately from the Theorem 4.3.10 and Proposition 9.4.3
below.
$\Box$

\medskip
\noindent{\bf 9.3.}
Assume that weak pairs
$(A, \omega_A),\; (B, \omega_B)$
are mirror symmetric. Use the identification of 9.2.1.
Consider the Clifford algebra $Cl(\La ,Q)$ (3.1). The total cohomology
groups $H^*(A,\Z)$ and $H^*(B,\Z)$ are naturally $Cl(\La ,Q)$--modules (3.7).
There exists a unique (up to $\pm 1$) isomorphism of these
$Cl(\La ,Q)$--modules
$$\beta :H^*(A,\Z)\stackrel{\sim}{\lto} H^*(B,\Z).
$$
This isomorphism either preserves the even and odd cohomology groups
or interchanges them (3.7). In the first case we called $\beta $ {\it even}
and in the second -- {\it odd}.

\medskip
\noindent{\bf 9.3.1.} Assume that the weak pair $(A,\omega _A)$ is mirror
symmetric to yet another weak pair $(C,\omega _C)$ by an isomorphism
$$\alpha ':\Ga_{\vphantom{\wh{A}}A}\oplus \Ga_{\wh{A}}\stackrel{\sim}{\lto}
\Ga_{\vphantom{\wh{C}}C}\oplus \Ga_{\wh{C}}.$$
Identify  $\Ga_{\vphantom{\wh{C}}C}\oplus \Ga_{\wh{C}}$ with $\Ga_{\vphantom{\wh{B}}B}\oplus \Ga_{\wh{B}}$ by
means of $\alpha $ and $\alpha '.$ We claim that $\dim(\Ga _{B,\Q}\cap
\Ga _{C,\Q})$ is even and hence $\Ga _{B,\Q}\sim \Ga _{C,\Q}$ in the sense
of 3.6.1. Indeed, both $V_B$ and $V_C$ are complex subspaces of $\La _{\R}$
with respect to the complex structure $I_{\omega_{A}}$ (8.4). Hence
$\dim _{\R}(V_B\cap V_C)$ is even.

\medskip
\noindent{\bf 9.3.2 Proposition.} {\it Let weak pairs $(A, \omega _A)$
and $(B,\omega _B)$ be mirror symmetric. Let $n=dimA=dimB.$ Then the
parity of $\beta $ is equal to the parity of $n.$}

\smallskip
\pr Recall that $\La ^{\cdot }V_A^*=H^{*}(A,\R)$ and
$\La ^{\cdot }V_B^*=H^{*}(B,\R)$ are both spinorial representations of
the Lie algebra ${\frak so}(\La _{\R},Q_{\R}),$ and $\beta _{\R}
:\La ^{\cdot }V_A^*\stackrel{\sim}{\lto}\La ^{\cdot }V_B^*$ is an
${\frak so}(\La _{\R},Q_{\R})$-morphism. Thus the parity of
$\beta _{\R}$ will not change if we replace $V_B$ by a subspace
$gV_B$ for some $g\in \SO(\La _{\R},Q_{\R}).$ Let $\omega _A=\varphi _1+
i\varphi _2,$ where $\varphi _1\in NS_A(\R),$ $\varphi _2\in NS_A(\R)^0.$
Choose $g=g^{-1}_{\varphi _1}\in G_{\varphi _1}\subset \SO(\La _{\R},
Q_{\R})$ as in 8.4.1 (7). Then $gV_B$ is a subspace of $\La _{\R}$ which
is preserved by the complex structure
$$I_{i\varphi _2}=\left(\begin{array}{cc}
                        o & -\varphi _2^{-1}\\
                        \varphi _2 & 0
                        \end{array}\right).$$
It suffices to find a maximal isotropic subspace $V\subset \La _{\R}$
which is preserved by $I_{i\varphi _2}$ (hence $V\sim gV_B$) and such that
$\dim _{\R}(V\cap V_A)=n.$ Let us construct such a subspace $V.$

Consider $\varphi _2$ as a skew-symmetric form on $V_A.$ There exists a basis
$e_1,...,e_n,e_{-1},...,e_{-n}$ of $V_A$ in which
$$\varphi _2 =\left(\begin{array}{rr}
                  0 & -\Delta \\
                  \Delta & 0
                  \end{array}\right),\quad
\Delta=\left( \begin{array}{cccc}
               \delta _1 & & & \\
                    & \delta _2 & & \\
                    & & ... & \\
                    & & & \delta _n
                    \end{array}\right).$$
Let $e_1^*,...,e_n^*,e_{-1}^*,...,e_{-n}^*$ be the dual basis of
$V_{\wh{A}}.$ Then $V:=<e_1,...,e_n,e_{-1}^*,...,e_{-n}^*>$
is clearly maximal isotropic subspace of $\La$ and $I_{i\varphi _2}$-invariant. This proves the
proposition.
$\Box$

\medskip
\noindent{\bf 9.3.3.} Let $A$ be an abelian variety. Recall that the total
cohomology $H^*(A,\Q)$ is a representation space of the $\Q $--algebraic
group $Hdg_{A,\Q}\times \overline{Spin(A)}$ (7.1.2).

\medskip
\noindent{\bf Proposition.} {\it Let algebraic pairs $(A,\omega _A)$ and
$(B\omega _B)$ be mirror symmetric. As in 9.3 above consider the canonical
(up to $\pm$) isomorphism of $Cl(\La ,Q)$--modules
$$\beta :H^*(A,\Z)\stackrel{\sim}{\lto}H^*(B,\Z).$$
Then $\beta _{\Q}$ induces the inclusions of algebraic groups
$$Hdg_{A,\Q}\subseteq \overline{Spin(B)},
\qquad
Hdg_{B,\Q}\subseteq \overline{Spin(A)}.$$}

\smallskip
\pr
Let us identify $H^*(A,\Z)$ and $H^*(B,\Z)$ by means of $\beta .$ It suffices to
prove one inclusion $Hdg_{A,\Q}\subset \overline{Spin(B)}.$

We know that the action of $Hdg_{A,\Q}$ on $H^*(A,\Q)=\La ^{\cdot}
\Ga _{\wh{A},\Q}$ is induced from its action on  $\Ga _{\wh{A},\Q}$ by
functoriality. It suffices to show that this action agrees with
the spinorial representation of $Spin(\La _{Q},Q_{\Q})$ when
$Hdg_{A,\Q}$ is considered as a subgroup of $\SO(\La _{Q},Q_{\Q}).$
This is done by identifying the weights of the Lie algebra
${\frak sl}(\Ga _{A,\Q})\subset {\frak so}(\La _{\Q},Q_{\Q})$ in the
spinorial representation and the natural representation on
$\La ^{\cdot}\Ga _{\wh{A},\Q}.$ $\Box$

\medskip
\noindent{\bf 9.4.}
Let
$(A, \omega)$
be a weak pair. Let $I\in U_{A,\Q}(\R),$ $I^2=-1.$
Then $I$ and $J_{A\ts\wh{A}}\in Hdg_{A,\Q}(\R)$ are two complex structures
on the  space
$$
\La_{\R}=V_{\vphantom{\wh{A}}A}\op V_{\wh{A}}.
$$
They commute and preserve the bilinear
form $Q_{A,\R}.$
Denote by
$c$
the product $J_{A\ts\wh{A}}I.$
The operator $c$ preserves the bilinear form
$Q_{A,\R}$ and  $c^2=1.$
Hence, the bilinear form
$Q_{A,\R}(c(\cdot), \cdot)$ is symmetric.
Denote by $E_{A,I}$ the corresponding quadratic form.
In case $I=I_{\omega}$ for some $\omega \in NS_A(\C)^0$ we denote
$E_{A,I_{\omega}}=E_{A,\omega }.$

\medskip
\noindent{\bf 9.4.1 Lemma.} {\it
Let $A$ be a complex torus.
Suppose $I\in \U_{A,\Q}(\R),$ $I^2=-1$
and $I$ has the form
$
\left(
\begin{smallmatrix}
I_{11}&I_{12}\\
I_{21}&I_{22}
\end{smallmatrix}
\right).
$
The following conditions are equivalent:

\begin{tabular}{ll}
1)& $I=I_{\omega}$ for some $\omega\in NS_{A}(\C)^0.$\\
2)& $I_{12}:V_{\wh{A}}\lto V_{\vphantom{\wh{A}}A}$ is invertible.                                   \\
3)& The restriction of bilinear  form
$Q_{A,\R}(c\cdot, \cdot)$
on $V_{\wh{A}}$ is non-degenerate.
\end{tabular}}

\smallskip
\pr $1)\Rightarrow 2)$ This immediately follows from the formula (\ref{pred}) for
$I_{\omega}$ in 8.4.

\noindent
$2)\Leftrightarrow 3)$ Let $x_1, x_2 \in V_{\wh{A}}.$ We have
$$
Q_{A,\R}(c(0,x_1), (0, x_2))=x_2(J_AI_{12}(x_1)).
$$
Thus, the restriction of $Q_{A,\R}(c(\cdot),\cdot)$ to $V_{\wh{A}}$
is non-degenerate iff
 $I_{12}$ is invertible.

\noindent
$2)\Rightarrow 1)$
Let $I_{12}$ is invertible.
Since $I^2=-1$ and $I\in
\U_{A,\Q}(\R)$
we have equalities
$$
\wh{I_{12}}=I_{12},\quad \wh{I_{21}}=I_{21}, \quad \wh{I_{11}}=-I_{22}
$$
Put
$\varphi_2=-I_{12}^{-1}$
and
$\varphi_1=I_{22} I_{12}^{-1}.$
Take
$\omega=\varphi_1 + i\varphi_2.$
It is easy to see that
$I_{\omega}=I.$
$\Box$

\medskip
\noindent{\bf 9.4.2 Lemma.} {\it
Let $A$ be a complex torus
and let $I\in \U_{A,\Q}(\R),$ $I^2=-1.$
Then $I$ coincides with $I_{\omega}$ for some
$\omega\in C^{+}_A$ ( resp. $C^{-}_A$) and, consequently,
$A$ is an abelian variety iff
the quadratic form
$E_{A,I}$
is positive (resp. negative) definite.}

\smallskip
\pr
$\Rightarrow$
Suppose
$I=I_{\omega}$
with
$\omega=\varphi_1 + i\varphi_2 \in NS_A(\C)^0.$ It is not hard to check that
$I_{\omega}=g_{\varphi_1}I_{i\varphi_2}g^{-1}_{\varphi_1},$
where
$g_{\varphi_1}=
\left(
\begin{smallmatrix}
1&0\\
\varphi_1&1
\end{smallmatrix}
\right).
$
Therefore, we have an equality
$$
E_{A,I}(\lambda)=Q_{A,\R}(c\lambda, \lambda)=
Q_{A,\R}(I_{i\varphi_2}J_{A}
(\lambda^{'}), \lambda^{'}),
$$
where
$\lambda^{'}=g^{-1}_{\varphi_1}(\lambda).$

Taking $\lambda'=(l, x)$ with $l\in V_A$ and $x\in V_{\wh{A}},$ we obtain
$$
E_{A,I}(\lambda)=
Q_{A,\R}((-\varphi_2^{-1}(J_{\wh{A}}x), \varphi_2(J_{A}l)),(l,x))=
\varphi_2(J_{A}l, l)
-x(\varphi_2^{-1}(J_{\wh{A}}x))=
\varphi_2(J_{A}l, l)-x(J_{A}\varphi_2^{-1}(x)).
$$
Denoting $\varphi_2^{-1}(x)$ by $m$ in  the last expression, we get
$$
E_{A,I}(\lambda)=
-\varphi_2(m, J_{A}m)+\varphi_2(J_{A}l, l)=
\varphi_2(J_{A}m, m)+\varphi_2(J_{A}l, l)
$$
This shows that the quadratic form $E_{A,J}$ is positive (resp. negative)
definite on $\La _{A,\R}$ iff
the form $\varphi_2(J_{A}\cdot ,\cdot)$ is positive  (resp. negative)
definite on $V_A,$ which is equivalent
to $\omega \in C^{+}_{A}$ (resp. $\omega \in C^-_A$).

$\Leftarrow$
Assume that the form
$E_{A,I}$
is definite.
Let
$I=
\left(
\begin{smallmatrix}
I_{11}&I_{12}\\
I_{21}&I_{22}
\end{smallmatrix}
\right).
$
 If
$I_{12}$ is degenerate , then  there is $0\neq x\in V_{\wh{A}}$ such that
$I_{12}(x)=0.$
Hence
$$
E_{A,I}((0, x))=Q_{A,\R}((J_{\vphantom{\wh{A}}A}I_{12}(x), J_{\wh{A}}I_{22}(x)), (0,x))=
Q_{A,\R}((0,x'),(0, x))=0
$$
This contradicts the definitness of
$E_{A,I}.$
By Lemma 9.4.1, as
$I_{12}$ is invertible,
$I=I_{\omega}$ for some $\omega \in NS_A(\C)^0.$
Then the previous argument shows that
$\omega$
belongs to
$C_{A}.$ Lemma is proved.
$\Box$

\medskip
\noindent{\bf 9.4.3 Proposition.} {\it
Let weak pairs
$(A, \omega_{A}),\;(B, \omega_{B})$
be mirror symmetric.
Suppose that the weak pair $(A, \omega_A)$ is in fact an algebraic pair.
Then
$(B, \omega_{B})$ is an algebraic pair too.}

\smallskip
\pr
Since
$(A, \omega_{A})$ is an algebraic pair we have $\omega_{A}\in C_{A}.$
Hence, by Lemma 9.4.2,
the quadratic form $E_{A,\omega_{A}}$ (9.4)  is definite.
Since the weak pairs
$(A, \omega_{A})$ and $(B, \omega_{B})$ are mirror symmetric
the quadratic forms $E_{A,\omega _A}$ and $E_{B,\omega _B}$ coincide
under an identification $\alpha.$
Therefore, the quadratic form
$E_{B,\omega _B}$ is definite.
Lemma 9.4.2
  implies that
$\omega_{B}\in C_{B}.$
Thus $B$ is an abelian variety and the  weak pair
$(B, \omega_B)$ is an algebraic pair.
$\Box$

\medskip
\noindent{\bf 9.4.4 Lemma.} {\it
Let
$(A, \omega_{A})$ be an algebraic pair and
$B$ be a complex torus.
Let
$
\alpha: \Ga _A\op \Ga _{\wh{A}} \stackrel{\sim}{\lto}
\Ga _B\op \Ga _{\wh{B}}
$
be an isomorphism which identifies the forms
$Q_{A}$ and
$Q_{B}$ and
$
\alpha_{\R}\cdot I_{\omega_{A}}= J_{B\ts\wh{B}}\cdot \alpha _{\R}.
$
 Then there exists
$\omega_B\in C_{B}$
(in particular $B$ is an abelian variety) such that
$
\alpha_{\R}\cdot J_{A\ts\wh{A}}= I_{\omega_{B}}\cdot \alpha _{\R}.
$
That is $\alpha$ establishes the mirror symmetry of
algebraic pairs
$(A, \omega_{A}),\; (B, \omega_{B}).$}

\smallskip
\pr
Since $(A,\omega _A)$ is an algebraic pair the form $E_{A,\omega _A}$ on
$V_{\vphantom{\wh{A}}A}\op V_{\wh{A}}$ is definite (9.4.2). Identify the spaces
$V_{\vphantom{\wh{A}}A} \op V_{\wh{A}}$ and $V_{\vphantom{\wh{B}}B}\op V_{\wh{B}}$ with the forms $Q_{A,\R}$
and $Q_{B,\R}$ by means of $\alpha _{\R}.$ Thus
$$I:=J_{A\ts\wh{A}}\in Hdg_{A,\Q}(\R)$$
and
$$I_{\omega_{A}}=J_{B\ts\wh{B}} \in Hdg_{B,\Q}(\R).$$

Since the $\Q$-closures of $J_{A\ts\wh{A}}$ and $I_{\omega_{A}}$
commute (5.2,4)) we find that $I\in \U_{B,\Q}(\R)$
(5.2,4)). Thus we may apply Lemma 9.4.2 to $I$ and $B.$ Note that
$\alpha _{\R}$ identifies the quadratic forms $E_{A,\omega_{A}}$ and
$E_{B,I}.$ Since the first one is definite it follows from Lemma
9.4.2 that $I =I_{\omega }$ for some $\omega =\omega _B \in C_B.$
This proves the lemma.
$\Box$

\medskip
\noindent{\bf 9.4.5 Lemma.} {\it Let
$\La=\Ga\op\Ga^*$ be a lattice with canonical symmetric bilinear form
$Q$ (3.1).  Let there be given an isotropic
decomposition
$\La=\La_1 \op \La_2.$
Suppose
$I\in \GL(\La_{\R})$ is a complex structure on
$\La_{\R},$ i.e. $I^2=-1,$ that satisfies the following assumptions
$$
\begin{array}{ll}
a)& I(\La_{i, \R})=\La_{i, \R}, \quad i=1,2\\
b)& I\in O(\La _{\R},Q_{\R})
\end{array}
$$
Then the complex tori
$$
B_1:=(\La_{1, \R}/\La_1, I_1),\quad
B_2:=(\La_{2, \R}/\La_2, I_2),
$$
where $I_i$ are the restriction of
$I$ on $\La_{i, \R},$ are dual to each other,
i.e.
$B_2\cong \wh{B_1}.$}

\smallskip
\pr
The form $Q$ induces the isomorphism
$t:\La_2\lto \La_1^*$ that takes
$\lambda_2\in \La_2$ to linear functional
$Q(\lambda_2, \cdot).$
It remains to show that for any
$v_1\in \La_{1, \R}, \; v_2\in \La_{2, \R}$ the following equality holds
$$
t(I_2 (v_2))(v_1)= t(v_2)(-I_1(v_1))
$$
(3.1).
This is straightforward using the inclusion
$J\in O(\La _{\R},Q_{\R}).$
$$
t(I_2(v_2))(v_1)=Q_{\R}(I_2(v_2), v_1)=
Q_{\R}(I_2^2(v_2), I_1(v_1))=
Q_{\R}(-v_2, I_1(v_1))=
t(v_2)(-I_1(v_1)).
$$
Lemma is proved.
$\Box$

\medskip
\noindent{\bf 9.4.6 Summary.} Let $(A,\omega _A)$ be an algebraic pair.
The last two lemmas give us a ``method'' for construction of a mirror
symmetric pair. Namely, put $\La =\Ga_{\vphantom{\wh{A}}A}\oplus \Ga _{\wh{A}}.$ It
suffices to find a $Q_A$--isotropic decomposition
$\La =\La _1\op \La _2$ such that $\La _{1\R},$ $\La _{2\R}$ are preserved
by $I_{\omega_{A}}.$ Then by Lemma 9.4.5 the tori
$B:=(\La _{1\R}/\La _1,I_{\omega_{A}})$ and $\wh{B}:=
(\La _{2\R}/\La _2,I_{\omega_{A}})$ are dual. Moreover, the natural
identification
$$\alpha =id:\Ga_{\vphantom{\wh{A}}A} \op \Ga_{\wh{A}}\stackrel{\sim}{\lto}\Ga_{\vphantom{\wh{A}}B}\op
\Ga_{\wh{B}}$$
($\Ga_{\vphantom{\wh{A}}B}=\La _1,$ $\Ga_{\wh{B}}=\La _2$) identifies the forms
$Q_A$ and $Q_B.$ Thus by Lemma 9.4.4 there exists $\omega _B\in C_B$ such
that the pairs $(A,\omega _A),$ $(B,\omega _B)$ are mirror symmetric.

\medskip
\noindent{\bf 9.5.}
It is not true that for every pair
$(A,\omega_A)$
there exists a mirror symmetric pair
$(B,\omega_B).$
The problem may occur if the group
$\U_{A, \Q}$
 is too big and
$\omega_A \in C_{A}$
is chosen too general.

\medskip
\noindent{\bf 9.5.1 Counterexample.}
Let
$E$
be an elliptic curve with complex multiplication.
We have
$e_0 =1, d=1$ (1.8)
and
$\U_{E\ts E, \Q}(\R)=\U(2,2)$ (5.3.2).
Therefore this group acts irreducibly on
$V_{\vphantom{\wh{E\ts E}}E\times E}\op V_{\wh{E\times E}}.$
Assume that for all
$\omega\in C_{E\ts E}$
there exists a mirror symmetric pair
$(B, \omega'),$ i.e.
there exists $\alpha $ as in 9.2 with the corresponding properties.
Let us use $\alpha $ to identify
$$\Ga_{\vphantom{\wh{E\ts E}}E\times E}\op \Ga _{\wh{E\times E}}=\Ga_{\vphantom{\wh{A}}B}\op \Ga _{\wh{B}}.$$
The Hodge group
$Hdg_{B, \Q}$
is the
$\Q$--closure of the compact torus
$\mu_{\omega}(\S^1)\subset
Aut (V_{\vphantom{\wh{A}}B}\op V_{\wh{B}}).$
For a general
$\omega$
this
$\Q$--closure is the group $\SU(2,2) \subset \U_{E\ts E, \Q}(\R).$
But the group $Hdg _{B,\Q}(\R)$ acts faithfully on
$V_{\vphantom{\wh{A}}B}\op V_{\wh{B}}$
and preserves each summand. This is a contradiction.

\medskip
\noindent{\bf 9.6.} However, each torus $A$ with a nonempty $NS_A^0$
has an $\omega _A\in NS_A(\C)^0$ such that the pair $(A,\omega _A)$
has a mirror symmetric pair. In particular for every abelian variety
$A$ there exists $\omega _A\in NS_A(\C)^0$ (even $\omega _A \in C_A$)
such that the pair $(A,\omega _A)$ has a mirror symmetric one. This
follows from the next proposition.

\medskip
\noindent{\bf 9.6.1 Proposition.} {\it Let $A$ be a complex torus of
dimension $n.$ Let $\varphi \in NS_A^0,$ i.e. $\varphi \in \Hom (A,\wh{A})$
is an isogeny. Let $\tau =a+ib\in\C,$ $b\neq 0.$ Consider the element
$\omega _A :=\tau \varphi \in NS_A(\C)^0$ and the weak pair $(A,\omega _A).$
Then there exist isogeneous elliptic curves $E_1,...,E_n$ and an element
$\omega _E\in NS_E(\C)^0,$ where $E=E_1\times ...\times E_n,$ such that
the weak pair $(E,\omega _E)$ is mirror symmetric to the weak pair
$(A,\omega _A).$ }

\smallskip
\pr There exists a basis $e_1,...,e_n,e_{-1},...,e_{-n}$ of $\Ga _A$
in which the bilinear form $\varphi $ has a matrix
$$\left(
  \begin{array}{rr}
   0 & \Delta \\
 -\Delta & 0
  \end{array}\right),
\quad
\mbox{where}
\quad
\Delta = \left(\begin{array}{ccc}
                \delta _1 & & 0 \\
                      & ... &   \\
                   0 &  & \delta _n
                   \end{array}\right),\quad \delta _i \in \Z , \quad
\delta _1|\delta _2|\delta _3...
$$

Let $e^*_1,...e^*_n,e^*_{-1},...,e^*_{-n}$ be the dual basis of
$\Ga _{\wh{A}}.$ Then the map $\varphi :\Ga _A\lto \Ga _{\wh{A}}$ is
$$\varphi : e_i\mapsto \delta _ie^*_{-i},\quad e_{-i}\mapsto -\delta _i
e^*_i.$$

Put
$$\Ga _i:=\Z e_i \oplus \Z e_{-i}^*,\quad \Ga _i^*:= \Z e_{-i}
\oplus \Z e_i^* ,\quad i=1,...,n.$$
Clearly the subgroups
$\La _1:=\oplus \Ga _i,$ $\La _2:=\oplus \Ga ^*_i$ are isotropic. We have
$$I_{\omega_{A}} =\left(\begin{array}{cc}
                    b^{-1}a & -b^{-1}\varphi ^{-1}\\
                    (b+ab^{-1}a)\varphi & -ab^{-1}
                     \end{array}\right).$$
This operator preserves each subspace $\Ga _{i\R},$ $\Ga _{i\R}^*.$
In particular it preserves $\La _{1,\R},$ $ \La _{2,\R}.$ Thus by Lemma
9.4.5 the tori $(\La _{1,\R}/\La _1,I_{\omega_{A}}),$
$(\La _{2,\R}/\La _2,I_{\omega_{A}})$ are dual. Moreover, it is clear that
the torus $(\La _{1,\R}/\La _1,I_{\omega_{A}})$ is a product of elliptic curves
$E_i:=(\Ga _{i\R}/\Ga _i,I_{\omega }),$ $i=1,...,n.$

For each $i$ the map
$$\Ga _1\lto \Ga _i,\quad e_1\mapsto e_i,\quad e_{-1}^*\mapsto \delta_1^{-1}
\delta _i
e_{-i}^*$$
commutes with $I_{\omega_{A}}$ and so is an isogeny of elliptic curves
$E_1$ and $E_i.$

If the weak pair $(A,\omega _A)$ were an algebraic pair, then by Lemma 9.4.4
there would exist an element $\omega _E\in C_E$ such that
the pairs $(A,\omega _A)$ and $(E,\omega _E)$ are mirror symmetric.
 Thus we are done in the algebraic case.
In general there may not exist the desired $\omega _E \in NS_E(\C)^0,$ so
it may be necessary to choose a different decomposition
$\La =\La _1\oplus \La _2.$

By Lemma 9.4.1 in order for $\omega _E$ to exist the symmetric
bilinear form
$Q_{\R}(I_{\omega_{A}}J_{A\ts\wh{A}}(\cdot),\cdot)$ must be nondegenerate on
$\La _{2,\R}.$ This is equivalent to the statement that
$J_{A\ts\wh{A}}\La _{2,\R}\cap \La _{2,\R}=0.$
Put $W:=\oplus \R e_{-i}.$
Then since $J_{A\ts\wh{A}}$ preserves the form $Q_{\R}$ the last equality is
equivalent to
$$J_{A\ts\wh{A}}W\cap W=0. \qquad \qquad \qquad (*)$$
We are free to apply elements of the symplectic group
$\Sp (\Ga _A,\varphi ;\Z)$ to the sublattice
$\oplus \Z e_{-i}$ to achieve the transversality condition (*) above.
Since this discrete group is Zariski dense in the corresponding group
of $\R$--points $\Sp (V_A, \varphi ;\R)$ and since $\Sp (V_A, \varphi ;\R)$
acts transitively on the collection of maximal $\varphi$--isotropic
subspaces of $V_A$ it suffices to prove the following lemma

\medskip
\noindent{\bf 9.6.2 Lemma.} {\it Let $V\cong \R^{2n}$ be a real vector space
with a nondegenerate symplectic form $\varphi .$ Let $J\in \End (V)$ be
a complex structure on $V,$ which preserves $\varphi .$ Then there
exists a maximal $\varphi $--isotropic subspace $W\subset V$ such that
$$JW\cap W=0.$$
}

\smallskip
\pr Note that the form $\varphi (J\cdot ,\cdot )$ on $V$ is symmetric and
nondegenerate.

We will choose elements $x_1,...,x_n\in V$ such that $W:=\oplus \R x_i$
has the desired properties by induction on $i.$ Choose $x_1$ such that
$\varphi (Jx_1,x_1)\neq 0.$ Let $V=<x_1,Jx_1>\oplus V'$ be a $\varphi $-
orthogonal decomposition. Note that $V'$ is $J$--invariant. Hence we may
replace $V$ by $V'$ and choose $x_2\in V'$ s.t. $\varphi (Jx_2,x_2)\neq 0.$
And so on. Cleary the resulting space
$W=\oplus \R x_i$ is maximal $\varphi $-isotropic and $JW\cap W=0.$
This proves the lemma and the Proposition 9.6.1.
$\Box$

\medskip
\noindent{\bf 9.6.3 Corollary.} {\it Let $A$ be an abelian variety
with $NS(A)\cong \Z.$ Then $C_A=NS_A(\C)^0$ and
for any $\omega _A\in C_A$ the algebraic
pair $(A,\omega _A)$ has a mirror symmetric one $(E,\omega _E)$ which
is like in Proposition 9.6.1 above.}

\smallskip
\pr
The first assertion is obvious. The second is a direct consequence of
Proposition 9.6.1. $\Box$

\sec{ The G--construction  and its application to mirror symmetry.}

\ssec
Let
$A$
be an abelian variety.
The ring
$R=\End(A)$
is finite over
$\Z$ such that
$D:=\End^0 (A)=R\ot_{\Z}{\Q}$
is semisimple with a positive definite anti-involution
$*:D\to D^{op}.$
The lattice
$\Ga= H_1 (A, \Z)$
carries a natural structure of a faithful
left $R$-module.

For a given
$D$
with an anti-involution
$*$
the construction by Gerritzen \cite{Ge} gives an abelian variety
$X$
with
$\End^0 (X)\cong D.$
Furthermore, it was proved in \cite{OZ} that any finite dimensional
algebra over
$\Q$
can be realized as
$\End^{0}(T)$
for some complex torus
$T.$

Now we are going to describe these constructions
with some corrections and modifications related to
the specifics  of our situation being that we work with
$R$
and not with
$D.$

\sssec
Let
$R$
be an order in
$\Q$-algebra
$D$
of  finite dimension over $\Q,$ i.e.
$D=R\ot_{\Z}{\Q}.$
Denote by
$\Ga$
a
left
$R$-module that satisfies the following conditions:
%
%
%
%
%
%
%
%
\begin{equation}\label{latt}
\begin{tabular}{ll}
1.& $\Ga$ is a lattice of even dimension, i.e.
$\Ga=\Z^{2n}$
as abelian group.\\
2.&$\Ga$ is a faithful module, i.e.
the homomorphism
$R\lto \End_{\Z}(\Ga)$
is an embedding.\\
3.& The order
$O(\Ga):=\{\; a\in D\; |\; a(\Ga)\subset \Ga\;\}$
coincides with
$R.$\\
4.&$\Ga$ is a direct sum of $R$-modules
$\Ga_1$ and $\Ga_2$
such that
there exists an isomorphism \\
&$e:M_1 \stackrel{\sim}{\to} M_2$
of $D$--modules
$M_i=\Ga_i \ot_{\Z}\Q, \; (i=1,2).$
\end{tabular}
\end{equation}
%
%
%
%
Let us fix an isomorphism $e$ and, henceforth, we will often
identify
$D$--modules
$
M_1$
and
$M_2$
with respect to
$e$ and will use
notation $M.$

Denote by
$C$
the algebra
$\End_{D}(M).$
The algebras
$D$
and
$C$
are two subalgebras of
$\End_{\Q}(M)\cong M(n, \Q).$
Moreover $C$ is a centralizer of $D$
in $M(n, \Q).$
On the other side, it is true that the centralizer of
$C$
contains
$D$
but does not necessarily coincide with
$D.$

Let
$J$
be an element of
$M(2, C\ot \R)=(\End_{D}M^{\op 2})\ot_{\Q} \R$
such that
$J^2=-Id.$
Any such
$J$
defines a complex structure on
the real space
$\Ga_{\R}=(M_1\op M_2)_{\R}$
and, as consequence,
defines a complex torus
\begin{equation}
A_{J}:=(\Ga_{\R}/\Ga, J).
\end{equation}
\lm{Lemma} For any complex torus
$A_{J},$ defined above,
the endomorphism ring
$\End(A_{J})$
contains $R$ as subring.
\elm
\pr It is clear, because
$\End(A_{J})$ coincides with maximal subring of
$\End_{\Z}(\Ga)$
that commutes with
$J.$ And each
$r\in R\subset D$
commutes with any element from
$M(2, C\ot \R).$
$\Box$

\sssec
The set of operators
$J\in M(2,C\ot \R)$
with
$J^2=-Id$
is
not empty.
Actually, any  element
$q=q_1 +iq_2\in C\ot_{\Q}\C$
with nondegenerate
$q_2$
defines the operator
\begin{equation}\label{comstr}
J(q)=
\begin{pmatrix}
q_2^{-1}q_1&-q^{-1}_{2}\\
q_{2}+q_1 q_2^{-1} q_1&-q_1 q_2^{-1}
\end{pmatrix}
=
\begin{pmatrix}
1&0\\
q_{1}&1
\end{pmatrix}
\begin{pmatrix}
0&-q^{-1}_{2}\\
q_{2}&0
\end{pmatrix}
\begin{pmatrix}
1&0\\
-q_{1}&1
\end{pmatrix}
\end{equation}
It is easy to check that
$J(q)^2=-Id.$
\lm{Proposition}\label{coinend}
Suppose that the centralizer of
$C$ in $M(n, \Q)$ coincides with $D.$
Then there exists
$q\in C\ot{\C}$
such that
$\End(A_{J(q)})=R.$
\elm
\pr
Let
$\{e_1=1,...,e_t\}$
be a basis of $\Q$-vector space
$C.$
Set
$q_1=0$
and
$q_2=\sum\limits_{i=1}^{t} r_i e_i,$
where
$r_i\in \R.$

Let us describe the algebra
$\End^0 (A_{J(q)}).$
It consists of elements of
$\End_{\Q}(M^{\op 2})$ that commute with $J.$
Each such element is a $2\ts 2$ matrix that satisfies
the condition
$$
\begin{pmatrix}
0& q_2^{-1}\\
-q_2& 0
\end{pmatrix}
\begin{pmatrix}
a& b\\
c& d
\end{pmatrix}
\begin{pmatrix}
0& -q^{-1}_2\\
q_2& 0
\end{pmatrix}
=
\begin{pmatrix}
q^{-1}_2dq_2& -q_2^{-1}cq_2^{-1}\\
-q_2bq_2& q_2aq_2^{-1}
\end{pmatrix}
=
\begin{pmatrix}
a& b\\
c& d
\end{pmatrix}
$$
This implies that
$c= -q_2bq_2$
and
$q_2a=dq_2.$

Let us take $q_2$ such that
$(r_1,...,r_t)$
are  algebraic independent over
$\Q.$
For such
$q_2$
we obtain
$c=0=b$ and
$ae_i = e_i d$
for all
$i.$
As $e_1=1$
we have
$a=d$ and
$a$ commutes with
any element of $C.$
It follows that
$a$
belongs to
centralizer of
$C$
which coincides with
$D.$
Thus for general
$q_2$
the algebra
$\End^0 (A(q))$
is isomorphic to
$D.$

Now, the algebra
$\End(A_{J(q)})$ is an order in
$\End^0 (A_{J(q)})=D$ that takes
the lattice
$\Ga$ to itself.
By condition 3) of (\ref{latt})
it coincides with
$R.$
$\Box$
\ex{Remark}
Note that if
$\Ga_{i}\cong R$
as left $R$-modules  for $i=1,2$
then the lattice
$\Ga=\Ga_{1}\op\Ga_{2}$
satisfies
the condition (\ref{latt}). In this case
$C\cong D^{op}$ and,
in addition, the centralizer of
$C$
coincides with
$D.$
This way  by Proposition \ref{coinend} for any
$R$
there exists a complex torus
with $R$ as the endomorphism ring.
\eex

\ssec
Let
$\Ga$
be a lattice as in (\ref{latt}).
Suppose that, in addition, $\Ga$ satisfies the following extra condition:
\begin{equation}\label{extra}
\begin{tabular}{ll}
5.& There is a rational symmetric form
$g:M\ts M\to\Q$ such that
the algebra $D$ \\
& is invariant with respect to an anti-involution $\sim$ on $\End_{\Q}(M),$ defined by the rule\\ &$\wt{\alpha}=g^{-1}\alpha^{t}g.$
\end{tabular}
\end{equation}
Consider
a rational skew-symmetric form
$s$ on
$\Ga_{\Q}$ given as
\begin{equation}\label{skew}
s\langle (m_1, m_2), (n_1, n_2)\rangle=
g(m_2, e(n_1))-g(e(m_1), n_2)
\end{equation}
Both subspaces
$M_1=\Ga_{1}\ot \Q$ and $M_2=\Ga_2\ot\Q$
are isotropic with respect to $s.$
Since $C$ is a centralizer of $D$ in
$\End_{\Q}(M)$  it
is invariant
with respect to
$\sim$ too.
Denote by
$S(D)$ and
$S(C)$
the subspaces of symmetric elements
$\wt{\alpha}=\alpha$
of
$D$
and
$C$
respectively.

Take an element
$q=q_1+iq_2\in S(C)\ot\C$
with nondegenerate
$q_2$
and consider the complex structure
$J(q)$ defined by the
formula (\ref{comstr}).
Since
$q$ is symmetric element a computation shows that the form
$s_{\R}$ is $J(q)$-invariant.
A some multiple of $s$ is integral.
Hence it
determines a line bundle on the complex torus
$A_{J(q)}=(\Ga_{\R}/\Ga, J(q)).$
Further, any element
$d\in S(D)$ defines the rational form
$$
s\langle(dm_1, d m_2),(n_1, n_2)\rangle
$$
that is skew-symmetric. This correspondence gives the embedding
$i_{s}:S(D)\hookrightarrow NS_{A}(\Q).$
\lm{Definition}
We say  that a weak pair $(A_{J(q)},\omega)$ is obtained by
{\bf{G}}--construction if the complex torus $A_{J(q)}$ is
$(\Ga_{\R}/\Ga, J(q)),$
where $\Ga$ satisfies contions (\ref{latt}) and (\ref{extra}), $q\in S(C)\ot \C,$ and the form $\omega$ belongs to
$i_{s}(S(D))\ot \C.$
\elm

Now we show that under some assumption we can guarantee
the existence of such a form $s$ and moreover the complex torus
$A_{J(q)}$
would be an abelian variety.

\lm{Lemma}\label{antiin}
Let the algebra
$D$ be  semi-simple  with
a positive definite anti-involution
$*:D\to D^{op}.$
Then for any $D$-module $M$ there exists a
positive definite symmetric
form $g: M\ts M\lto \Q$
such that the anti-involution
$\sim$ on $\End_{\Q}(M),$ defined by the rule
$\wt{a}=g^{-1}a^t g,$
under the  restriction on
$D$
coincides with
$*.$
\elm
\pr
Actually, since
$D$
is semi-simple
each module is projective.
Hence
$M$
can be considered as submodule
of left free module
$D^{\op k}$
for some natural number $k.$
Let us define a bilinear form
$u: D^{\op k}\ts D^{\op k}\lto \Q$
by formula:
$$
u(x,y):= \sum_{i=1}^{k} \left[Tr(x_i^* y_i)+Tr(y_i^* x_i)\right]
\quad\mbox{where}\quad x=(x_1,...,x_k), y=(y_1,...,y_k)\in D^{\op k}
$$
The bilinear form
$u$
is symmetric and positive definite.
Denote by
$g$
the restriction of
$u$ on
$M.$
It is a positive definite symmetric bilinear form on
$M.$
For any endomorphism
$a\in\End_{\Q}(M)$
there exists an adjoint with respect to
$g$
endomorphism
$\wt{\alpha}$ given by rule
$$
g(m, a(n))= g(\wt{a}(m), n),\quad\mbox{i.e.}
\quad \wt{a}=g^{-1}a^{t}g
$$
It is clear that
$\sim$
is anti-involution
on
$\End_{\Q}(W)$
and
$\wt{d}=d^*$
for any
$d\in D.$
$\Box$
\lm{Proposition}
Let  a lattice $\Ga$ be as in (\ref{latt}) and
the algebra $D$ be a semi-simple with a positive definite anti-involution $*.$ Then
there exists a neighbourhood $W\subset S(C)\ot\C$ of the point $i$ such that
for any $q\in W$ the complex torus
$A_{J(q)}$ is an abelian variety.

If, in addition, the algebra
$C$ is symmetrically generated (i.e. is generated by the subspace
$S(C)$).
Then there exists $q\in W$
such that
$\End(A_{J(q)})= R.$
\elm
\pr
By Lemma \ref{antiin} there is a positive definite symmetric form $g$ on $M.$
Consider the skew-symmetric form $s$ defined as in (\ref{skew}). The form $s_{\R}$ on
$M_{\R}^{\op 2}$
is invariant under the action of $J(q),$
and, consequently, a some multiple of $s$ gives a line bundle
on
$A_{J(q)}.$
If this line bundle is ample, then
$A_{J(q)}$
is an abelian variety.
The ampleness is equivalent to
positive definiteness of the symmetric form
$s_{\R}\langle J(q)(\cdot), \cdot\rangle.$
For
$q_1=0$ and $q_2=1$
we have
$$
s_{\R}\langle J(q)(m_1 , m_2), (n_1 , n_2)\rangle=
s_{\R}\langle (-m_2, m_1), (n_1, n_2)\rangle=
g(m_1, n_1) + g(m_2, n_2)
$$
Hence the positive definiteness  of $g$ implies  the positive definiteness of
the form
$s_{\R}\langle J(q)(\cdot),\cdot\rangle$
for $q=i$ and, consequently, for any $q$ from some neighbourhood
$W$
of $i.$

Further, since $D$ is a semi-simple algebra
the cenralizer of $C$ coincides with $D.$
By assumption $C$ is  generated by $S(C)$
and , consequently, $D$ is the centralizer of $S(C).$
Now the existence of $q$ for which $\End A_{J(q)}=R$
is proved  the same way as in Proposition \ref{coinend}.
$\Box$

\ex{Remark} Note that if the endomorphism ring
$\End A_{J(q)}$  coincides with $R$ then the embedding
$i_{s}: S(D)\hookrightarrow NS_{A_J}(\Q)$ is an isomorphism.
\eex

\ssec
Now we are going to give an alternative description of weak pairs
which are obtained by the {\bf{G}}--construction.
\lm{Definition}\label{wb}
We say that a weak pair $(A, \omega)$ is {\sf well-becoming},
if there is a decomposition
$\Ga=H_{1}(A, \Z)=\Ga_{1}\op\Ga_{2}$
such that the following conditions hold:
\begin{enumerate}
\item The subspaces $\Ga_{k\R} (k=1,2)$ are isotropic with
respect to $\omega=\varphi_1 +i\varphi_2,$ i.e.
$\omega\mid_{\Ga_{k\R}}\equiv 0.$
\item Maps $J_{21}:\Ga_{1\R}\to\Ga_{2\R}$ and
$J_{12}:\Ga_{2\R}\to\Ga_{1\R},$ which are the components of the complex
structure
$J=
\left(
\begin{smallmatrix}
J_{11}& J_{12}\\
J_{21}& J_{22}
\end{smallmatrix}
\right)$ on $A,$
 are invertible.
\end{enumerate}
\elm
\ex{Remark}
The second condition of the previous definition  is equivalent
to saying that the restrictions of holomorphic $n$-form to
the tori $\Ga_{k\R}/\Ga_{k}$ are non--zero.
Moreover, these tori are special Lagragian submanifolds of $A.$
\eex
The following technical results are needed in the sequel.
\lm{Definition}
Let $M$ be a vector space over $\Q$ and let $r_{1},\dots,r_{k}$ be elements
of real space $M_{\R}.$ The minimal vector subspace $W\subset M$
sucht that  $r_{1},\dots, r_{k}\in W_{\R}$ is called
$\Q$--envelope of the set $\{r_{1},\dots,r_{k}\}.$
\elm
\lm{Lemma} \label{tech}
Let $(A, \omega)$ be a weak pair.
Suppose there is a decomposition
$\Ga=H_{1}(A, \Z)=\Ga_{1}\op\Ga_{2}$ such that
the subspaces $\Ga_{k\R}$ are $\omega$--isotropic and
the map $J_{12}:\Ga_{2\R}\lto\Ga_{1\R}$ is invertible.
Then the pair $(A, \omega)$ is well-becoming.
\elm
\pr
Consider the subspace $U\in \Hom_{\Q}(\Ga_{1\Q}, \Ga_{2\Q})$
consisting of such maps that satisfy the condition
$
s'\langle l, e(m)\rangle=s'\langle m, e(l)\rangle
$
with
$
l,m\in \Ga_{1\Q}
$
for any
$s'\in NS_{A}\ot\Q.$
One can see that $J_{12}^{-1}$ belongs to $U_{\R}.$
Hence, there is
invertible $e\in U.$
Fix one such $e.$

Take a sublattice
$\Ga_1'=(\gamma_1,Ne(\gamma_1))$
for sufficient large $N,$
where
$\gamma_1\in\Ga_1.$
Consider the decomposition
$\Ga=\Ga_1'\op\Ga_2.$
It easy to check that
$\Ga_{1\R}'$ is isotropic with respect to
$\omega.$
In addition, maps $J'_{21}:\Ga'_{1\R}\to\Ga_{2\R}$ and
$J'_{12}:\Ga_{2\R}\to\Ga'_{1\R}$ will be invertible.
Hence the pair $(A_{J(q)}, \omega)$ is well-becoming.
$\Box$
\lm{Proposition}\label{equival}
The weak pair
$(A, \omega)$
is well-becoming if and only if
it can be obtained by the {\bf{G}}--construction.
\elm
\pr
$\Leftarrow$
Suppose that a weak pair
$(A_{J(q)}, \omega)$
is obtained by
the {\bf{G}}--construction.
The lattice
$H_{1}(A_{J(q)}, \Z)=\Ga=\Ga_1\op\Ga_2$
satisfies
the conditions (\ref{latt}).
Moreover, there is a skew-symmetric form
$s$
with isotropic sublattices
$\Ga_i$ for $i=1,2.$
This imples that
$\Ga_{i\Q}$
are isotropic with respect
to any element from
$i_{s}(S(D)).$
As
$\omega$
belongs to
$i_s(S(D))\ot\C,$
the real spaces
$\Ga_{i\R}$
are isotropic with respect to
$\omega.$
Further, it follows from the formula (\ref{comstr}) for
$J(q)$
that the map
$J_{21}$
is invertible.
Hence, by Lemma \ref{tech} the pair $(A_{J(q)}, \omega)$ is well-becoming.

$\Rightarrow$
Set $\omega=\varphi_1 + i\varphi_2.$
Denote by
$W$
the
$\Q$-envelope of
$\{\varphi_{1}, \varphi_{2}\}.$
Let
$D\subset \End(\Ga_{\Q})$
be a
$\Q$-algebra that is generated by all elements
of the form
$s_{2}^{-1}s_{1}$
with
$s_{1}, s_{2}\in W.$
Any element
$s_{2}^{-1}s_{1}$
sends the space
$\Ga_{i\Q}$
to itself, because
$\Ga_{i\Q}$
are maximal isotropic
with respect to any non-degenerate
$s_{2}\in W.$
Therefore the rational spaces
$\Ga_{i\Q}$
are  modules over
$D.$
Take some integral element
$s\in W\bigcap NS_{A}^0.$
Define an anti-involution
$\sim$ on $D$ by rule
$\wt{d}=s^{-1}d^{t}s.$
It is clear that an element
$s^{-1}s'$
belongs
to the subspace of symmetric elements
 $S(D)$ for any $s'.$
As in previous lemma let us consider the subspace
$U\in \Hom_{\Q}(\Ga_{1\Q}, \Ga_{2\Q})$
consisting of such maps that satisfy the condition
$
s'\langle l, e(m)\rangle=s'\langle l, e(m)\rangle
$
with
$
l,m\in \Ga_{1\Q}
$
for any
$s'\in NS_{A}\ot\Q.$
Since $J_{12}^{-1}$ belongs to $U_{\R},$
there
is an invertible $e\in U.$
It follows from definition of $D$ that
$e$ gives an isomorphism between $D$-modules $\Ga_{1\Q}$
and $\Ga_{2\Q}.$

Let
$R$ be a maximal order in
$D$ that takes
$\Ga$ to itself.
It follows from above that
$\Ga=\Ga_{1}\op\Ga_{2}$ satisfies the
conditions (\ref{latt}).
By construction,
$\Ga_{1}$
and
$\Ga_{2}$
are isotropic by means of
$s.$
Identifying $\Ga_{1\Q}$ and
$\Ga_{2\Q}$ with respect to $e,$
we can consider
$J_{A}$
as element of
$M(2, C\ot\R),$ where $C$ is a centralizer of $D.$
In addition, $J_{A}$ preserves the form
$s.$
Let it has a form
$
\left(
\begin{smallmatrix}
J_{11}&J_{12}\\
J_{21}&J_{22}
\end{smallmatrix}
\right)
$
with respect to decomposition $\Ga=\Ga_{1}\op\Ga_{2}.$
As in the proof of Lemma 9.4.1, we have
$$
\wt{J_{21}}=J_{21},\quad \wt{J_{12}}=J_{12}, \quad \wt{J_{11}}=-J_{22}.
$$
Put
$q_2=-e^{-1}J_{12}^{-1}$
and
$q_1=e^{-1}J_{22} J_{12}^{-1}.$
It is easy to see that
$q=q_1+iq_2$ belongs to $S(C)\ot\C$ and
$J_{A}=J(q),$ where
$J(q)$ is defined by the formula (\ref{comstr}).

Finally, any $\omega'\in W_{\C}$ belongs to
$i_{s}(S(D))\ot\C\subseteq NS_{A}(\C)$ because $W\subseteq i_{s}(S(D)).$
Thus $(A, \omega')$ is obtained by
{\bf{G}}--consruction. Proposition is proved.
$\Box$

We claim that any well-becoming pair has a mirror symmetric pair.
First, let us give a construction.
\ssec {\bf Construction.}
Let $(A_{J}, \omega)$ be a well-becoming pair.
As above by
$\La$
denote
$H_1(A_{J}\ts\wh{A_{J}}, \Z).$
There is a decomposition
$$
\La=\Ga\op\Ga^{*}= \Ga_1\op\Ga_2\op\Ga_1^*\op\Ga_2^*.
$$
The form $\omega$ defines an element
$I_{\omega}\in SO_{\Q}(Q,\R)$ by formula (\ref{pred}),
where $Q$ is canonical symmetric bilinear form on $\La$ as in (3.1).
By definition
the subspaces
$\Ga_{k, \R}$
are isotropic with respect to a
$\omega.$
This implies that
$(\Ga_1^{*}\op \Ga_2)_{\R}$ and
$(\Ga_1\op \Ga_2^*)_{\R}$ are
the
$I_{\omega}$--invariant subspaces.
By
$I$
and
$I'$
denote the restriction of the linear operator
$I_{\omega}$
on the subspaces
$(\Ga_1^{*}\op \Ga_2)_{\R}$
and
$(\Ga_1\op \Ga_2^*)_{\R}$
respectively.
Denote the sublattice
$(\Ga_1^{*}\op \Ga_2)$ by $\Si.$
Define a complex torus $B_{I}$ by the rule
$B_{I}=(\Si_{\R}/\Si, I).$
\lm{Theorem}\label{mirsym}
Any well-becoming pair $(A_{J}, \omega)$
has a mirror symmetric pair
$(B_{I}, \theta),$ where the complex torus
$B_{I}$ is constructed above.
Moreover, the pair $(B_{I}, \theta)$ is well-becoming too.
\elm
\pr
First, take $\Si^{*}=\Ga_{1}\op \Ga_{2}^*$ and consider the complex torus
$$
C_{I'}=(\Si^*_{\R}/\Si^*, I')
$$
By Lemma 9.4.5, the torus
$C_{I'}$
is isomorphic to
$\wh{B_{I}}.$
Therefore, there is a canonical identification of lattices
$H_1(B_{I}\ts\wh{B_{I}}, \Z)$ and $\La.$
Moreover, by construction, the element $J_{B\ts\wh{B}}$ coincides
with
$I_{\omega}.$
Hence, to define mirror symmetric pair we have to find
$\theta\in NS_{B_{I}}(\C)^0$
such that
$I_{\theta}=J_{A\ts\wh{A}}.$

Denote by
$
\left(
\begin{smallmatrix}
\alpha&\beta\\
\gamma&\delta
\end{smallmatrix}
\right)
$
a matrix of $J_{A\ts\wh{A}}$
with respect to the decomposition
$\La=\Si\op\Si^*.$
Since
$J_{A}=
\left(
\begin{smallmatrix}
J_{11}&J_{12}\\
J_{21}&J_{22}
\end{smallmatrix}
\right)
$
the map
$\beta$
sends
$(a,b)\in \Si^{*}_{\R}=(\Ga_1\op \Ga_2^*)_{\R}$
to
$(-J_{21}^{*}(b),J_{21}(a))\in \Si_{\R}=(\Ga_1^*\op \Ga_2^*)_{\R}.$
As
$J_{21}$ is invertible  $\beta$ is invertible too.
Hence, by Lemma 9.4.1,
there is
$\theta=\psi_1+i\psi_{2}\in NS_{B_I}(\C)^0$
such that
$I_{\theta}=J_{A\ts\wh{A}}.$
Thus, there is a mirror symmetric pair $(B_{I}, \theta).$

Further, since
$\psi_2=-\beta^{-1}$
and
$\psi_1=\delta \beta^{-1}$
we have that $\psi_{k}(\Ga_{2\R})=\Ga_{1\R}$
and
$\psi_{k}(\Ga_{1\R}^*)=\Ga_{2\R}^*,$ for $k=1,2.$
Hence the subspaces
$\Ga_{1\R},
\Ga_{2\R}^*$ are isotropic with respect to
$\theta.$
In addition, the map $I_{21}:\Ga_{1\R}^*\to \Ga_{2\R}$ coincides
with the restriction on $\Ga_{1\R}$ of the map
$-\varphi_{2}^{-1}.$ Hence $I_{21}$ is invertible and
the pair $(B_{I}, \theta)$ is well-becoming by Lemma 10.11.
This completes the proof.
$\Box$
\lm{Corollary}
Let $A_{J(q)}=(\Ga_{\R}/\Ga, J(q))$ be a complex torus
obtained by {\bf{G}}--constriction with
$R=\End A_{J(q)}.$
Then for any
$\omega=\varphi_1+i\varphi_2\in NS_{A_{J(q)}}(\C)$
with non-degenerate
$\varphi_2$
the pair
$(A_{J(q)}, \omega)$ has a mirror symmetric pair.
\elm
\pr Since
$\End A_{J(q)}=R$
the map
$i_s: S(D)\lto NS_{A_{J(q)}}(\Q)$
is an isomorphism.
Hence, by Lemma \ref{equival}, any weak pair
$(A_{J(q)}, \omega)$
is well-becoming. By Theorem \ref{mirsym}, it has a mirror symmetric
pair
$(B_{I}, \theta).$
$\Box$
\lm{Corollary} For any abelian variety $A$ there exists
an element $\omega\in C_{A}$ such that the pair
$(A, \omega)$
has a mirror symmetric pair.
\elm
\pr
By Theorem \ref{mirsym}, it is sufficient to show that there is
$\omega$ such that the pair $(A, \omega)$ is well-becoming.
Since $A$ is an abelian variety there is a non-degenerate
 skew-symmetric
form $s$ on $\Ga=H_1(A, \Z)$ such that
the form $s\langle J(\cdot), \cdot\rangle$ is positive definite.
There exists   isotropic with respect to
$s$
decomposition  $\Ga=\Ga_1\op\Ga_2.$
Since the form $s\langle J(\cdot), \cdot\rangle$
is positive definite the linear maps $J_{12}$ and $J_{21}$ are invertible.
Thus, the pair $(A, \omega)$ with
$\omega=c s$ is well-becoming for any $c\in\C^{*}.$
Therefore, it has a mirror symmetric pair.
$\Box$

Notice that this corollary is a particular case of (9.6.1).
\th{Theorem}
Let $A$ be an abelian variety of dimension $n$ and
$\omega=\varphi_{1}+i\varphi_{2}\in C_{A}.$
Let $W\subset NS_{A}(\Q)$ be an $\Q$--envelope of $\{\varphi_{1},\varphi_{2}\}.$
Suppose that algebraic pair  $(A,t\omega)$ has a mirror symmetric pair
for any $t\in \R^{*}.$
Then for any $\omega'\in C_{A}\bigcap W_{\C}$
the pair $(A,\omega')$ is well-becoming.
\eth
\pr
As above, consider operators
$$
I_{t\omega}:=
\begin{pmatrix}
\varphi_2^{-1}\varphi_1&-t^{-1}\varphi_2^{-1}\\
t(\varphi_2 +\varphi_1\varphi_2^{-1}\varphi_1)&-\varphi_1\varphi_2^{-1}
\end{pmatrix}
$$
that act on the real space $\La_{\R}= (\Ga\op\Ga^{*})_{\R}.$
By assumption, any pair $(A, t\omega)$ has a mirror symmetric one.
Therefore there is a decomposition $\La=\Si_{1t}\op\Si_{2t}$ for any
$t\in \R^{*}$ such that both $\Si_{kt}$ and  are $Q$--isotropic and
$\Si_{kt\R}$
are
$I_{t\omega}$--invariant
($k=1,2$).
Since the set of all sublattices in a lattice is countable
there is a decomposition
$\La=\Si_{1}\op\Si_{2}$
such that $\Si_{k\R}$ are
$I_{t\omega}$--invariant
for infinite set of $t.$
Hence these subspaces are invariant with respect to the
operators
$$
\gamma_1=\begin{pmatrix}
0&\varphi_2^{-1}\\
0&0
\end{pmatrix},
\qquad
\gamma_2=\begin{pmatrix}
\varphi_2^{-1}\varphi_1&0\\
0&-\varphi_1\varphi_2^{-1}
\end{pmatrix},
\qquad
\gamma_3=\begin{pmatrix}
0&0\\
\varphi_2+\varphi_{1}\varphi_{2}^{-1}\varphi_{1}&0
\end{pmatrix}.
$$
It follows from Lemmas 9.4.1 and 9.4.2 that the map
$\psi=\varphi_{2}+\varphi_{1}\varphi_{2}^{-1}
\varphi_{1}$
from
$\Ga_{\R}$ to $\Ga_{\R}^{*}$
is invertible.

Denote by
$\Pi$ and $\Xi$
the projections of sublattice
$\Si_{1}$
on $\Ga$
and
$\Ga^*$
respectively.
It is clear that
$\Si_{1}\subseteq \Pi\op\Xi.$
Since
$\Si_{1}$
is invariant under the action of $\gamma_1$ and $\gamma_3$
we have
$\gamma_{1}\gamma_{3}(\Si_{1\R})\subset \Si_{1\R}.$
One can see that
$\gamma_{1}\gamma_{3}(\Si_{1\R})=(\psi\varphi^{-1}_{2}(\Pi_{\R}),0)$
This imples that
$\psi\varphi^{-1}_{2}(\Pi_{\R})\subseteq\Pi_{\R}.$
As $\psi$ is invertible, we have
that
$\psi\varphi^{-1}_{2}(\Pi_{\R})=\Pi_{\R}$
 and
$(\Pi_{\R},0)\subset\Si_{1\R}.$
By the same argument,
$(0,\Xi_{\R})\subset\Si_{1\R}.$
Hence
$\Si_{1\R}=\Pi_{\R}\op\Xi_{\R}$
and, consequently,
$\Si_{1}$
is a sublattice
of finite index
in
$\Pi\op\Xi.$
But, as the sublattice
$\Si_{1}$
is a direct summand of
$\La$
it coincides with
$\Pi\op\Xi.$
Moreover, it is easy to see that
$\dim\Pi=\dim\Xi=n,$ because
there are inclusions
$\psi(\Pi_{\R})\subseteq\Xi_{\R}$
and
$\varphi^{-1}(\Xi_{\R})\subseteq
\Pi_{\R}$
which are actually equalities.

Futher, take
$x, y\in \Pi_{\R}.$
There is
$l\in\Xi_{\R}$
such that
$x=\varphi_{2}^{-1}(l).$
Since
$\Si_{1}$
is isotropic with respect to the
form
$Q$
we obtain that
$\varphi_{2}(x, y)=l(y)=0.$
Similarly, we have
$\varphi_{1}(x, y)=\varphi_{1}\varphi^{-1}_{2}(l)(y)=0.$
Hence
$\Pi_{\R}$
is isotropic with respect to
$\varphi_{2}$ and $\varphi_{1}.$
This yields that the lattice
$\Pi$
is isotropic with respect to
any element from
$W$
that is the
$\Q$--envelope of
$\{\varphi_{1},\varphi_{2}\}.$

Now put
$\Ga_{1}:=\Pi.$
Produce a sublattice
$\Ga_2\subset \Ga$
from
$\Si_{2}$
in the
same way as
$\Ga_{1}$
from
$\Si_{1}.$
It is clear that
$\Ga=\Ga_{1}\op\Ga_{2}.$
By construction,
$\Ga_{1}$
and
$\Ga_{2}$
are isotropic with respect to any
element of
$W.$
Let $s\in W\bigcap C_{A}$ be a form that corresponding to
an ample line bundle on $A.$
Hence the operator $J_{A}$ preserves $s_{\R}$
and the symmetric form $s\langle J(\cdot), \cdot\rangle$
is positive definite.
It can be checked as in  Lemma 9.4.2 that $J_{12}, J_{21}$
are invertible.
Thus
$(A, \omega')$
is well-becoming.
$\Box$

\ssec
Let
$(A, \omega)$
and
$(B, \theta)$ be two mirror symmetric well--becoming pairs, as
in Theorem \ref{mirsym}. Then
$H_1(A, \Z) =  \Ga = \Ga_1 \op \Ga_2$ and
$H_1(B, \Z) =  \Si = \Ga_1^* \op \Ga_2.$
According to (3.7), both cohomology lattices
$H^{*}(A,\Z)$ and
$H^{*}(B, \Z)$  carry the structure of
$Cl(\La, Q)$--modules.

Recall  that by Proposition 9.3.3 there exists a unique (up to
$ \pm 1$) isomorphism of
$Cl(\La, Q)$--modules
$\beta : H^{*}(A, \Z) \stackrel{\sim}{\lto} H^{*}(B, \Z).$
It is clear that any such isomorphism can be represented as
\begin{equation}\label{transform}
v_{\xi}(\cdot) = p_{2*}(\xi \cup p_1^*(\cdot))
\end{equation}
 for some class
$\xi \in H^*(A \ts B, \Z).$

We want to show that this class
$\xi$ is in fact the Chern character of some complex line bundle on
a real subtorus of the product
$A \times B.$ As we shall see, this torus is of real dimension
$3n,$ where
$n = \dim_{\C} A = \dim_{\C}  B.$

\sssec
Let us fix bases
$\langle l_1, \dots , l_n\rangle$ of
$\Ga_1$ and
$\langle l_{n+1}, \dots , l_{2n}\rangle$ of
$\Ga_2.$ Let the dual bases of
$\Ga^*_1$ and
$\Ga_2^*$ be
$\langle x_1, \dots , x_n\rangle$  and
$\langle x_{n+1}, \dots , x_{2n}\rangle,$ respectively.

Let, as in (3.2),
$I_A$ and
$I_B$ be the following left ideals of
$Cl(\La, Q):$
$$ I_A = Cl(\La, Q)
\cdot l_1 \cdots l_{2n},
\qquad I_B = Cl(\La, Q)
\cdot l_{n+1} \cdots l_{2n} x_1 \cdots x_n.$$

Since
$H^1(A, \Z) \cong \Ga^*,$ we
identify
$H^*(A, \Z)$ with the ideal
$I_A$ as
$Cl(\La, Q)$--modules by the following rule:
$ x_{i_1} \cup \cdots \cup x_{i_k} \in H^k(A, \Z)$ goes to
$ x_{i_1} \cdots x_{i_k} \cdot l_1 \cdots l_{2n} \in I_A.$
In the same manner, since
$H^1(B,\Z) \cong \Si^* = \Ga_1 \op \Ga_2^*,$ we identify
$H^*(B, \Z)$ with the ideal
$I_B$ as
$Cl(\La, Q)$--modules.

Both
$I_A$ and
$I_B$ are isomorphic irreducible
$Cl(\La, Q)$--modules.
Therefore, there is a unique (up to
$\pm 1$) isomorphism between them.
It is given by the right multiplication
$I_A \cdot x_1 \cdots x_n = I_B.$

All the above gives us a formula for cohomology map
$\beta$
in terms of the
bases chosen above:
\begin{equation}\label{bet}
 x_{S} \cup x_R
\mapsto
(-1)^{\epsilon}
 x_R \cup l_{\overline{S}},
\end{equation}
with
\begin{equation}\label{znak}
\epsilon = |S| \cdot |R| + \sum\limits_{i \in S} (i-1),
\end{equation}
Here
$S\subset \{1,\dots , n\},\quad R\subset \{ n+1, \dots, 2n\},$
$\overline{S}$ is  the complement subset,
and
$x_S = x_{i_1} \cup \cdots \cup x_{i_{|S|}},$ where
$i_1 < i_2 < \cdots <i_{|S|}$ are the elements of
$S,$ etc.

\ssec
Now that we have an explicit formula for the map
$\beta$ in the chosen bases, our next step is to
compute the class
$\xi \in H^*(A \times B, \Z),$ that represents
$\beta$ in the form (\ref{transform}).
In fact, we will define a real
subtorus in
$A \ts B$ with a complex line bundle on it, and then show
that the direct image of the Chern character of this line bundle
is the desired
$\xi.$
First,  take this real torus to be
$T:=\Pi_{\R}/\Pi,$ where
$\Pi = \Ga_1 \op \Ga_2 \op \Ga_1^*.$
Clearly, since
$\Pi$ is a sublattice of
$\La$ one can consider
$T$ as a subtorus both of
$A \ts \wh{A}$  and
$B \ts \wh{B},$ because the first homology lattices of these two
are identified with
$\La.$

\lm{Lemma} The restrictions
$P_{A}\left|_T\right.$ and
$P_B^{-1}\left|_T\right.$ are isomorphic as complex line bundles on $T.$
\elm
\pr It follows from the well-known  fact that the complex line bundle is determined by its first Chern class (see for example \cite{Hus} or \cite{Ka}).
Using Lemma 4.2.3.1, it is easy to check that the first Chern classes of these
complex line bundles  coincide.
$\Box$

\sssec
Put
$L:=\left. P_{A}\right|_{T}.$
Lemma 4.2.3.1 gives us the formula for
$c_1(P_A).$ So we have
$$
c_1(L)=\sum\limits_{i=1}^{n} x_{i}\cup l_{i},
$$
where $x_{i}$ and $l_{i}$ are elements
of
$H^{1}(T, \Z)\cong \Pi^{*}= \Ga_{1}^{*}\op \Ga_{2}^{*}\op \Ga_{1}.$

Note that
$T$ has projections onto
$A$ and
$B$ and thus is naturally embedded into the product
$j: T \hookrightarrow A \ts B.$

\sssec
Let now
$\tau$ be a class in
$H^n(A \ts B, \Z)$ such that
if we evaluate the formula (\ref{transform}) on the class
$\tau,$ the resulting transformation sends
the monomials of the form
$x_R \cup x_n \cup \cdots \cup x_1 \in H^*(A, \Z)$ to
$x_R \in H^*(B, \Z),$
and any other monomials  to zero (here
$R \subset \{n+1, \dots, 2n\}).$
The class $\tau$ corresponds  to the
homology class
$[T]$ in
$H_{3n}(A \ts B, \Z)$
under Poincare isomorphism and suitable choice of an orientation on $T.$

Any class in
$H^*(T, \Z)$ is the restriction of a class  from
$H^*(A\ts B, \Z).$
In particular,
$c_1(L)$ is the restriction of the class
$D \in H^2(A\ts B , \Z)$ given by the following formula:
$$
D=\sum\limits_{i=1}^{n} p^*_1(x_{i}) \cup p_2^*(l_{i}).
$$

Therefore, one has
$$
j_*(ch(L))= j_*j^*(\exp(D))=\tau\cup\exp (D),
$$
where $j_{*}$ is map from $H^{*}(T, \Z)$ to
$H^{*+n}(A\ts B, \Z)$ that is obtained from
evident homology map $H_{*}(T, \Z)\lto H_{*}(A\ts B, \Z)$
under Poincare isomorphism.

\th{Proposition} The map
$\beta$ coincides with the map given by the formula (\ref{transform})
with
$\xi=j_*(ch(L^{(-1)^{(n-1)}}).$
\eth
\pr
Let
$S$ and
$R$ be index subsets of
$\{1, \dots , n\}$ and
$\{n+1, \dots , 2n\},$ respectively. Take a monomial
$x_S \cup x_R$ in the cohomologies of
$A.$
We already know the image of this element under
$\beta$ (see (\ref{bet})). Now we are going to compute its image under
the map
given  by
$v_{\xi}.$
\begin{equation}
\label{transform-by-L}
v_{\xi}(x_S \cup x_R)=
p_{2*}(\tau\cup \exp((-1)^{n-1}D) \cup p_{2}^*(x_{S}\cup x_{R})).
\end{equation}
Computing the Chern character
$\exp((-1)^{n-1}D)$ analogously  to Lemma 4.2.3.2, we see that
$$
\exp((-1)^{n-1} D)=\sum\limits_{Q\subset\{1,\cdots, n\}}(-1)^{(n-1)|Q|}
p_{1}^{*}({}_{Q}x) \cup p_{2}^{*}(l_{Q})
$$
where ${}_{Q}x$ is a product $x_{i_{|Q|}}\cup \cdots\cup x_{i_1}$
by set $Q=\{i_1,\dots, i_{|Q|}\}$ in the descending order.
Therefore, the entire expression on the right--hand side of
(\ref{transform-by-L})
rewrites as
$$
\begin{array}{rl}
v_{\xi}(x_{S} \cup x_{R})=
(-1)^{(n-1)|\overline{S}|}&\hspace{-8pt}p_{2*}(\tau\cup p_{1}^*(x_{S})\cup p_{1}^{*}(x_{R})
\cup
p_{1}^{*}( {}_{\overline{S}}x )\cup
p_{2}^{*}(l_{\overline{S}}))= \\
(-1)^{(n-1)|\overline{S}|+|S||R|}&\hspace{-8pt}p_{2*}(\tau\cup p_{1}^{*}(x_{R}\cup
x_{S}\cup{}_{\overline{S}}x))\cup p_{2}^{*}( l_{\overline{S}}))=\\
(-1)^{(n-1)|\overline{S}|+|S||R|+\sum\limits_{i\in S}(n-i)}&\hspace{-8pt}
p_{2*}(\tau\cup p_{1}^*(x_{R}\cup x_{n}\cup\cdots\cup x_{1}) \cup p_{2}^*(
l_{\overline{S}}))=\\
(-1)^{(n-1)(|\overline{S}|+|S|)+|S||R|+\sum\limits_{i\in S}(i-1)}&\hspace{-8pt}
p_{2*}(\tau\cup p_{1}^{*}(x_{R}\cup x_{n}\cup \cdots \cup x_1))
\cup l_{\overline{S}}=\\
(-1)^{(n-1)n +\epsilon}&\hspace{-8pt} x_{R}\cup l_{\overline{S}}=\\
(-1)^{\epsilon}&\hspace{-8pt} x_{R}\cup l_{\overline{S}}
\end{array}
$$
in  $H^*(B, \Z)$ and $\epsilon$ was defined by formula
(\ref{znak}). Thus, the maps $v_{\xi}$ and $\beta$ coincide.
$\Box$

\sect{APPENDIX (the proof of Theorem 8.2).}

\noindent{\bf A.1 Theorem.} {\it Let $A$ be an abelian variety.
\begin{list}{\alph{tmp})}%
{\usecounter{tmp}}
\item
 The action of $\U_{A,\Q}(\R)$ on
$C_A $
is well defined
\item
 $C^+_A$ and $C^-_A$ are single orbits under this action.
\item
 The stabilizer of a point in
$C_A$
is a maximal compact subgroup of
$\U_{A,\Q}(\R)$
\end{list}}
\pr Let us prove the assertion about $C_A^+.$ The case of $C^-_A$
is similar.
It will be more convenient to work with a different open subset of
$NS_A(\C). $
Namely, recall that every element
$z \in NS_A(\R)^0$
defines the corresponding Rosati involution on
$\End(A)  \otimes  \R:$
$$a \mapsto a' :=z^{-1} \wh{a} z.$$

Put
$NS_A(\R)^+:= \{z \in NS_A(\R) \mid \forall a \in \End(A)  \otimes \R
\quad Tr(aa')>0 \}.$
That is
$NS_A(\R)^+$ is an open subset of
$NS_A(\R)$ which consists of elements
$z$ that define a positive Rosati involution.

\medskip
\noindent{\bf A.2 Lemma.} {\it The ample cone
$C_A^a$
is a connected component of
$NS_A(\R)^+.$ Hence $C^+_A$ is a connected component of
$NS_A(\R)+iNS_A(\R)^+.$}

\smallskip
\pr We know [Mum1] that an ample class
$z \in NS_A$
belongs to $NS_A(\R)^+.$ Hence
$C^a_A \subset NS^+_A.$

Every
$z \in NS_A(\R)$
defines a symmetric bilinear form
$b_z$ on
$V_A$ as follows
$$
b_z(x,y):=z(J_Ax)(y).
$$
By the Lefschetz theorem ([Mum1]) a class
$z \in NS_A$ is ample iff
$b_z$ is positive definite.  Thus we obtain another characterization of
the ample cone
$$C_A^a = \{z \in NS_A(\R) \mid b_z  \ is \ positive \ definite
\}.$$
It follows that for
$z \in \partial C_A^a$ the form
$b_z$ is degenerate, i.e. the map
$z$ has nontrivial kernel. Thus
$
\partial C_A^a \subset \partial NS_A(\R)^+.
$
Therefore,
$C_A^a$ is a connected component of
$NS_A(\R)^+.$ $\Box$

In the next few lemmas we study the
$\U_{A,\Q}(\R)$-action on the set $NS_A(\R)+iNS_A(\R)^+.$

\medskip
\noindent{\bf A.3 Lemma.} {\it Let
$g=  \left(
\begin{smallmatrix}a&b\\c&d
\end{smallmatrix}
\right) \in \U_{A,\Q}(\R),$ $ \omega \in NS_A(\R) + iNS_A(\R)^+.$
Then the endomorphism
$(a+b \omega) \in \End(V_A) $ is invertible, i.e. the element
$g\omega \in NS_A(\C) $ is well defined.}

\smallskip
\pr Let
$\omega = \eta  + iz,\quad  \eta \in NS_A(\R),\quad  z \in NS_A(\R)^+.$
The lemma
follows from the following assertion: for every nonzero
$x \in \End(V_A\otimes _{\R} \C)$
$$
Im(Tr(z^{-1}  \wh{\bar x} (\wh c + \wh {\bar \omega} \wh d)(a +
b \omega ) x ) ) < 0.
$$
The assertion is proved by a straightforward calculation which we omit .

\medskip
\noindent{\bf A.4 Remark.} Let
$\omega = \eta  + iz \in NS_A(\C).$ Then
$\left(
\begin{smallmatrix}1&0\\-\eta&1
\end{smallmatrix}
\right) \in \U_{A,\Q}(\R)$ and
$\left(
\begin{smallmatrix}1&0\\-\eta&1
\end{smallmatrix}
\right) \omega = iz.$ Therefore all
$\U_{A,\Q}(\R)$-orbits in
$NS_A(\C)$ are invariant under translation along
$NS_A(\R).$

\medskip
\noindent{\bf A.5 Lemma.}  {\it Let
$\omega \in NS_A(\R) + iNS_A(\R)^+.$ Let
$K_{\omega}$ be its stabilizer in
$\U_{A,\Q}(\R). $ Then
$K_{\omega}$ is a maximal compact subgroup of
$\U_{A,\Q}(\R).$}

\smallskip
\pr By the previous remark we may assume that
$\omega = iz.$ Denote by
$'$ the Rosati involution defined by
$z.$ We have
$$
K_{\omega} = \left  \{
\begin{pmatrix}a&b\\-zbz&zaz^{-1}
\end{pmatrix}
 \in \U_{A,\Q}(\R)
\right|\left. \;  z^{-1} \wh a z a + \wh b z b z = 1, \   z^{-1} \wh a z b = \wh b z a
z^{-1}  \right\}
$$

Consider a map
$
\te :  \U_{A,\Q}(\R) \to M(2, \End(A)\otimes\R)^*,$
given as
$
\te
\begin{pmatrix}a&b\\c&d
\end{pmatrix}
=
\begin{pmatrix}a&bz\\z^{-1}c&z^{-1}dz
\end{pmatrix}.
$

Then
$$
Im \te = \left  \{
\begin{pmatrix}\alpha&\beta\\\gamma&\delta
\end{pmatrix} \right| \left.
\begin{pmatrix}\alpha&\beta\\\gamma&\delta
\end{pmatrix}^{-1}
= \begin{pmatrix}\delta'&-\beta'\\-\gamma'&\alpha'
\end{pmatrix} \right \}
\quad
\mbox{and}
\quad
\te (K_{\omega}) = \left \{ \begin{pmatrix}\alpha&\beta\\
-\beta&\alpha
\end{pmatrix}
\right| \left. \begin{matrix}\alpha' \alpha  + \beta' \beta =Id,&\\
\alpha'\beta = \beta' \alpha&
\end{matrix}
\right\}
$$
Extend the Rosati involution
$'$
to
$M(2, \End(A)) \otimes \R$
by the formula
$$
\left(\begin{smallmatrix}\alpha&\beta\\\gamma&\delta
\end{smallmatrix}\right)'=
\left(\begin{smallmatrix}\alpha'&\gamma'\\\beta'&\delta'
\end{smallmatrix}
\right)
$$
and consider the symmetric bilinear form on
$M(2, \End(A)) \otimes \R$
$$
\Psi(X,Y) =Tr(X'Y), \quad
\mbox{i.e.}
\quad
\Psi(
\begin{pmatrix}\alpha&\beta\\\gamma&\delta
\end{pmatrix},
\begin{pmatrix}a&b\\c&d
\end{pmatrix})
=Tr(
\begin{pmatrix}\alpha'&\gamma'\\\beta'&\delta'
\end{pmatrix}
\begin{pmatrix}a&b\\c&d
\end{pmatrix}) =
Tr(\alpha' a +\gamma' c + \beta' b + \delta' d).
$$
This form is positive definite.

Consider the action of the group
$ M(2, \End(A) \otimes \R )^*$
on
$ M(2, \End(A) \otimes \R )$ by left multiplication. Let
$\U(\Psi) \subset M(2, \End(A) \otimes \R )^*$ be the subgroup which
preserves the form  $\Psi.$ We claim  that
$$
\te(\U_{A,\Q}(\R))  \cap \U(\Psi) = \te(K_{\omega}).
$$
Indeed,
$$
\U(\Psi)= \{u \in M(2, \End(A) \otimes \R )^* \mid u'u = Id \}.
$$
Thus
$\left(\begin{smallmatrix}a&b\\c&d
\end{smallmatrix}
\right)  \in  \te(\U_{A,\Q}(\R) )\cap \U(\Psi) $
iff
$\left(\begin{smallmatrix}a'&b'\\c'&d'
\end{smallmatrix}\right) =
\left(\begin{smallmatrix}d'&-b'\\-c'&a'
\end{smallmatrix}\right).$
This means that
$a=d,$ $b=-c,$ $a'a+b'b=1,$ $a'b=b'a,$
which in turn means that
$\left(\begin{smallmatrix}a&b\\c&d
\end{smallmatrix}
\right)  \in \te(K_{\omega}).$
Thus
$K_{\omega}$ is compact.

Notice that the adjoint operator for the left
multiplication by
$m \in M(2, \End(A) \otimes \R)$ with respect to
$\Psi$ is the left multiplication by
$m'.$ The group
$\te(\U_{A,\Q}(\R))$ is self adjoint, i.e.
$t \in  \te(\U_{A,\Q}(\R)) \Rightarrow  t' \in \te(\U_{A,\Q}(\R)).$
Put
$$
P_\omega = \{g \in  \U_{A,\Q}(\R)  \mid \te(g) \; is \;  positive \;
self\; adjoint
\}.
$$
The lemma now follows from the following claim.

\medskip
\noindent{\bf A.6 Claim.} {\it
\begin{list}{\alph{tmp})}%
{\usecounter{tmp}}
\item
 The multiplication map
$
K_{\omega} \times P_\omega \to \U_{A,\Q}(\R),$
such that
$(b,p) \mapsto bp
,$
is a diffeomorphism
\item
$K_{\omega}$
is a maximal compact subgroup of
$\U_{A,\Q}(\R).$
\end{list}
}

\smallskip
Indeed, since
$\te(\U_{A,\Q}(\R))$
is self-adjoint, a) follows from theorem 1 in sect.2, ch.5 in [OV].
The assertion b) follows from problem 2 in sect.2, ch.5 in [OV]. $\Box$

\medskip
\noindent{\bf A.7 Lemma.} {\it The set
$NS_A(\R)+i NS_A (\R)^+ $ is
$\U_{A,\Q}(\R)$-invariant.}

\smallskip
\pr Let
$\omega \in NS_A(\R) + iNS_A(\R)^+, g \in \U_{A,\Q}(\R).$
We need to show that
$g\omega \in NS_A(\R) + iNS_A(\R)^+.$
By Remark A.4 above we may asssume that
$\omega =iz, \quad  z \in NS_A(\R)^+,$ and
$g\omega = ip,\quad  p \in NS_A(\R).$ Then we need to show that

(i)
$p$ is invertible,

(ii) the Rosati involution defined by
$p $ is positive definite (i.e.
$p \in NS_A(\R)^+.$)

Denote by  $'$ the Rosati involution defined by
$z.$ Put
$p=zk$ for
$k \in \End(A) \otimes \R.$ We have
$k'=k.$ Let $\lambda \in \U_{A,\Q}(\R)$ be such that $\lambda (iz)=ip.$
Then
$\lambda$ has the form
$$
\lambda = \begin{pmatrix}a&b\\-zkbz&zkaz^{-1}
\end{pmatrix}.
$$
Under the map
$
\te : \U_{A,\Q}(\R) \to M(2, \End(A) \otimes \R),
$
defined in the proof of Lemma A.5,
$\lambda$ goes to
$$
\te(\lambda) =
\begin{pmatrix}\alpha&\beta\\-k\beta&k\alpha
\end{pmatrix},
$$
where
$$
1)\; \alpha\alpha' k + \beta\beta' k = Id
\quad
\mbox{and}
\quad
2)\; -\alpha \beta' + \beta \alpha' = 0
$$

1) implies that
$k$ is invertible , hence
$p$ is such. This proves (i). Let us prove (ii). Note that
$$
\begin{pmatrix}k&0\\0&\wh k^{-1}
\end{pmatrix} \in \U_{A,\Q}(\R)
\quad
\mbox{and}
\quad
\begin{pmatrix}k&0\\0&\wh k^{-1}
\end{pmatrix}ip = izk^{-1}.
$$

Replacing
$izk$ by
$izk^{-1}$
and applying the previous argument we find
$s,t \in \End(A)\otimes \R$
such that
$
ss'k^{-1} + tt'k^{-1} =1,\quad -st' +ts' =0,
$
that is,
$$
3)\; k=ss' +tt',
\quad
\mbox{and}
\quad
4)\; -st'+ts'=0
$$

The Rosati involution defined by
$p$ is
$$
a \mapsto p^{-1} \wh  a p = k^{-1} a' k, a \in \End(A) \otimes \R.
$$
We need to show that
$Tr(k^{-1}a'ka)>0$ for all
$0 \ne  a \in \End(A) \otimes \R .$ Or, equivalently, that the following quadratic form on
$M(2, \End(A) \otimes \R)$ is positive definite
$$
Tr(
\begin{pmatrix}k^{-1}&0\\0&k^{-1}
\end{pmatrix}
X'
\begin{pmatrix}k&0\\
0&k
\end{pmatrix}X),\quad X \in M(2, \End(A) \otimes \R).
$$
Put
$Y=\left(
\begin{smallmatrix}\alpha&\beta\\-\beta&\alpha
\end{smallmatrix}\right)$ and
$ Z=
\left(
\begin{smallmatrix}s&t\\-t&s
\end{smallmatrix}
\right).$
Then
$$
1),2) \Leftrightarrow YY' =\begin{pmatrix}k^{-1}&0\\0&k^{-1}
\end{pmatrix},
\quad
\mbox{and}
\quad
3),4) \Leftrightarrow ZZ' =\begin{pmatrix}k&0\\0&k
\end{pmatrix}.
$$
Hence
$$
Tr(\begin{pmatrix}k^{-1}&0\\0&k^{-1}
\end{pmatrix}X'\begin{pmatrix}k&0\\0&k
\end{pmatrix}X)=
Tr(YY'X'ZZ'X)=Tr(Y'X'ZZ'XY)=Tr((Z'XY)'(Z'XY))>0
$$
This proves (ii) and the Lemma.

\medskip
\noindent{\bf A.8 Lemma.} {\it The set
$NS_A(\R) + iNS_A(\R)^+$
consists of finitely many
$\U_{A,\Q}(\R)$-orbits.}

\medskip
\noindent{\bf A.9 Corollary.} {\it The
$\U_{A,\Q}(\R)$-orbits in
$NS_A(\R) + iNS_A(\R)^+$ coincide with the connected components of this set.}

\smallskip
\noindent Proof of corollary. By the Corollary 5.3.5 the group
$\U_{A,\Q}(\R)$ is connected. Hence each
$\U_{A,\Q}(\R)$-orbit is contained in a connected component of
$NS_A(\R) + iNS_A(\R)^+.$
It follows from Lemma A.5 that all
$\U_{A,\Q}(\R)$-orbits in
$NS_A(\R) + iNS_A(\R)^+$ are isomorphic (all maximal compact
subgroups in a reductive Lie group are conjugate). Then by Lemma A.8 each
$\U_{A,\Q}(\R)$-orbit is an open subset in
$NS_A(\R) + iNS_A(\R)^+,$ hence must coincide with a
connected component of this set. This proves the corollary.

\medskip
\noindent{\bf A.10 Remark.} Theorem A.1 now follows from A.2, A.3, A.5 and
A.9. So it remains to prove Lemma A.7.

\medskip
\noindent {\it Proof of Lemma A.7.} By Remark A.4 above it suffices to show
that the set
$iNS_A(\R)^+$ is contained in a finite number  of
$\U_{A,\Q}(\R)$-orbits.

Fix
$z \in NS_A(\R)^+$ and let ${}\prime$ be the corresponding Rosati involution.
Let
$l \in Aut_A(\R) :=(\End(A) \otimes \R)^*.$ Then
$iz$ and
$izl'l$ are
$\U_{A,\Q}(\R)$--conjugate. Indeed,
$$
\begin{pmatrix}l^{-1}&0\\0&\wh{l}
\end{pmatrix} iz = i \wh{l} z l = i zl'l.
$$
Consider the space
$\End(A) \otimes \R$ with the positive definite symmetric form
$
\phi(a,b) :=Tr(ab')
$
and the adjoint action of the group
$Aut_A(\R)$:
$ad_s(a):=sas^{-1}.
$

The adjoint with respect to
$\phi$ of the operator
$ad_s$ is
$ad_{s'}.$ Hence the following are equivalent

(i)
$ad_s$ is self adjoint;

(ii)
$s=s'c$ for some
$c$ in the center of
$Aut_A(\R).$

Let now
$izk \in  iNS_A(\R)^+$ for some
$k \in Aut_A(\R).$

\medskip
\noindent{\bf A.11 Claim.} {\it There exist
$l \in Aut_A(\R)$ and
$c \in Z(Aut_A(\R))$ such that
$c'=c, c^2=1$ and
$k=l'lc.$}

\smallskip

Indeed, we have
$\widehat{zk} = zk, $ i.e.
$k'=k.$ Thus
$ad_k$ is self adjoint with respect to
$\phi.$ The positivity of the Rosati involution defined by
$zk$
$$
a \mapsto  k^{-1}a'k
$$
is equivalent to the positivity of the self adjoint operator
$ad_k.$ Hence the operator
$ad_k$ has a square root (which is also positive self adjoint).
  But this square root is also in the image of the homomorphism
$ad$ (see ch. 5, sect 2, thm. 1 in [OV]). Thus there exists
$l \in Aut_A(\R)$ such that
$ad_l =ad_{l'}$ and
$(ad_l)^2=ad_k.$ We get
$l'lc=k$ for some
$c\in Z(Aut_A(\R)).$ Moreover,
$$
l'lc=k=k'=(l'lc)'=c'l'l.
$$
Thus
$c'=c.$

It follows from the discussion of
$\End(A)$ in section 1.8 that the $\prime$-invariant part of the center
$Z(\End(A) \otimes \R)$ is isomorphic to
$\R \times \dots \times \R.$ therefore, multiplying
$c$ by a positive real number
$r$ (and
$l$ by
$r^{-1/2}$) we may assume that
$c^2=1.$ This proves the claim.

\medskip
Note that there are finitely many
$c \in Z(Aut_A(\R))$ such that
$c'=c$ and
$c^2=1,$ say,
$c_1, \dots, c_m.$ Then by the argument in the beginning of the proof of
the lemma
$iNS_A(\R)^+$ is contained in the union of
$\U_{A,\Q}(\R)$-orbits of
$izc_1, izc_2, \dots, izc_m.$ This proves Lemma A.8 and completes the
proof of Theorem A.1.

\end{document}